\documentclass{article}

\usepackage{amsmath}
\usepackage{amssymb}

\title{On the representability of actions\\in the topos context}
\author{Francis Borceux}
\date{{\em To George Janelidze, on his sixtieth birthday}}

\immediate\write0{After the last statement, type the command \string\par.}
\newtheorem{theorem}{Theorem}%[section]
\newtheorem{definition}[theorem]{Definition}
\newtheorem{proposition}[theorem]{Proposition}

\newtheorem{corollary}[theorem]{Corollary}

\newtheorem{example}[theorem]{Example}

\newtheorem{counterexample}[theorem]{Counterexample}

\newtheorem{exercise}[theorem]{\hspace*{-4pt}}

\def\ex#1\par{\begin{exercise}\em#1\end{exercise}}

\outer\def\proclaim#1. #2\par{\medbreak
\noindent{\bf#1\enspace}{\sl#2}\par
\ifdim\lastskip<\medskipamount
\removelastskip\penalty55\medskip\fi}

\newcommand{\proof}{\medbreak\noindent%
{\it Proof}\hspace{1em plus .5em}}

\newcommand{\Halmos}{\mbox{$\square$}}

\newcommand{\qed}{\mbox{}\hfill\mbox{}%
\Halmos\medbreak}

\newcommand{\qedmanual}[1]%
{\makebox[0pt]{\hspace{#1mm}\Halmos}%
\typeout{*******}%
\typeout{WARNING : manually adjusted%
"end of proof" symbol}
\typeout{*******}}

\newcommand{\BbbN}{{\mathbb N}}

\newcommand{\calB}{{\mathcal B}}
\newcommand{\calC}{{\mathcal C}}

\newcommand{\calE}{{\mathcal E}}

\newcommand{\calL}{{\mathcal L}}

\newcommand{\calN}{{\mathcal N}}

\newcommand{\calP}{{\mathcal P}}

\newcommand{\calS}{{\mathcal S}}

\newif{\ifLabels}

\newcommand{\BreakAlign}{\bgroup\allowdisplaybreaks}

\newcommand{\eBreakAlign}{\egroup%
\typeout{*****}%
\typeout{WARNING: Breaks in alignment}%
\typeout{*****}}

\catcode`\@=11

\newcommand{\@makefigcaption}[2]{%
\mbox{}\hfill#1#2\hfill\mbox{}}

\renewcommand{\figure}{%
\let\@makecaption\@makefigcaption\@float{figure}}

\newcounter{FigNumber}[section]

\newcommand{\adjust}{0}
\newcommand{\topadjust}{0}

\newenvironment{Tikzpicture}[1]{%
\addtocounter{FigNumber}{1}
\immediate\write\@auxout{\string\newlabel{#1}%
{{\thesection.\theFigNumber}{thepage}}}%
\ifLabels\global\def\FigLabel{~~\protect\framebox{#1}}%
\else\global\def\FigLabel{}\fi%
\begin{figure}[tp]%
\vspace*{-\topadjust mm}%
\begin{center}%
\begin{tikzpicture}}%
{\end{tikzpicture}%
\vspace*{-\adjust mm}%
\end{center}%
\caption{\FigLabel}%
\end{figure}%
}

\newcommand{\EPS}[2]%
{\begin{Tikzpicture}{#1}%
\draw (0,0) node {\epsfig{file=#2}};%
\end{Tikzpicture}}%

\catcode`\@=12

{%
$$%
\arraycolsep=1.5pt%
\left\{%
\begin{array}{rl}%
}%
{%
\end{array}%
\right.%
\arraycolsep=5pt%
$$%
\hspace*{-4pt}%
}

\newcommand{\id}{\mathsf{id}}
\newcommand{\Ker}{\mathsf{Ker}\,}
\newcommand{\Coker}{\mathsf{Coker}\,}
\newcommand{\Set}{\mathsf{Set}}
\newcommand{\Or}{\mathsf{or}}

\newcommand{\One}{{\mathbf1}}
\let\card=\sharp
\newcommand{\EStar}{\calE_{\star}}
\newcommand{\EStarOp}{\calE_{\star}^{\rm op}}
\newcommand{\Act}{\mathsf{Act}}
\newcommand{\Aut}{\mathsf{Aut}}
\catcode`\@=11

\newcount\printerresolution
\global\printerresolution=300

\newcount\headsize
\global\headsize=450

\font\dotfont = lcircle10 at 3pt

\newcount\epihead
\global\epihead=200

\newcount\defaultscale
\def\setdefaultscale#1{\global\defaultscale=#1}
\setdefaultscale{100}

\newcount\actualtextarrowspace
\newcount\actualtextarrowlength

\newcommand{\computetextparameters}%
{\global\actualtextarrowspace=\textarrowlength%
\global\advance\actualtextarrowspace by 3%
\global\actualtextarrowlength=\textarrowlength%
\global\multiply\actualtextarrowlength by 100}

\newcount\textarrowlength
\def\settextarrowlength#1{\global\textarrowlength=#1%
\computetextparameters}
\settextarrowlength{20}

\newcount\actualdisplayarrowspace
\newcount\actualdisplayarrowlength

\newcommand{\computedisplayparameters}%
{\global\actualdisplayarrowspace=\displayarrowlength%
\global\advance\actualdisplayarrowspace by 3%
\global\actualdisplayarrowlength=\displayarrowlength%
\global\multiply\actualdisplayarrowlength by 100}

\newcount\displayarrowlength
\def\setdisplayarrowlength#1{\global\displayarrowlength=#1%
\computedisplayparameters}
\setdisplayarrowlength{40}

\def\@ifnexttok#1#2#3{\let\@tempe #1\def\@tempa{#2}\def\@tempb{#3}%
\futurelet\@tempc\@ifntok}
\def\@ifntok{\ifx \@tempc \@tempe\let\@tempd\@tempa\else\let\@tempd\@tempb\fi%
\@tempd}

\def\DOT{{\dotfont q}}

\newcount\basicstep
\global\divide\basicstep by \printerresolution

\newcount\numberofsteps
\numberofsteps=\headsize
\global\divide\numberofsteps by \basicstep

\newcommand{\makehead}[3]{%
\begin{picture}(0,0)%
\multiput(0,0)(#1,#2){#3}{\DOT}%
\multiput(0,0)(-#2,#1){#3}{\DOT}%
\end{picture}}

\newcount\xstep
\newcount\ystep

\newsavebox{\northhead}
\savebox{\northhead}{%
\xstep=-\basicstep%
\multiply\xstep by 7071%
\divide\xstep by 10000%
\ystep=\xstep%
\makehead{\xstep}{\ystep}{\numberofsteps}}
\newcommand{\nhead}{\usebox{\northhead}}

\newsavebox{\easthead}
\savebox{\easthead}{%
\xstep=-\basicstep%
\multiply\xstep by 7071%
\divide\xstep by 10000%
\ystep=-\xstep%
\makehead{\xstep}{\ystep}{\numberofsteps}}
\newcommand{\ehead}{\usebox{\easthead}}

\newsavebox{\southhead}
\savebox{\southhead}{%
\xstep=\basicstep%
\multiply\xstep by 7071%
\divide\xstep by 10000%
\ystep=\xstep%
\makehead{\xstep}{\ystep}{\numberofsteps}}
\newcommand{\shead}{\usebox{\southhead}}

\newsavebox{\westhead}
\savebox{\westhead}{%
\xstep=\basicstep%
\multiply\xstep by 7071%
\divide\xstep by 10000%
\ystep=-\xstep%
\makehead{\xstep}{\ystep}{\numberofsteps}}
\newcommand{\whead}{\usebox{\westhead}}

\newsavebox{\northwesthead}
\savebox{\northwesthead}{%
\makehead{0}{-\basicstep}{\numberofsteps}}
\newcommand{\nwhead}{\usebox{\northwesthead}}

\newsavebox{\northeasthead}
\savebox{\northeasthead}{%
\makehead{-\basicstep}{0}{\numberofsteps}}
\newcommand{\nehead}{\usebox{\northeasthead}}

\newsavebox{\southwesthead}
\savebox{\southwesthead}{%
\makehead{\basicstep}{0}{\numberofsteps}}
\newcommand{\swhead}{\usebox{\southwesthead}}

\newsavebox{\southeasthead}
\savebox{\southeasthead}{%
\makehead{0}{\basicstep}{\numberofsteps}}
\newcommand{\sehead}{\usebox{\southeasthead}}

\newsavebox{\eastnortheasthead}
\savebox{\eastnortheasthead}{%
\xstep=-\basicstep%
\multiply\xstep by 9486%
\divide\xstep by 10000%
\ystep=\xstep%
\divide\ystep by -3%
\makehead{\xstep}{\ystep}{\numberofsteps}}
\newcommand{\enehead}{\usebox{\eastnortheasthead}}

\newsavebox{\northnortheasthead}
\savebox{\northnortheasthead}{%
\xstep=-\basicstep%
\multiply\xstep by 9486%
\divide\xstep by 10000%
\ystep=\xstep%
\divide\ystep by 3%
\makehead{\xstep}{\ystep}{\numberofsteps}}
\newcommand{\nnehead}{\usebox{\northnortheasthead}}

\newsavebox{\southsouthwesthead}
\savebox{\southsouthwesthead}{%
\xstep=\basicstep%
\multiply\xstep by 9486%
\divide\xstep by 10000%
\ystep=\xstep%
\divide\ystep by 3%
\makehead{\xstep}{\ystep}{\numberofsteps}}
\newcommand{\sswhead}{\usebox{\southsouthwesthead}}

\newsavebox{\westsouthwesthead}
\savebox{\westsouthwesthead}{%
\xstep=\basicstep%
\multiply\xstep by 9486%
\divide\xstep by 10000%
\ystep=\xstep%
\divide\ystep by -3%
\makehead{\xstep}{\ystep}{\numberofsteps}}
\newcommand{\wswhead}{\usebox{\westsouthwesthead}}

\newsavebox{\westnorthwesthead}
\savebox{\westnorthwesthead}{%
\xstep=\basicstep%
\multiply\xstep by 3162%
\divide\xstep by 10000%
\ystep=\xstep%
\multiply\ystep by -3%
\makehead{\xstep}{\ystep}{\numberofsteps}}
\newcommand{\wnwhead}{\usebox{\westnorthwesthead}}

\newsavebox{\eastsoutheasthead}
\savebox{\eastsoutheasthead}{%
\xstep=-\basicstep%
\multiply\xstep by 3162%
\divide\xstep by 10000%
\ystep=\xstep%
\multiply\ystep by -3%
\makehead{\xstep}{\ystep}{\numberofsteps}}
\newcommand{\esehead}{\usebox{\eastsoutheasthead}}

\newsavebox{\northnorthwesthead}
\savebox{\northnorthwesthead}{%
\xstep=-\basicstep%
\multiply\xstep by 3162%
\divide\xstep by 10000%
\ystep=\xstep%
\multiply\ystep by 3%
\makehead{\xstep}{\ystep}{\numberofsteps}}
\newcommand{\nnwhead}{\usebox{\northnorthwesthead}}

\newsavebox{\southsoutheasthead}
\savebox{\southsoutheasthead}{%
\xstep=\basicstep%
\multiply\xstep by 3162%
\divide\xstep by 10000%
\ystep=\xstep%
\multiply\ystep by 3%
\makehead{\xstep}{\ystep}{\numberofsteps}}
\newcommand{\ssehead}{\usebox{\southsoutheasthead}}

\newsavebox{\easteastnortheasthead}
\savebox{\easteastnortheasthead}{%
\xstep=-\basicstep%
\multiply\xstep by 8944%
\divide\xstep by 10000%
\ystep=\xstep%
\divide\ystep by -2%
\makehead{\xstep}{\ystep}{\numberofsteps}}
\newcommand{\eenehead}{\usebox{\easteastnortheasthead}}

\newsavebox{\northnorthnortheasthead}
\savebox{\northnorthnortheasthead}{%
\xstep=-\basicstep%
\multiply\xstep by 8944%
\divide\xstep by 10000%
\ystep=\xstep%
\divide\ystep by 2%
\makehead{\xstep}{\ystep}{\numberofsteps}}
\newcommand{\nnnehead}{\usebox{\northnorthnortheasthead}}

\newsavebox{\southsouthsouthwesthead}
\savebox{\southsouthsouthwesthead}{%
\xstep=\basicstep%
\multiply\xstep by 8944%
\divide\xstep by 10000%
\ystep=\xstep%
\divide\ystep by 2%
\makehead{\xstep}{\ystep}{\numberofsteps}}
\newcommand{\ssswhead}{\usebox{\southsouthsouthwesthead}}

\newsavebox{\westwestsouthwesthead}
\savebox{\westwestsouthwesthead}{%
\xstep=\basicstep%
\multiply\xstep by 8944%
\divide\xstep by 10000%
\ystep=\xstep%
\divide\ystep by -2%
\makehead{\xstep}{\ystep}{\numberofsteps}}
\newcommand{\wwswhead}{\usebox{\westwestsouthwesthead}}

\newsavebox{\westwestnorthwesthead}
\savebox{\westwestnorthwesthead}{%
\xstep=\basicstep%
\multiply\xstep by 4472%
\divide\xstep by 10000%
\ystep=\xstep%
\multiply\ystep by -2%
\makehead{\xstep}{\ystep}{\numberofsteps}}
\newcommand{\wwnwhead}{\usebox{\westwestnorthwesthead}}

\newsavebox{\easteastsoutheasthead}
\savebox{\easteastsoutheasthead}{%
\xstep=-\basicstep%
\multiply\xstep by 4472%
\divide\xstep by 10000%
\ystep=\xstep%
\multiply\ystep by -2%
\makehead{\xstep}{\ystep}{\numberofsteps}}
\newcommand{\eesehead}{\usebox{\easteastsoutheasthead}}

\newsavebox{\northnorthnorthwesthead}
\savebox{\northnorthnorthwesthead}{%
\xstep=-\basicstep%
\multiply\xstep by 4472%
\divide\xstep by 10000%
\ystep=\xstep%
\multiply\ystep by 2%
\makehead{\xstep}{\ystep}{\numberofsteps}}
\newcommand{\nnnwhead}{\usebox{\northnorthnorthwesthead}}

\newsavebox{\southsouthsoutheasthead}
\savebox{\southsouthsoutheasthead}{%
\xstep=\basicstep%
\multiply\xstep by 4472%
\divide\xstep by 10000%
\ystep=\xstep%
\multiply\ystep by 2%
\makehead{\xstep}{\ystep}{\numberofsteps}}
\newcommand{\sssehead}{\usebox{\southsouthsoutheasthead}}

\newsavebox{\northeasteastnortheasthead}
\savebox{\northeasteastnortheasthead}{%
\xstep=-\basicstep%
\multiply\xstep by 9806%
\divide\xstep by 10000%
\ystep=\xstep%
\divide\ystep by -5%
\makehead{\xstep}{\ystep}{\numberofsteps}}
\newcommand{\neenehead}{\usebox{\northeasteastnortheasthead}}

\newsavebox{\northeastnorthnortheasthead}
\savebox{\northeastnorthnortheasthead}{%
\xstep=-\basicstep%
\multiply\xstep by 9806%
\divide\xstep by 10000%
\ystep=\xstep%
\divide\ystep by 5%
\makehead{\xstep}{\ystep}{\numberofsteps}}
\newcommand{\nennehead}{\usebox{\northeastnorthnortheasthead}}

\newsavebox{\southwestsouthsouthwesthead}
\savebox{\southwestsouthsouthwesthead}{%
\xstep=\basicstep%
\multiply\xstep by 9806%
\divide\xstep by 10000%
\ystep=\xstep%
\divide\ystep by 5%
\makehead{\xstep}{\ystep}{\numberofsteps}}
\newcommand{\swsswhead}{\usebox{\southwestsouthsouthwesthead}}

\newsavebox{\southwestwestsouthwesthead}
\savebox{\southwestwestsouthwesthead}{%
\xstep=\basicstep%
\multiply\xstep by 9806%
\divide\xstep by 10000%
\ystep=\xstep%
\divide\ystep by -5%
\makehead{\xstep}{\ystep}{\numberofsteps}}
\newcommand{\swwswhead}{\usebox{\southwestwestsouthwesthead}}

\newsavebox{\northwestwestnorthwesthead}
\savebox{\northwestwestnorthwesthead}{%
\xstep=\basicstep%
\multiply\xstep by 1961%
\divide\xstep by 10000%
\ystep=\xstep%
\multiply\ystep by -5%
\makehead{\xstep}{\ystep}{\numberofsteps}}
\newcommand{\nwwnwhead}{\usebox{\northwestwestnorthwesthead}}

\newsavebox{\southeasteastsoutheasthead}
\savebox{\southeasteastsoutheasthead}{%
\xstep=-\basicstep%
\multiply\xstep by 1961%
\divide\xstep by 10000%
\ystep=\xstep%
\multiply\ystep by -5%
\makehead{\xstep}{\ystep}{\numberofsteps}}
\newcommand{\seesehead}{\usebox{\southeasteastsoutheasthead}}

\newsavebox{\northwestnorthnorthwesthead}
\savebox{\northwestnorthnorthwesthead}{%
\xstep=-\basicstep%
\multiply\xstep by 1961%
\divide\xstep by 10000%
\ystep=\xstep%
\multiply\ystep by 5%
\makehead{\xstep}{\ystep}{\numberofsteps}}
\newcommand{\nwnnwhead}{\usebox{\northwestnorthnorthwesthead}}

\newsavebox{\southeastsouthsoutheasthead}
\savebox{\southeastsouthsoutheasthead}{%
\xstep=\basicstep%
\multiply\xstep by 1961%
\divide\xstep by 10000%
\ystep=\xstep%
\multiply\ystep by 5%
\makehead{\xstep}{\ystep}{\numberofsteps}}
\newcommand{\sessehead}{\usebox{\southeastsouthsoutheasthead}}

\newsavebox{\isomorphismmark}
\newcommand{\isomark}[1]{\savebox{\isomorphismmark}{#1}}
\isomark{$\cong$}
\newcommand{\isosign}{\usebox{\isomorphismmark}}

\newif\ifuserdist
\newsavebox{\distributormark}
\newcommand{\distmark}[1]{\ifx#1\distcircle\userdistfalse\else%
\userdisttrue\savebox{\distributormark}{#1}\fi}
\newsavebox{\distributorcircle}
\savebox{\distributorcircle}{\begin{picture}(0,0)%
\put(0,0){\circle{4}}\end{picture}}
\newcommand{\distsign}{\ifuserdist\usebox{\distributormark}\else%
\usebox{\distributorcircle}\fi}
\userdistfalse

\newcount\monolength
\newcount\epilength
\newcount\bimolength
\newcount\secondmonolength
\newcount\secondepilength

\newcount\monotail%
\global\monotail=\headsize%
\global\multiply\monotail by 7070%
\global\divide\monotail by 10000%

\newcount\truemonotail%
\newcount\trueepihead%

\def\truetail{\truemonotail=\monotail%
\multiply\truemonotail by 100%
\divide\truemonotail by \SCALE}

\def\truehead{\trueepihead=\epihead%
\multiply\trueepihead by 100%
\divide\trueepihead by \SCALE}

\newcount\Monotail%
\global\Monotail=\headsize%
\global\divide\Monotail by 2%

\newcount\Epihead%
\global\Epihead=\epihead%
\global\multiply\Epihead by 7070%
\global\divide\Epihead by 10000%

\newcount\Truemonotail%
\newcount\Trueepihead%

\def\Truetail{\Truemonotail=\Monotail%
\multiply\Truemonotail by 100%
\divide\Truemonotail by \SCALE}%

\def\Truehead{\Trueepihead=\Epihead%
\multiply\Trueepihead by 100%
\divide\Trueepihead by \SCALE}

\newcount\MonoTail%
\global\MonoTail=\headsize%
\global\multiply\MonoTail by 6324%
\global\divide\MonoTail by 10000%

\newcount\EpiHead%
\global\EpiHead=\epihead%
\global\multiply\EpiHead by 8944%
\global\divide\EpiHead by 10000%

\newcount\TrueMonoTail%
\newcount\TrueEpiHead%

\def\TrueTail{\TrueMonoTail=\MonoTail%
\multiply\TrueMonoTail by 100%
\divide\TrueMonoTail by \SCALE}%

\def\TrueHead{\TrueEpiHead=\EpiHead%
\multiply\TrueEpiHead by 100%
\divide\TrueEpiHead by \SCALE}

\newcount\monotaiL%
\global\monotaiL=\headsize%
\global\multiply\monotaiL by 3162%
\global\divide\monotaiL by 10000%

\newcount\epiheaD%
\global\epiheaD=\epihead%
\global\multiply\epiheaD by 4472%
\global\divide\epiheaD by 10000%

\newcount\truemonotaiL%
\newcount\trueepiheaD%

\def\truetaiL{\truemonotaiL=\monotaiL%
\multiply\truemonotaiL by 100%
\divide\truemonotaiL by \SCALE}%

\def\trueheaD{\trueepiheaD=\epiheaD%
\multiply\trueepiheaD by 100%
\divide\trueepiheaD by \SCALE}

\newcount\SCALE%

\newcount\NUMBER%

\newcount\LINE%

\newcount\COLUMN%

\newcount\WIDTH%

\newcount\SOURCE%

\newcount\ARROW%

\newcount\TARGET%

\newcount\ARROWLENGTH%

\newcount\NUMBEROFDOTS%

\newcounter{x}%

\newcounter{y}%

\newcounter{z}%

\newcounter{horizontal}%

\newcounter{vertical}%

\newskip\itemlength%

\newskip\firstitem%

\newskip\seconditem%

\newcommand{\printarrow}{}%

\newcount\X%

\newcount\Y%

\newcount\Z%

\newcommand{\truex}[1]{%
\NUMBER=#1%
\multiply\NUMBER by 100%
\divide\NUMBER by \SCALE%
\setcounter{x}{\NUMBER}}%

\newcommand{\truey}[1]{%
\NUMBER=#1%
\multiply\NUMBER by 100%
\divide\NUMBER by \SCALE%
\setcounter{y}{\NUMBER}}%

\newcommand{\truez}[1]{%
\NUMBER=#1%
\multiply\NUMBER by 100%
\divide\NUMBER by \SCALE%
\setcounter{z}{\NUMBER}}%

\newcommand{\changecounters}[1]{%
\SOURCE=\ARROW%
\ARROW=\TARGET%
\settowidth{\itemlength}{#1}%
\ifdim \itemlength > 2800\unitlength%
\addtolength{\itemlength}{-2800\unitlength}%
\TARGET=\itemlength%
\divide\TARGET by 1310%
\multiply\TARGET by 100%
\divide\TARGET by \SCALE%
\else%
\TARGET=0%
\fi%
\ARROWLENGTH=5000%
\advance\ARROWLENGTH by -\SOURCE%
\advance\ARROWLENGTH by -\TARGET%
\divide\ARROWLENGTH by 100%
\advance\SOURCE by -\TARGET}%

\newcommand{\initialize}[1]{%
\LINE=0%
\COLUMN=0%
\WIDTH=0%
\ARROW=0%
\TARGET=0%
\changecounters{#1}%
\renewcommand{\printarrow}{#1}%
\begin{center}%
\vspace{2pt}%
\begin{picture}(0,0)}%

\newcommand{\DIAGV}[2]{%
\SCALE=#1%
\setlength{\unitlength}{655sp}%
\multiply\unitlength by \SCALE%
\divide\unitlength by 100%
\initialize{\mbox{$#2$}}}%

\newcommand{\n}[1]{%
\changecounters{\mbox{$#1$}}%
\put(\COLUMN,\LINE){\makebox(0,0){\printarrow}}%
\thinlines%
\renewcommand{\printarrow}{\mbox{$#1$}}%
\advance\COLUMN by 4000}%

\newcommand{\nn}[1]{%
\put(\COLUMN,\LINE){\makebox(0,0){\printarrow}}%
\thinlines%
\ifnum \WIDTH < \COLUMN%
\WIDTH=\COLUMN%
\else%
\fi%
\advance\LINE by -4000%
\COLUMN=0%
\ARROW=0%
\TARGET=0%
\changecounters{\mbox{$#1$}}%
\renewcommand{\printarrow}{\mbox{$#1$}}}%

\newcommand{\conclude}{%
\put(\COLUMN,\LINE){\makebox(0,0){\printarrow}}%
\thinlines%
\ifnum \WIDTH < \COLUMN%
\WIDTH=\COLUMN%
\else%
\fi%
\setcounter{horizontal}{\WIDTH}%
\setcounter{vertical}{-\LINE}%
\end{picture}}%

\newcommand{\diag}{%
\conclude%
\raisebox{0pt}[0pt][\value{vertical}\unitlength]{}%
\hspace*{\value{horizontal}\unitlength}%
\vspace{12pt}%
\end{center}%
\setlength{\unitlength}{1pt}%
}%

\newcommand{\diagv}[3]{%
\conclude%
\NUMBER=#1%
\rule{0pt}{\NUMBER pt}%
\hspace*{-#2pt}%
\raisebox{0pt}[0pt][\value{vertical}\unitlength]{}%
\hspace*{\value{horizontal}\unitlength}%6
\NUMBER=#3%
\advance\NUMBER by 12%
\vspace*{\NUMBER pt}%
\end{center}%
\setlength{\unitlength}{1pt}%
}%

\def\movename(#1,#2)#3{%
\hspace{#1pt}%
\raisebox{#2pt}[5pt][2pt]{\raisebox{#2pt}{$#3$}}%
\hspace{-#1pt}}%

\def\movearrow(#1,#2)#3{%
\makebox[0pt]{%
\hspace{#1pt}\hspace{#1pt}%
\raisebox{#2pt}[0pt][0pt]{\raisebox{#2pt}{$#3$}}}}%

\def\movevertex(#1,#2)#3{%
\mbox{\hspace{#1pt}%
\raisebox{#2pt}{\raisebox{#2pt}{$#3$}}%
\hspace{-#1pt}}}%

\newcommand{\movevertexleft}[1]{\makebox[2800\unitlength][r]{$#1$}}

\newcommand{\movevertexright}[1]{\makebox[2800\unitlength][l]{$#1$}}

\newcommand{\crosslength}[2]{%
\settowidth{\firstitem}{#1}%
\settowidth{\seconditem}{#2}%
\ifdim\firstitem < \seconditem%
\itemlength=\seconditem%
\else%
\itemlength=\firstitem%
\fi%
\divide\itemlength by 2%
\hspace{\itemlength}}%

\newcommand{\cross}[2]{%
\crosslength{\mbox{$#1$}}{\mbox{$#2$}}%
\begin{picture}(0,0)%
\put(0,0){\makebox(0,0){$#1$}}%
\thinlines%
\put(0,0){\makebox(0,0){$#2$}}%
\thinlines%
\end{picture}%
\crosslength{\mbox{$#1$}}{\mbox{$#2$}}}%

\def\basicDIAG#1?{\DIAGV{\defaultscale}{#1}\@ifnexttok?{\finishline}{\basicn}}
\def\basicDIAGV[#1]#2?{\DIAGV{#1}{#2}\@ifnexttok?{\finishline}{\basicn}}
\def\basicn#1?{\n{#1}\@ifnexttok?{\finishline}{\basicn}}
\def\basicnn#1?{\nn{#1}\@ifnexttok?{\finishline}{\basicn}}
\def\finishline#1{\@ifnextchar\end{\diag}%
{\@ifnextchar\spacing{\relax}{\basicnn}}}
\def\spacing(#1,#2,#3){\diagv{#1}{#2}{#3}}

\newif\ifcaption%

\newenvironment{diagram}{%
\iffloatdiag\relax\else
\global\def\diagramcaption##1{%
\global\captiontrue%
\global\def\@diagcaption{##1}}%
\global\def\@diagcaption{}\fi%
\@ifnextchar[{\basicDIAGV}{\basicDIAG}}%
{\iffloatdiag\relax\else%
\ifcaption
\begin{center}\mbox{}\@diagcaption\end{center}%
\else\relax\fi\fi\global\captionfalse}

\global\captionfalse

\gdef\@diaglabel{Diagram}
\gdef\diagramlabel#1{\gdef\@diaglabel{#1}}

\newcounter{Diagram}
\def\theDiagram{\@arabic\c@Diagram}
\def\fps@Diagram{tpbh}
\def\ftype@Diagram{1}
\def\ext@Diagram{lof}
\def\fnum@Diagram{\@diaglabel\ \theDiagram}
\def\Diagram{\@float{Diagram}}
\let\endDiagram\end@float
\@namedef{Diagram*}{\@dblfloat{Diagram}}
\@namedef{endDiagram*}{\end@dblfloat}

\newif\iffloatdiag

\newenvironment{floatingdiagram}{%
\global\floatdiagtrue%
\long\def\@makecaption##1##2{
 \vskip 0pt
 \setbox\@tempboxa\hbox{##1##2}
 \ifdim \wd\@tempboxa >\hsize \unhbox\@tempboxa\par \else \hbox
to\hsize{\hfil\box\@tempboxa\hfil}
 \fi}%
\global\def\diagramcaption##1{\global\def\@diagcaption{: ##1}}%
\global\def\@diagcaption{}%
\begin{Diagram}[htp]%
\vspace*{5pt plus 5pt}%
\begin{diagram}}%
{\end{diagram}\caption{\@diagcaption}\end{Diagram}%
\global\floatdiagfalse}

{\end{diagram}\caption{\@diagcaption}\end{Diagram}%
\global\floatdiagfalse%
\typeout{*****}%
\typeout{WARNING: Floating diagram forced to be on the page}%
\typeout{*****}
}

{\end{diagram}\caption{\@diagcaption}\end{Diagram}%
\global\floatdiagfalse%
}

\global\floatdiagfalse

\def\dlabel#1{\immediate\write\@auxout{\string\newlabel{#1}%
{{\theChapter.\@arabic\c@Diagram}{\thepage}}}}

\newcommand{\TUP}[1]{\raisebox{0pt}[0pt][3pt]{}#1}

\newcommand{\TDOWN}[1]{\raisebox{0pt}[6pt][0pt]{}#1}

\newcommand{\tlowername}[2]%
{$\stackrel{\makebox[1pt]{#1}}%
{\begin{picture}(0,0)%
\put(0,0){\makebox(0,6)[t]{\makebox[1pt]{$\scriptstyle#2$}}}%
\end{picture}}$}%

\newcommand{\tcase}[1]{%
\setlength{\unitlength}{0.01pt}%
\makebox[\actualtextarrowspace pt]%
{\raisebox{2.5pt}{#1{\actualtextarrowlength}}}%
\setlength{\unitlength}{1pt}}%

\newcommand{\Tcase}[2]{%
\setlength{\unitlength}{0.01pt}%
\makebox[\actualtextarrowspace pt]%
{\raisebox{2.5pt}{$\stackrel{\scriptstyle #2}{#1{\actualtextarrowlength}}$}}%
\setlength{\unitlength}{1pt}}%

\newcommand{\tbicase}[1]{%
\setlength{\unitlength}{0.01pt}%
\makebox[\actualtextarrowspace pt]%
{\raisebox{1pt}{#1{\actualtextarrowlength}}}%
\setlength{\unitlength}{1pt}}%

\newcommand{\Tbicase}[3]{%
\setlength{\unitlength}{0.01pt}%
\makebox[\actualtextarrowspace pt]%
{\raisebox{-1pt}%
{$\stackrel{\scriptstyle #2}%
{\mbox{\tlowername{#1{\actualtextarrowlength}}%
{\scriptstyle #3}}}$}}%
\setlength{\unitlength}{1pt}}%

\newcommand{\ttricase}[1]{%
\setlength{\unitlength}{0.01pt}%
\makebox[\actualtextarrowspace pt]%
{\raisebox{4pt}{#1{\actualtextarrowlength}%
{\rule{0pt}{6pt}}{\rule{0pt}{6pt}}{\rule{0pt}{6pt}}}}%
\setlength{\unitlength}{1pt}}%

\newcommand{\Ttricase}[4]{%
\setlength{\unitlength}{0.01pt}%
\makebox[\actualtextarrowspace pt]%
{\raisebox{4pt}%
{#1{\actualtextarrowlength}{\rule{0pt}{6pt}\scriptstyle #2}%
{\rule{0pt}{6pt}\scriptstyle #3}{\rule{0pt}{6pt}\scriptstyle #4}}}%
\setlength{\unitlength}{1pt}}%

\newcommand{\tmulticase}[1]{%
\setlength{\unitlength}{0.01pt}%
\makebox[\actualtextarrowspace pt]%
{\raisebox{1pt}{#1{\actualtextarrowlength}}}%
\setlength{\unitlength}{1pt}}%

\newcommand{\DUP}[1]{\raisebox{0pt}[0pt][4pt]{}#1}

\newcommand{\DDOWN}[1]{\raisebox{0pt}[9pt][0pt]{}#1}

\newcommand{\dlowername}[2]%
{$\stackrel{\makebox[1pt]{#1}}%
{\begin{picture}(0,0)%
\put(0,0){\makebox(0,6)[t]{\makebox[1pt]{$\textstyle#2$}}}%
\end{picture}}$}%

\newcommand{\dcase}[1]{%
\setlength{\unitlength}{0.01pt}%
\makebox[\actualdisplayarrowspace pt]%
{\raisebox{2.5pt}{#1{\actualdisplayarrowlength}}}%
\setlength{\unitlength}{1pt}}%

\newcommand{\Dcase}[2]{%
\setlength{\unitlength}{0.01pt}%
\makebox[\actualdisplayarrowspace pt]%
{\raisebox{2.5pt}{$\stackrel{\textstyle #2}{#1{\actualdisplayarrowlength}}$}}%
\setlength{\unitlength}{1pt}}%

\newcommand{\dbicase}[1]{%
\setlength{\unitlength}{0.01pt}%
\makebox[\actualdisplayarrowspace pt]%
{\raisebox{1pt}{#1{\actualdisplayarrowlength}}}%
\setlength{\unitlength}{1pt}}%

\newcommand{\Dbicase}[3]{%
\setlength{\unitlength}{0.01pt}%
\makebox[\actualdisplayarrowspace pt]%
{\raisebox{-1pt}%
{$\stackrel{\textstyle #2}%
{\mbox{\tlowername{#1{\actualdisplayarrowlength}}%
{\textstyle #3}}}$}}%
\setlength{\unitlength}{1pt}}%

\newcommand{\dtricase}[1]{%
\setlength{\unitlength}{0.01pt}%
\makebox[\actualdisplayarrowspace pt]%
{\raisebox{4pt}{#1{\actualdisplayarrowlength}%
{\rule{0pt}{6pt}}{\rule{0pt}{6pt}}{\rule{0pt}{6pt}}}}%
\setlength{\unitlength}{1pt}}%

\newcommand{\Dtricase}[4]{%
\setlength{\unitlength}{0.01pt}%
\makebox[\actualdisplayarrowspace pt]%
{\raisebox{4pt}%
{#1{\actualdisplayarrowlength}{\rule{0pt}{6pt}\textstyle #2}%
{\rule{0pt}{6pt}\textstyle #3}{\rule{0pt}{6pt}\textstyle #4}}}%
\setlength{\unitlength}{1pt}}%

\newcommand{\dmulticase}[1]{%
\setlength{\unitlength}{0.01pt}%
\makebox[\actualdisplayarrowspace pt]%
{\raisebox{1pt}{#1{\actualdisplayarrowlength}}}%
\setlength{\unitlength}{1pt}}%

\newcommand{\AR}[1]%
{\begin{picture}(#1,0)%
\put(0,0){\line(1,0){#1}}%
\put(#1,0){\ehead}%
\end{picture}}%

\newcommand{\DIST}[1]%
{\begin{picture}(#1,0)%
\put(0,0){\line(1,0){#1}}%
\put(#1,0){\ehead}%
\NUMBER=#1%
\divide\NUMBER by 2%
\put(\NUMBER,0){\distsign}%
\end{picture}}%

\newcommand{\DOTAR}[1]%
{\NUMBEROFDOTS=#1%
\divide\NUMBEROFDOTS by 300%
\advance\NUMBEROFDOTS by 1%
\begin{picture}(#1,0)%
\multiput(0,0)(300,0){\NUMBEROFDOTS}{\circle*{100}}%
\put(#1,0){\ehead}%
\end{picture}}%

\newcommand{\MONO}[1]%
{\monolength=#1%
\advance\monolength by -\monotail%
\begin{picture}(#1,0)%
\put(\monotail,0){\line(1,0){\monolength}}%
\put(#1,0){\ehead}%
\put(\monotail,0){\ehead}%
\end{picture}}%

\newcommand{\EPI}[1]%
{\epilength=#1%
\advance\epilength by -\epihead%
\begin{picture}(#1,0)(-#1,0)%
\put(-#1,0){\line(1,0){\epilength}}%
\put(-\epihead,0){\ehead}%
\put(0,0){\ehead}%
\end{picture}}%

\newcommand{\BIMO}[1]%
{\monolength=#1%
\advance\monolength by -\monotail%
\epilength=\monolength%
\advance\epilength by -\epihead%
\begin{picture}(#1,0)(-#1,0)%
\put(-\monolength,0){\line(1,0){\epilength}}%
\put(-\monolength,0){\ehead}%
\put(-\epihead,0){\ehead}%
\put(0,0){\ehead}%
\end{picture}}%

\newcommand{\BIAR}[1]%
{\begin{picture}(#1,700)%
\put(0,0){\line(1,0){#1}}%
\put(#1,0){\ehead}%
\put(0,700){\line(1,0){#1}}%
\put(#1,700){\ehead}%
\end{picture}}%

\newcommand{\BIDIST}[1]%
{\begin{picture}(#1,700)%
\put(0,0){\line(1,0){#1}}%
\put(#1,0){\ehead}%
\put(0,700){\line(1,0){#1}}%
\put(#1,700){\ehead}%
\NUMBER=#1%
\divide\NUMBER by 2%
\put(\NUMBER,0){\distsign}%
\put(\NUMBER,700){\distsign}%
\end{picture}}%

\newcommand{\EQL}[1]%
{\begin{picture}(#1,0)%
\put(0,100){\line(1,0){#1}}%
\put(0,-100){\line(1,0){#1}}%
\end{picture}}%

\newcommand{\LLINE}[1]%
{\begin{picture}(#1,0)%
\put(0,0){\line(1,0){#1}}%
\end{picture}}%

\newcommand{\ADJAR}[1]%
{\begin{picture}(#1,700)%
\put(0,0){\line(1,0){#1}}%
\put(#1,0){\ehead}%
\put(#1,700){\line(-1,0){#1}}%
\put(0,700){\whead}
\end{picture}}%

\newcommand{\ADJDIST}[1]%
{\begin{picture}(#1,700)%
\put(0,0){\line(1,0){#1}}%
\put(#1,0){\ehead}%
\put(#1,700){\line(-1,0){#1}}%
\put(0,700){\whead}
\NUMBER=#1%
\divide\NUMBER by 2%
\put(\NUMBER,0){\distsign}%
\put(\NUMBER,700){\distsign}%
\end{picture}}%

\newcommand{\TRIAR}[4]%
{%
\begin{tabular}{c}%
\mbox{$#2$}\\ \AR{#1}\\%
\mbox{$#3$}\\ \AR{#1}\\%
\mbox{$#4$}\\ \AR{#1}%
\end{tabular}%
}%

\newcommand{\TRIADJAR}[4]%
{%
\begin{tabular}{c}%
\mbox{$#2$}\\ \BKAR{#1}\\%
\mbox{$#3$}\\ \AR{#1}\\%
\mbox{$#4$}\\ \BKAR{#1}%
\end{tabular}%
}%

\newcommand{\QUADRIAR}[1]%
{%
\begin{tabular}{c}%
\AR{#1}\\%
\rule{0pt}{3pt}\\ \AR{#1}\\%
\rule{0pt}{3pt}\\ \AR{#1}\\%
\rule{0pt}{3pt}\\ \AR{#1}%
\end{tabular}%
}%

\newcommand{\QUADRIADJAR}[1]%
{%
\begin{tabular}{c}%
\BKAR{#1}\\%
\rule{0pt}{3pt}\\ \AR{#1}\\%
\rule{0pt}{3pt}\\ \BKAR{#1}\\%
\rule{0pt}{3pt}\\ \AR{#1}%
\end{tabular}%
}%

\newcommand{\QUINTIAR}[1]%
{%
\begin{tabular}{c}%
\AR{#1}\\%
\rule{0pt}{2pt}\\ \AR{#1}\\%
\rule{0pt}{2pt}\\ \AR{#1}\\%
\rule{0pt}{2pt}\\ \AR{#1}\\%
\rule{0pt}{2pt}\\ \AR{#1}%
\end{tabular%
}%

\newcommand{\QUINTIADJAR}[1]%
{%
\begin{tabular}{c}%
\AR{#1}\\%
\rule{0pt}{2pt}\\ \BKAR{#1}\\%
\rule{0pt}{2pt}\\ \AR{#1}\\%
\rule{0pt}{2pt}\\ \BKAR{#1}\\%
\rule{0pt}{2pt}\\ \AR{#1}%
\end{tabular}}%
}%

\newcommand{\ar}{\ifinner\tcase{\AR}\else\dcase{\AR}\fi}%

\newcommand{\Ar}[1]{\ifinner\Tcase{\AR}{#1}\else\Dcase{\AR}{#1}\fi}%

\newcommand{\dist}{\ifinner\tcase{\DIST}\else\dcase{\DIST}\fi}%

\newcommand{\Dist}[1]{\ifinner\Tcase{\DIST}{\TUP{#1}}%
\else\Dcase{\DIST}{\TUP{#1}}\fi}%

\newcommand{\dotar}{\ifinner\tcase{\DOTAR}\else\dcase{\DOTAR}\fi}%

\newcommand{\Dotar}[1]{\ifinner\Tcase{\DOTAR}{#1}%
\else\Dcase{\DOTAR}{#1}\fi}%

\newcommand{\mono}{\ifinner\tcase{\MONO}\else\dcase{\MONO}\fi}%

\newcommand{\Mono}[1]{\ifinner\Tcase{\MONO}{#1}\else\Dcase{\MONO}{#1}\fi}%

\newcommand{\epi}{\ifinner\tcase{\EPI}\else\dcase{\EPI}\fi}%

\newcommand{\Epi}[1]{\ifinner\Tcase{\EPI}{#1}\else\Dcase{\EPI}{#1}\fi}%

\newcommand{\bimo}{\ifinner\tcase{\BIMO}\else\dcase{\BIMO}\fi}%

\newcommand{\Bimo}[1]{\ifinner\Tcase{\BIMO}{#1}%
\else\Dcase{\BIMO}{#1}\fi}%

\newcommand{\iso}{\ifinner\Tcase{\AR}{\isosign}\else\Dcase{\AR}{\isosign}\fi}%

\newcommand{\Iso}[1]{\ifinner\Tcase{\AR}{\isosign{#1}}%
\else\Dcase{\AR}{\isosign{#1}}\fi}%

\newcommand{\biar}{\ifinner\tbicase{\BIAR}\else\dbicase{\BIAR}\fi}%

\newcommand{\Biar}[2]{\ifinner\Tbicase{\BIAR}{#1}{#2}%
\else\Dbicase{\BIAR}{#1}{#2}\fi}%

\newcommand{\bidist}{\ifinner\tbicase{\BIDIST}\else\dbicase{\BIDIST}\fi}%

\newcommand{\Bidist}[2]{\ifinner\Tbicase{\BIDIST}{\TUP{#1}}{\TDOWN{#2}}%
\else\Dbicase{\BIDIST}{\DUP{#1}}{\DDOWN{#2}}\fi}%

\newcommand{\eql}{\ifinner\tcase{\EQL}\else\dcase{\EQL}\fi}%

\newcommand{\Eql}[1]{\ifinner\Tcase{\EQL}{\TUP{#1}}%
\else\Dcase{\EQL}{\DUP{#1}}\fi}%

\newcommand{\lline}{\ifinner\tcase{\LLINE}\else\dcase{\LLINE}\fi}%

\newcommand{\Line}[1]{\ifinner\Tcase{\LLINE}{\TUP{#1}}%
\else\Dcase{\LLINE}{\DUP{#1}}\fi}%

\newcommand{\adjar}{\ifinner\tbicase{\ADJAR}\else\dbicase{\ADJAR}\fi}%

\newcommand{\Adjar}[2]{\ifinner\Tbicase{\ADJAR}{#1}{#2}%
\else\Dbicase{\ADJAR}{#1}{#2}\fi}%

\newcommand{\adjdist}{\ifinner\tbicase{\ADJDIST}\else\dbicase{\ADJDIST}\fi}%

\newcommand{\Adjdist}[2]{\ifinner\Tbicase{\ADJDIST}{\TUP{#1}}{\TDOWN{#2}}%
\else\Dbicase{\ADJDIST}{\DUP{#1}}{\DDOWN{#2}}\fi}%

\newcommand{\triar}{\ifinner\ttricase{\TRIAR}\else\dtricase{\TRIAR}\fi}%

\newcommand{\Triar}[3]{\ifinner\Ttricase{\TRIAR}{#1}{#2}{#3}%
\else\Dtricase{\TRIAR}{#1}{#2}{#3}\fi}%

\newcommand{\triadjar}{\ifinner\ttricase{\TRIADJAR}%
\else\dtricase{\TRIADJAR}\fi}%

\newcommand{\Triadjar}[3]{\ifinner\Ttricase{\TRIADJAR}{#1}{#2}{#3}%
\else\Dtricase{\TRIADJAR}{#1}{#2}{#3}\fi}%

\newcommand{\quadriar}{\ifinner\tquadricase{\QUADRIAR}%
\else\dmulticase{\QUADRIAR}\fi}%

\newcommand{\quadriadjar}{\ifinner\tmulticase{\QUADRIADJAR}%
\else\dmulticase{\QUADRIADJAR}\fi}%

\newcommand{\quintiar}{\ifinner\tmulticase{\QUINTIAR}%
\else\dmulticase{\QUINTIAR}\fi}%

\newcommand{\quintiadjar}{\ifinner\tmulticase{\QUINTIADJAR}%
\else\dmulticase{\QUINTIADJAR}\fi}%

\newcommand{\BKAR}[1]%
{\begin{picture}(#1,0)%
\put(#1,0){\line(-1,0){#1}}%
\put(0,0){\whead}%
\end{picture}}%

\newcommand{\BKDIST}[1]%
{\begin{picture}(#1,0)%
\put(#1,0){\line(-1,0){#1}}%
\put(0,0){\whead}%
\NUMBER=#1%
\divide\NUMBER by 2%
\put(\NUMBER,0){\distsign}%
\end{picture}}%

\newcommand{\BKDOTAR}[1]%
{\NUMBEROFDOTS=#1%
\divide\NUMBEROFDOTS by 300%
\advance\NUMBEROFDOTS by 1%
\begin{picture}(#1,0)%
\multiput(#1,0)(-300,0){\NUMBEROFDOTS}{\circle*{100}}%
\put(0,0){\whead}%
\end{picture}}%

\newcommand{\BKMONO}[1]%
{\monolength=#1%
\advance\monolength by -\monotail%
\begin{picture}(#1,0)(-#1,0)%
\put(-\monotail,0){\line(-1,0){\monolength}}%
\put(-\monotail,0){\whead}%
\put(-#1,0){\whead}%
\end{picture}}%

\newcommand{\BKEPI}[1]%
{\epilength=#1%
\advance\epilength by -\epihead%
\begin{picture}(#1,0)%
\put(#1,0){\line(-1,0){\epilength}}%
\put(\epihead,0){\whead}%
\put(0,0){\whead}%
\end{picture}}%

\newcommand{\BKBIMO}[1]%
{\monolength=#1%
\advance\monolength by -\monotail%
\epilength=\monolength%
\advance\epilength by -\epihead%
\begin{picture}(#1,0)%
\put(\monolength,0){\line(-1,0){\epilength}}%
\put(\monolength,0){\whead}%
\put(\epihead,0){\whead}%
\put(0,0){\whead}%
\end{picture}}%

\newcommand{\BKBIAR}[1]%
{\begin{picture}(#1,700)%
\put(#1,0){\line(-1,0){#1}}%
\put(0,0){\whead}%
\put(#1,700){\line(-1,0){#1}}%
\put(0,700){\whead}%
\end{picture}}%

\newcommand{\BKBIDIST}[1]%
{\begin{picture}(#1,700)%
\put(#1,0){\line(-1,0){#1}}%
\put(0,0){\whead}%
\put(#1,700){\line(-1,0){#1}}%
\put(0,700){\whead}%
\NUMBER=#1%
\divide\NUMBER by 2%
\put(\NUMBER,0){\distsign}%
\put(\NUMBER,700){\distsign}%
\end{picture}}%

\newcommand{\BKADJAR}[1]%
{\begin{picture}(#1,700)%
\put(0,700){\line(1,0){#1}}%
\put(#1,700){\ehead}%
\put(#1,0){\line(-1,0){#1}}%
\put(0,0){\whead}%
\end{picture}}%

\newcommand{\BKADJDIST}[1]%
{\begin{picture}(#1,700)%
\put(0,700){\line(1,0){#1}}%
\put(#1,700){\ehead}%
\put(#1,0){\line(-1,0){#1}}%
\put(0,0){\whead}%
\NUMBER=#1%
\divide\NUMBER by 2%
\put(\NUMBER,0){\distsign}%
\put(\NUMBER,700){\distsign}%
\end{picture}}%

\newcommand{\BKTRIAR}[4]%
{%
\begin{tabular}{c}%
\mbox{$#2$}\\ \BKAR{#1}\\%
\mbox{$#3$}\\ \BKAR{#1}\\%
\mbox{$#4$}\\ \BKAR{#1}%
\end{tabular}%
}%

\newcommand{\BKTRIADJAR}[4]%
{%
\begin{tabular}{c}%
\mbox{$#2$}\\ \AR{#1}\\%
\mbox{$#3$}\\ \BKAR{#1}\\%
\mbox{$#4$}\\ \AR{#1}%
\end{tabular}%
}%

\newcommand{\BKQUADRIAR}[1]%
{%
\begin{tabular}{c}%
\BKAR{#1}\\%
\rule{0pt}{3pt}\\ \BKAR{#1}\\%
\rule{0pt}{3pt}\\ \BKAR{#1}\\%
\rule{0pt}{3pt}\\ \BKAR{#1}%
\end{tabular}%
}%

\newcommand{\BKQUADRIADJAR}[1]%
{%
\begin{tabular}{c}%
\AR{#1}\\%
\rule{0pt}{3pt}\\ \BKAR{#1}\\%
\rule{0pt}{3pt}\\ \AR{#1}\\%
\rule{0pt}{3pt}\\ \BKAR{#1}%
\end{tabular}%
}%

\newcommand{\BKQUINTIAR}[1]%
{%
\begin{tabular}{c}%
\BKAR{#1}\\%
\rule{0pt}{2pt}\\ \BKAR{#1}\\%
\rule{0pt}{2pt}\\ \BKAR{#1}\\%
\rule{0pt}{2pt}\\ \BKAR{#1}\\%
\rule{0pt}{2pt}\\ \BKAR{#1}%
\end{tabular}%
}%

\newcommand{\BKQUINTIADJAR}[1]%
{%
\begin{tabular}{c}%
\BKAR{#1}\\%
\rule{0pt}{2pt}\\ \AR{#1}\\%
\rule{0pt}{2pt}\\ \BKAR{#1}\\%
\rule{0pt}{2pt}\\ \AR{#1}\\%
\rule{0pt}{2pt}\\ \BKAR{#1}%
\end{tabular}%
}%

\newcommand{\bkar}{\ifinner\tcase{\BKAR}\else\dcase{\BKAR}\fi}%

\newcommand{\Bkar}[1]{\ifinner\Tcase{\BKAR}{#1}\else\Dcase{\BKAR}{#1}\fi}%

\newcommand{\bkdist}{\ifinner\tcase{\BKDIST}\else\dcase{\BKDIST}\fi}%

\newcommand{\Bkdist}[1]{\ifinner\Tcase{\BKDIST}{\TUP{#1}}%
\else\Dcase{\BKDIST}{\TUP{#1}}\fi}%

\newcommand{\bkdotar}{\ifinner\tcase{\BKDOTAR}\else\dcase{\BKDOTAR}\fi}%

\newcommand{\Bkdotar}[1]{\ifinner\Tcase{\BKDOTAR}{#1}%
\else\Dcase{\BKDOTAR}{#1}\fi}%

\newcommand{\bkmono}{\ifinner\tcase{\BKMONO}\else\dcase{\BKMONO}\fi}%

\newcommand{\Bkmono}[1]{\ifinner\Tcase{\BKMONO}{#1}%
\else\Dcase{\BKMONO}{#1}\fi}%

\newcommand{\bkepi}{\ifinner\tcase{\BKEPI}\else\dcase{\BKEPI}\fi}%

\newcommand{\Bkepi}[1]{\ifinner\Tcase{\BKEPI}{#1}%
\else\Dcase{\BKEPI}{#1}\fi}%

\newcommand{\bkbimo}{\ifinner\tcase{\BKBIMO}\else\dcase{\BKBIMO}\fi}%

\newcommand{\Bkbimo}[1]{\ifinner\Tcase{\BKBIMO}{\hspace{9pt}#1}%
\else\Dcase{\BKBIMO}{\hspace{9pt}#1}\fi}%

\newcommand{\bkiso}{\ifinner\Tcase{\BKAR}{\isosign}%
\else\Dcase{\BKAR}{\isosign}\fi}%

\newcommand{\Bkiso}[1]{\ifinner\Tcase{\BKAR}{\isosign{#1}}%
\else\Dcase{\BKAR}{\isosign{#1}}\fi}%

\newcommand{\bkbiar}{\ifinner\tbicase{\BKBIAR}\else\dbicase{\BKBIAR}\fi}%

\newcommand{\Bkbiar}[2]{\ifinner\Tbicase{\BKBIAR}{#1}{#2}%
\else\Dbicase{\BKBIAR}{#1}{#2}\fi}%

\newcommand{\bkbidist}{\ifinner\tbicase{\BKBIDIST}%
\else\dbicase{\BKBIDIST}\fi}%

\newcommand{\Bkbidist}[2]{\ifinner\Tbicase{\BKBIDIST}{\TUP{#1}}{\TDOWN{#2}}%
\else\Tbicase{\BKBIDIST}{\DUP{#1}}{\DDOWN{#2}}\fi}%

\newcommand{\bkadjar}{\ifinner\tbicase{\BKADJAR}%
\else\dbicase{\BKADJAR}\fi}%

\newcommand{\Bkadjar}[2]{\ifinner\Tbicase{\BKADJAR}{#1}{#2}%
\else\Dbicase{\BKADJAR}{#1}{#2}\fi}%

\newcommand{\bkadjdist}{\ifinner\tbicase{\BKADJDIST}%
\else\dbicase{\BKADJDIST}\fi}%

\newcommand{\Bkadjdist}[2]{\ifinner\Tbicase{\BKADJDIST}{\TUP{#1}}{\TDOWN{#2}}%
\else\Dbicase{\BKADJDIST}{\TUP{#1}}{\TDOWN{#2}}\fi}%

\newcommand{\bktriar}{\ifinner\ttricase{\BKTRIAR}\else\dtricase{\BKTRIAR}\fi}%

\newcommand{\Bktriar}[3]{\ifinner\Ttricase{\BKTRIAR}{#1}{#2}{#3}%
\else\Dtricase{\BKTRIAR}{#1}{#2}{#3}\fi}%

\newcommand{\bktriadjar}{\ifinner\ttricase{\BKTRIADJAR}%
\else\dtricase{\BKTRIADJAR}\fi}%

\newcommand{\Bktriadjar}[3]{\ifinner\Ttricase{\BKTRIADJAR}{#1}{#2}{#3}%
\else\Dtricase{\BKTRIADJAR}{#1}{#2}{#3}\fi}%

\newcommand{\bkquadriar}{\ifinner\tquadricase{\BKQUADRIAR}%
\else\dmulticase{\BKQUADRIAR}\fi}%

\newcommand{\bkquadriadjar}{\ifinner\tmulticase{\BKQUADRIADJAR}%
\else\dmulticase{\BKQUADRIADJAR}\fi}%

\newcommand{\bkquintiar}{\ifinner\tmulticase{\BKQUINTIAR}%
\else\dmulticase{\BKQUINTIAR}\fi}%

\newcommand{\bkquintiadjar}{\ifinner\tmulticase{\BKQUINTIADJAR}%
\else\dmulticase{\BKQUINTIADJAR}\fi}%

\def\maketopbox#1{\vbox to 0\unitlength%
{\vfill\hbox to 0\unitlength{\hss#1\hss}\vss}}

\def\makebottombox#1{\vbox to 0\unitlength%
{\vss\hbox to 0\unitlength{\hss#1\hss}\vfill}}

\newcommand{\hcase}[2]%
{\makebox[0pt]%
{\raisebox{0pt}[0pt][0pt]%
{\begin{picture}(0,0)%
\put(0,0){\makebox(0,0){#1{#2}}}%
\end{picture}%
}}}%

\newcommand{\Hcase}[3]%
{\makebox[0pt]%
{\raisebox{0pt}[0pt][0pt]%
{\begin{picture}(0,0)%
\put(0,0){\makebox(0,0){#1{#3}}}%
\truex{200}%
\put(0,\value{x}){\makebottombox{$#2$}}%
\end{picture}%
}}}%

\newcommand{\hcasE}[3]%
{\makebox[0pt]%
{\raisebox{0pt}[0pt][0pt]%
{\begin{picture}(0,0)%
\put(0,0){\makebox(0,0){#1{#3}}}%
\truex{200}%
\put(0,-\value{x}){\maketopbox{$#2$}}%
\end{picture}%
}}}%

\newcommand{\Hisocase}[4]%
{\makebox[0pt]
{\raisebox{0pt}[0pt][0pt]%
{\begin{picture}(0,0)%
\put(0,0){\makebox(0,0){#1{#4}}}%
\truex{200}%
\put(0,\value{x}){\makebottombox{$#2$}}%
\put(0,-\value{x}){\maketopbox{$#3$}}%
\end{picture}%
}}}%

\newcommand{\hbicase}[2]%
{\makebox[0pt]%
{\raisebox{0pt}[0pt][0pt]%
{\begin{picture}(0,0)%
\put(0,0){\makebox(0,0){#1{#2}}}%
\end{picture}%
}}}%

\newcommand{\Hbicase}[4]%
{\makebox[0pt]
{\raisebox{0pt}[0pt][0pt]%
{\begin{picture}(0,0)%
\put(0,0){\makebox(0,0){#1{#4}}}%
\truex{700}%
\put(0,\value{x}){\makebottombox{$#2$}}%
\put(0,-\value{x}){\maketopbox{$#3$}}%
\end{picture}%
}}}%

\newcommand{\htricase}[2]%
{\makebox[0pt]%
{\raisebox{-12pt}[0pt][0pt]{#1{#2}{}{}{}{}}}}%

\newcommand{\Htricase}[6]%
{\makebox[0pt]%
{\raisebox{-12pt}[0pt][0pt]{#1{#6}{#2}{#3}{#4}{#5}}}}%

\newcommand{\hmulticase}[2]%
{\makebox[0pt]{\raisebox{-12pt}[0pt][0pt]{#1{#2}}}}%

\newcommand{\EAR}[1]%
{\begin{picture}(#1,0)%
\put(0,0){\line(1,0){#1}}%
\put(#1,0){\ehead}%
\end{picture}}%

\newcommand{\EDIST}[1]%
{\begin{picture}(#1,0)%
\put(0,0){\line(1,0){#1}}%
\put(#1,0){\ehead}%
\truex{400}
\NUMBER=#1%
\divide\NUMBER by 2%
\put(\NUMBER,0){\distsign}
\end{picture}}%

\newcommand{\EDOTAR}[1]%
{\truex{100}\truey{300}%
\NUMBEROFDOTS=#1%
\divide\NUMBEROFDOTS by \value{y}%
\advance\NUMBEROFDOTS by 1%
\begin{picture}(#1,0)%
\multiput(0,0)(\value{y},0){\NUMBEROFDOTS}%
{\circle*{\value{x}}}%
\put(#1,0){\ehead}%
\end{picture}}%

\newcommand{\EMONO}[1]%
{\truetail
\monolength=#1%
\advance\monolength by -\truemonotail%
\begin{picture}(#1,0)%
\put(\truemonotail,0){\line(1,0){\monolength}}%
\put(#1,0){\ehead}%
\put(\truemonotail,0){\ehead}%
\end{picture}}%

\newcommand{\EEPI}[1]%
{\truehead%
\epilength=#1%
\advance\epilength by -\trueepihead%
\begin{picture}(#1,0)(-#1,0)%
\put(-#1,0){\line(1,0){\epilength}}%
\put(-\trueepihead,0){\ehead}%
\put(0,0){\ehead}%
\end{picture}}%

\newcommand{\EBIMO}[1]%
{\truehead\truetail%
\monolength=#1%
\advance\monolength by -\truemonotail%
\epilength=\monolength%
\advance\epilength by -\trueepihead%
\begin{picture}(#1,0)(-#1,0)%
\put(-\monolength,0){\line(1,0){\epilength}}%
\put(-\monolength,0){\ehead}%
\put(-\trueepihead,0){\ehead}%
\put(0,0){\ehead}%
\end{picture}}%

\newcommand{\EBIAR}[1]%
{\truex{700}%
\begin{picture}(#1,\value{x})%
\put(0,0){\line(1,0){#1}}%
\put(#1,0){\ehead}%
\put(0,\value{x}){\line(1,0){#1}}%
\put(#1,\value{x}){\ehead}%
\end{picture}}%

\newcommand{\EBIDIST}[1]%
{\truex{700}%
\begin{picture}(#1,\value{x})%
\put(0,0){\line(1,0){#1}}%
\put(#1,0){\ehead}%
\put(0,\value{x}){\line(1,0){#1}}%
\put(#1,\value{x}){\ehead}%
\NUMBER=#1%
\divide\NUMBER by 2%
\put(\NUMBER,0){\distsign}
\put(\NUMBER,\value{x}){\distsign}%
\end{picture}}%

\newcommand{\EEQL}[1]%
{\begin{picture}(#1,0)%
\truex{200}%
\put(0,\value{x}){\line(1,0){#1}}%
\put(0,0){\line(1,0){#1}}%
\end{picture}}%

\newcommand{\ELINE}[1]%
{\begin{picture}(#1,0)%
\truex{100}%
\put(0,\value{x}){\line(1,0){#1}}%
\end{picture}}%

\newcommand{\EADJAR}[1]%
{\truex{700}%
\begin{picture}(#1,\value{x})%
\put(0,0){\line(1,0){#1}}%
\put(#1,0){\ehead}%
\put(#1,\value{x}){\line(-1,0){#1}}%
\put(0,\value{x}){\whead}%
\end{picture}}%

\newcommand{\EADJDIST}[1]%
{\truex{700}%
\begin{picture}(#1,\value{x})%
\put(0,0){\line(1,0){#1}}%
\put(#1,0){\ehead}%
\put(#1,\value{x}){\line(-1,0){#1}}%
\put(0,\value{x}){\whead}%
\NUMBER=#1%
\divide\NUMBER by 2%
\put(\NUMBER,0){\distsign}
\put(\NUMBER,\value{x}){\distsign}%
\end{picture}}%

\newcommand{\ETRIAR}[5]%
{\truex{1200}%
\X=\value{x}%
\multiply\X by 2%
\begin{picture}(#1,\X)%
\put(0,0){\line(1,0){#1}}%
\put(#1,0){\ehead}%
\put(0,\value{x}){\line(1,0){#1}}%
\put(#1,\value{x}){\ehead}%
\put(0,\X){\line(1,0){#1}}%
\put(#1,\X){\ehead}%
\X=#1%
\divide\X by 2%
\truex{-600}%
\put(\X,\value{x}){\makebox(0,0){$#5$}}%
\truex{600}%
\put(\X,\value{x}){\makebox(0,0){$#4$}}%
\truex{1800}%
\put(\X,\value{x}){\makebox(0,0){$#3$}}%
\truex{3000}%
\put(\X,\value{x}){\makebox(0,0){$#2$}}%
\end{picture}}%

\newcommand{\ETRIADJAR}[5]%
{\truex{1200}%
\X=\value{x}%
\multiply\X by 2%
\begin{picture}(#1,\X)%
\put(#1,0){\line(-1,0){#1}}%
\put(0,0){\whead}%
\put(0,\value{x}){\line(1,0){#1}}%
\put(#1,\value{x}){\ehead}%
\put(#1,\X){\line(-1,0){#1}}%
\put(0,\X){\whead}%
\X=#1%
\divide\X by 2%
\truex{-600}%
\put(\X,\value{x}){\makebox(0,0){$#5$}}%
\truex{600}%
\put(\X,\value{x}){\makebox(0,0){$#4$}}%
\truex{1800}%
\put(\X,\value{x}){\makebox(0,0){$#3$}}%
\truex{3000}%
\put(\X,\value{x}){\makebox(0,0){$#2$}}%
\end{picture}}%

\newcommand{\EQUADRIAR}[1]%
{\truex{800}%
\truey{1600}%
\X=\value{x}%
\multiply\X by 3%
\begin{picture}(#1,\X)%
\put(0,0){\line(1,0){#1}}%
\put(#1,0){\ehead}%
\put(0,\value{x}){\line(1,0){#1}}%
\put(#1,\value{x}){\ehead}%
\put(0,\value{y}){\line(1,0){#1}}%
\put(#1,\value{y}){\ehead}%
\put(0,\X){\line(1,0){#1}}%
\put(#1,\X){\ehead}%
\end{picture}}%

\newcommand{\EQUADRIADJAR}[1]%
{\truex{800}%
\truey{1600}%
\X=\value{x}%
\multiply\X by 3%
\begin{picture}(#1,\X)%
\put(0,0){\line(1,0){#1}}%
\put(#1,0){\ehead}%
\put(#1,\value{x}){\line(-1,0){#1}}%
\put(0,\value{x}){\whead}%
\put(0,\value{y}){\line(1,0){#1}}%
\put(#1,\value{y}){\ehead}%
\put(#1,\X){\line(-1,0){#1}}%
\put(0,\X){\whead}%
\end{picture}}%

\newcommand{\EQUINTIAR}[1]%
{\truex{600}%
\truey{1200}%
\truez{1800}%
\X=\value{x}%
\multiply\X by 4%
\begin{picture}(#1,\X)%
\put(0,0){\line(1,0){#1}}%
\put(#1,0){\ehead}%
\put(0,\value{x}){\line(1,0){#1}}%
\put(#1,\value{x}){\ehead}%
\put(0,\value{y}){\line(1,0){#1}}%
\put(#1,\value{y}){\ehead}%
\put(0,\value{z}){\line(1,0){#1}}%
\put(#1,\value{z}){\ehead}%
\put(0,\X){\line(1,0){#1}}%
\put(#1,\X){\ehead}%
\end{picture}}%

\newcommand{\EQUINTIADJAR}[1]%
{\truex{600}%
\truey{1200}%
\truez{1800}%
\X=\value{x}%
\multiply\X by 4%
\begin{picture}(#1,\X)%
\put(0,0){\line(1,0){#1}}%
\put(#1,0){\ehead}%
\put(#1,\value{x}){\line(-1,0){#1}}%
\put(0,\value{x}){\whead}%
\put(0,\value{y}){\line(1,0){#1}}%
\put(#1,\value{y}){\ehead}%
\put(#1,\value{z}){\line(-1,0){#1}}%
\put(0,\value{z}){\whead}%
\put(0,\X){\line(1,0){#1}}%
\put(#1,\X){\ehead}%
\end{picture}}%

\def\basicear[#1]{%
\Z=#1%
\multiply \Z by 100%
\hcase{\EAR}{\Z}}%

\newcommand{\ear}{\@ifnextchar[{\basicear}%
{\hspace{\SOURCE\unitlength}\basicear[\ARROWLENGTH]}}%

\def\basicEar[#1]#2{%
\Z=#1%
\multiply \Z by 100%
\Hcase{\EAR}{#2}{\Z}}%

\newcommand{\Ear}{\@ifnextchar[{\basicEar}%
{\hspace{\SOURCE\unitlength}\basicEar[\ARROWLENGTH]}}%

\def\basiceaR[#1]#2{%
\Z=#1%
\multiply \Z by 100%
\hcasE{\EAR}{#2}{\Z}}%

\newcommand{\eaR}{\@ifnextchar[{\basiceaR}%
{\hspace{\SOURCE\unitlength}\basiceaR[\ARROWLENGTH]}}%

\def\basicedist[#1]{%
\Z=#1%
\multiply \Z by 100%
\hcase{\EDIST}{\Z}}%

\newcommand{\edist}{\@ifnextchar[{\basicedist}%
{\hspace{\SOURCE\unitlength}\basicedist[\ARROWLENGTH]}}%

\def\basicEdist[#1]#2{%
\Z=#1%
\multiply \Z by 100%
\Hcase{\EDIST}{\DUP{#2}}{\Z}}%

\newcommand{\Edist}{\@ifnextchar[{\basicEdist}%
{\hspace{\SOURCE\unitlength}\basicEdist[\ARROWLENGTH]}}%

\def\basicedisT[#1]#2{%
\Z=#1%
\multiply \Z by 100%
\hcasE{\EDIST}{\DDOWN{#2}}{\Z}}%

\newcommand{\edisT}{\@ifnextchar[{\basicedisT}%
{\hspace{\SOURCE\unitlength}\basicedisT[\ARROWLENGTH]}}%

\def\basicedotar[#1]{%
\Z=#1%
\multiply \Z by 100%
\hcase{\EDOTAR}{\Z}}%

\newcommand{\edotar}{\@ifnextchar[{\basicedotar}%
{\hspace{\SOURCE\unitlength}\basicedotar[\ARROWLENGTH]}}%

\def\basicEdotar[#1]#2{%
\Z=#1%
\multiply \Z by 100%
\Hcase{\EDOTAR}{#2}{\Z}}%

\newcommand{\Edotar}{\@ifnextchar[{\basicEdotar}%
{\hspace{\SOURCE\unitlength}\basicEdotar[\ARROWLENGTH]}}%

\def\basicedotaR[#1]#2{%
\Z=#1%
\multiply \Z by 100%
\hcasE{\EDOTAR}{#2}{\Z}}%

\newcommand{\edotaR}{\@ifnextchar[{\basicedotaR}%
{\hspace{\SOURCE\unitlength}\basicedotaR[\ARROWLENGTH]}}%

\def\basicemono[#1]{%
\Z=#1%
\multiply \Z by 100%
\hcase{\EMONO}{\Z}}%

\newcommand{\emono}{\@ifnextchar[{\basicemono}%
{\hspace{\SOURCE\unitlength}\basicemono[\ARROWLENGTH]}}%

\def\basicEmono[#1]#2{%
\Z=#1%
\multiply \Z by 100%
\Hcase{\EMONO}{#2}{\Z}}%

\newcommand{\Emono}{\@ifnextchar[{\basicEmono}%
{\hspace{\SOURCE\unitlength}\basicEmono[\ARROWLENGTH]}}%

\def\basicemonO[#1]#2{%
\Z=#1%
\multiply \Z by 100%
\hcasE{\EMONO}{#2}{\Z}}%

\newcommand{\emonO}{\@ifnextchar[{\basicemonO}%
{\hspace{\SOURCE\unitlength}\basicemonO[\ARROWLENGTH]}}%

\def\basiceepi[#1]{%
\Z=#1%
\multiply \Z by 100%
\hcase{\EEPI}{\Z}}%

\newcommand{\eepi}{\@ifnextchar[{\basiceepi}%
{\hspace{\SOURCE\unitlength}\basiceepi[\ARROWLENGTH]}}%

\def\basicEepi[#1]#2{%
\Z=#1%
\multiply \Z by 100%
\Hcase{\EEPI}{#2}{\Z}}%

\newcommand{\Eepi}{\@ifnextchar[{\basicEepi}%
{\hspace{\SOURCE\unitlength}\basicEepi[\ARROWLENGTH]}}%

\def\basiceepI[#1]#2{%
\Z=#1%
\multiply \Z by 100%
\hcasE{\EEPI}{#2}{\Z}}%

\newcommand{\eepI}{\@ifnextchar[{\basiceepI}%
{\hspace{\SOURCE\unitlength}\basiceepI[\ARROWLENGTH]}}%

\def\basicebimo[#1]{%
\Z=#1%
\multiply \Z by 100%
\hcase{\EBIMO}{\Z}}%

\newcommand{\ebimo}{\@ifnextchar[{\basicebimo}%
{\hspace{\SOURCE\unitlength}\basicebimo[\ARROWLENGTH]}}%

\def\basicEbimo[#1]#2{%
\Z=#1%
\multiply \Z by 100%
\Hcase{\EBIMO}{#2}{\Z}}%

\newcommand{\Ebimo}{\@ifnextchar[{\basicEbimo}%
{\hspace{\SOURCE\unitlength}\basicEbimo[\ARROWLENGTH]}}%

\def\basicebimO[#1]#2{%
\Z=#1%
\multiply \Z by 100%
\hcasE{\EBIMO}{#2}{\Z}}%

\newcommand{\ebimO}{\@ifnextchar[{\basicebimO}%
{\hspace{\SOURCE\unitlength}\basicebimO[\ARROWLENGTH]}}%

\def\basiceiso[#1]{%
\Z=#1%
\multiply \Z by 100%
\Hisocase{\EAR}{\isosign}{}{\Z}}%

\newcommand{\eiso}{\@ifnextchar[{\basiceiso}%
{\hspace{\SOURCE\unitlength}\basiceiso[\ARROWLENGTH]}}%

\def\basicEiso[#1]#2{%
\Z=#1%
\multiply \Z by 100%
\Hisocase{\EAR}{#2}{\isosign}{\Z}}%

\newcommand{\Eiso}{\@ifnextchar[{\basicEiso}%
{\hspace{\SOURCE\unitlength}\basicEiso[\ARROWLENGTH]}}%

\def\basiceisO[#1]#2{%
\Z=#1%
\multiply \Z by 100%
\Hisocase{\EAR}{\isosign}{#2}{\Z}}%

\newcommand{\eisO}{\@ifnextchar[{\basiceisO}%
{\hspace{\SOURCE\unitlength}\basiceisO[\ARROWLENGTH]}}%

\def\basiceeql[#1]{%
\Z=#1%
\multiply \Z by 100%
\hcase{\EEQL}{\Z}}%

\newcommand{\eeql}{\@ifnextchar[{\basiceeql}%
{\hspace{\SOURCE\unitlength}\basiceeql[\ARROWLENGTH]}}%

\def\basicEeql[#1]#2{%
\Z=#1%
\multiply \Z by 100%
\Hcase{\EEQL}{\DUP{#2}}{\Z}}%

\newcommand{\Eeql}{\@ifnextchar[{\basicEeql}%
{\hspace{\SOURCE\unitlength}\basicEeql[\ARROWLENGTH]}}%

\def\basiceeqL[#1]#2{%
\Z=#1%
\multiply \Z by 100%
\hcasE{\EEQL}{#2}{\Z}}%

\newcommand{\eeqL}{\@ifnextchar[{\basiceeqL}%
{\hspace{\SOURCE\unitlength}\basiceeqL[\ARROWLENGTH]}}%

\def\basiceline[#1]{%
\Z=#1%
\multiply \Z by 100%
\hcase{\ELINE}{\Z}}%

\newcommand{\eline}{\@ifnextchar[{\basiceline}%
{\hspace{\SOURCE\unitlength}\basiceline[\ARROWLENGTH]}}%

\def\basicEline[#1]#2{%
\Z=#1%
\multiply \Z by 100%
\Hcase{\ELINE}{\DUP{#2}}{\Z}}%

\newcommand{\Eline}{\@ifnextchar[{\basicEline}%
{\hspace{\SOURCE\unitlength}\basicEline[\ARROWLENGTH]}}%

\def\basicelinE[#1]#2{%
\Z=#1%
\multiply \Z by 100%
\hcasE{\ELINE}{#2}{\Z}}%

\newcommand{\elinE}{\@ifnextchar[{\basicelinE}%
{\hspace{\SOURCE\unitlength}\basicelinE[\ARROWLENGTH]}}%

\def\basicebiar[#1]{%
\Z=#1%
\multiply \Z by 100%
\hbicase{\EBIAR}{\Z}}%

\newcommand{\ebiar}{\@ifnextchar[{\basicebiar}%
{\hspace{\SOURCE\unitlength}\basicebiar[\ARROWLENGTH]}}%

\def\basicEbiar[#1]#2#3{%
\Z=#1%
\multiply \Z by 100%
\Hbicase{\EBIAR}{#2}{#3}{\Z}}%

\newcommand{\Ebiar}{\@ifnextchar[{\basicEbiar}%
{\hspace{\SOURCE\unitlength}\basicEbiar[\ARROWLENGTH]}}%

\def\basicebidist[#1]{%
\Z=#1%
\multiply \Z by 100%
\hbicase{\EBIDIST}{\Z}}%

\newcommand{\ebidist}{\@ifnextchar[{\basicebidist}%
{\hspace{\SOURCE\unitlength}\basicebidist[\ARROWLENGTH]}}%

\def\basicEbidist[#1]#2#3{%
\Z=#1%
\multiply \Z by 100%
\Hbicase{\EBIDIST}{\DUP{#2}}{\DDOWN{#3}}{\Z}}%

\newcommand{\Ebidist}{\@ifnextchar[{\basicEbidist}%
{\hspace{\SOURCE\unitlength}\basicEbidist[\ARROWLENGTH]}}%

\def\basiceadjar[#1]{%
\Z=#1%
\multiply \Z by 100%
\hbicase{\EADJAR}{\Z}}%

\newcommand{\eadjar}{\@ifnextchar[{\basiceadjar}%
{\hspace{\SOURCE\unitlength}\basiceadjar[\ARROWLENGTH]}}%

\def\basicEadjar[#1]#2#3{%
\Z=#1%
\multiply \Z by 100%
\Hbicase{\EADJAR}{#2}{#3}{\Z}}%

\newcommand{\Eadjar}{\@ifnextchar[{\basicEadjar}%
{\hspace{\SOURCE\unitlength}\basicEadjar[\ARROWLENGTH]}}%

\def\basiceadjdist[#1]{%
\Z=#1%
\multiply \Z by 100%
\hbicase{\EADJDIST}{\Z}}%

\newcommand{\eadjdist}{\@ifnextchar[{\basiceadjdist}%
{\hspace{\SOURCE\unitlength}\basiceadjdist[\ARROWLENGTH]}}%

\def\basicEadjdist[#1]#2#3{%
\Z=#1%
\multiply \Z by 100%
\Hbicase{\EADJDIST}{\DUP{#2}}{\DDOWN{#3}}{\Z}}%

\newcommand{\Eadjdist}{\@ifnextchar[{\basicEadjdist}%
{\hspace{\SOURCE\unitlength}\basicEadjdist[\ARROWLENGTH]}}%

\def\basicetriar[#1]{%
\Z=#1%
\multiply \Z by 100%
\htricase{\ETRIAR}{\Z}}%

\newcommand{\etriar}{\@ifnextchar[{\basicetriar}%
{\hspace{\SOURCE\unitlength}\basicetriar[\ARROWLENGTH]}}%

\def\basicEtriar[#1]#2#3#4{%
\Z=#1%
\multiply \Z by 100%
\Htricase{\ETRIAR}{#2}{#3}{#4}{}{\Z}}%

\newcommand{\Etriar}{\@ifnextchar[{\basicEtriar}%
{\hspace{\SOURCE\unitlength}\basicEtriar[\ARROWLENGTH]}}%

\def\basicetriaR[#1]#2#3#4{%
\Z=#1%
\multiply \Z by 100%
\Htricase{\ETRIAR}{}{#2}{#3}{#4}{\Z}}%

\newcommand{\etriaR}{\@ifnextchar[{\basicetriaR}%
{\hspace{\SOURCE\unitlength}\basicetriaR[\ARROWLENGTH]}}%

\def\basicetriadjar[#1]{%
\Z=#1%
\multiply \Z by 100%
\htricase{\ETRIADJAR}{\Z}}%

\newcommand{\etriadjar}{\@ifnextchar[{\basicetriadjar}%
{\hspace{\SOURCE\unitlength}\basicetriadjar[\ARROWLENGTH]}}%

\def\basicEtriadjar[#1]#2#3#4{%
\Z=#1%
\multiply \Z by 100%
\Htricase{\ETRIADJAR}{#2}{#3}{#4}{}{\Z}}%

\newcommand{\Etriadjar}{\@ifnextchar[{\basicEtriadjar}%
{\hspace{\SOURCE\unitlength}\basicEtriadjar[\ARROWLENGTH]}}%

\def\basicetriadjaR[#1]#2#3#4{%
\Z=#1%
\multiply \Z by 100%
\Htricase{\ETRIADJAR}{}{#2}{#3}{#4}{\Z}}%

\newcommand{\etriadjaR}{\@ifnextchar[{\basicetriadjaR}%
{\hspace{\SOURCE\unitlength}\basicetriadjaR[\ARROWLENGTH]}}%

\def\basicequadriar[#1]{%
\Z=#1%
\multiply \Z by 100%
\hmulticase{\EQUADRIAR}{\Z}}%

\newcommand{\equadriar}{\@ifnextchar[{\basicequadriar}%
{\hspace{\SOURCE\unitlength}\basicequadriar[\ARROWLENGTH]}}%

\def\basicequadriadjar[#1]{%
\Z=#1%
\multiply \Z by 100%
\hmulticase{\EQUADRIADJAR}{\Z}}%

\newcommand{\equadriadjar}{\@ifnextchar[{\basicequadriadjar}%
{\hspace{\SOURCE\unitlength}\basicequadriadjar[\ARROWLENGTH]}}%

\def\basicequintiar[#1]{%
\Z=#1%
\multiply \Z by 100%
\hmulticase{\EQUINTIAR}{\Z}}%

\newcommand{\equintiar}{\@ifnextchar[{\basicequintiar}%
{\hspace{\SOURCE\unitlength}\basicequintiar[\ARROWLENGTH]}}%

\def\basicequintiadjar[#1]{%
\Z=#1%
\multiply \Z by 100%
\hmulticase{\EQUINTIADJAR}{\Z}}%

\newcommand{\equintiadjar}{\@ifnextchar[{\basicequintiadjar}%
{\hspace{\SOURCE\unitlength}\basicequintiadjar[\ARROWLENGTH]}}%

\newcommand{\WAR}[1]%
{\begin{picture}(#1,0)%
\put(#1,0){\line(-1,0){#1}}%
\put(0,0){\whead}%
\end{picture}}%

\newcommand{\WDIST}[1]%
{\begin{picture}(#1,0)%
\put(#1,0){\line(-1,0){#1}}%
\put(0,0){\whead}%
\NUMBER=#1%
\divide\NUMBER by 2%
\put(\NUMBER,0){\distsign}%
\end{picture}}%

\newcommand{\WDOTAR}[1]%
{\truex{100}\truey{300}%
\NUMBEROFDOTS=#1%
\divide\NUMBEROFDOTS by \value{y}%
\advance\NUMBEROFDOTS by 1%
\begin{picture}(#1,0)%
\multiput(#1,0)(-\value{y},0){\NUMBEROFDOTS}%
{\circle*{\value{x}}}%
\put(0,0){\whead}%
\end{picture}}%

\newcommand{\WMONO}[1]%
{\truetail%
\monolength=#1%
\advance\monolength by -\truemonotail%
\begin{picture}(#1,0)(-#1,0)%
\put(-\truemonotail,0){\line(-1,0){\monolength}}%
\put(-\truemonotail,0){\whead}%
\put(-#1,0){\whead}%
\end{picture}}%

\newcommand{\WEPI}[1]%
{\truehead%
\epilength=#1%
\advance\epilength by -\trueepihead%
\begin{picture}(#1,0)%
\put(#1,0){\line(-1,0){\epilength}}%
\put(\trueepihead,0){\whead}%
\put(0,0){\whead}%
\end{picture}}%

\newcommand{\WBIMO}[1]%
{\truehead\truetail%
\monolength=#1
\advance\monolength by -\truemonotail%
\epilength=\monolength%
\advance\epilength by -\trueepihead%
\begin{picture}(#1,0)%
\put(\monolength,0){\line(-1,0){\epilength}}%
\put(\monolength,0){\whead}%
\put(\trueepihead,0){\whead}%
\put(0,0){\whead}%
\end{picture}}%

\newcommand{\WBIAR}[1]%
{\truex{700}%
\begin{picture}(#1,\value{x})%
\put(#1,0){\line(-1,0){#1}}%
\put(0,0){\whead}%
\put(#1,\value{x}){\line(-1,0){#1}}%
\put(0,\value{x}){\whead}%
\end{picture}}%

\newcommand{\WBIDIST}[1]%
{\truex{700}%
\begin{picture}(#1,\value{x})%
\put(#1,0){\line(-1,0){#1}}%
\put(0,0){\whead}%
\put(#1,\value{x}){\line(-1,0){#1}}%
\put(0,\value{x}){\whead}%
\NUMBER=#1%
\divide\NUMBER by 2%
\put(\NUMBER,0){\distsign}%
\put(\NUMBER,\value{x}){\distsign}%
\end{picture}}%

\newcommand{\WADJAR}[1]%
{\truex{700}%
\begin{picture}(#1,\value{x})%
\put(0,\value{x}){\line(1,0){#1}}%
\put(#1,\value{x}){\ehead}%
\put(#1,0){\line(-1,0){#1}}%
\put(0,0){\whead}%
\end{picture}}%

\newcommand{\WADJDIST}[1]%
{\truex{700}%
\begin{picture}(#1,\value{x})%
\put(0,\value{x}){\line(1,0){#1}}%
\put(#1,\value{x}){\ehead}%
\put(#1,0){\line(-1,0){#1}}%
\put(0,0){\whead}%
\NUMBER=#1%
\divide\NUMBER by 2%
\put(\NUMBER,0){\distsign}%
\put(\NUMBER,\value{x}){\distsign}%
\end{picture}}%

\newcommand{\WTRIAR}[5]%
{\truex{1200}%
\X=\value{x}%
\multiply\X by 2%
\begin{picture}(#1,\X)%
\put(#1,0){\line(-1,0){#1}}%
\put(0,0){\whead}%
\put(#1,\value{x}){\line(-1,0){#1}}%
\put(0,\value{x}){\whead}%
\put(#1,\X){\line(-1,0){#1}}%
\put(0,\X){\whead}%
\X=#1%
\divide\X by 2%
\truex{-600}%
\put(\X,\value{x}){\makebox(0,0){$#5$}}%
\truex{600}%
\put(\X,\value{x}){\makebox(0,0){$#4$}}%
\truex{1800}%
\put(\X,\value{x}){\makebox(0,0){$#3$}}%
\truex{3000}%
\put(\X,\value{x}){\makebox(0,0){$#2$}}%
\end{picture}}%

\newcommand{\WTRIADJAR}[5]%
{\truex{1200}%
\X=\value{x}%
\multiply\X by 2%
\begin{picture}(#1,\X)%
\put(0,0){\line(1,0){#1}}%
\put(#1,0){\ehead}%
\put(#1,\value{x}){\line(-1,0){#1}}%
\put(0,\value{x}){\whead}%
\put(0,\X){\line(1,0){#1}}%
\put(#1,\X){\ehead}%
\X=#1%
\divide\X by 2%
\truex{-600}%
\put(\X,\value{x}){\makebox(0,0){$#5$}}%
\truex{600}%
\put(\X,\value{x}){\makebox(0,0){$#4$}}%
\truex{1800}%
\put(\X,\value{x}){\makebox(0,0){$#3$}}%
\truex{3000}%
\put(\X,\value{x}){\makebox(0,0){$#2$}}%
\end{picture}}%

\newcommand{\WQUADRIAR}[1]%
{\truex{800}%
\truey{1600}%
\X=\value{x}%
\multiply\X by 3%
\begin{picture}(#1,\X)%
\put(#1,0){\line(-1,0){#1}}%
\put(0,0){\whead}%
\put(#1,\value{x}){\line(-1,0){#1}}%
\put(0,\value{x}){\whead}%
\put(#1,\value{y}){\line(-1,0){#1}}%
\put(0,\value{y}){\whead}%
\put(#1,\X){\line(-1,0){#1}}%
\put(0,\X){\whead}%
\end{picture}}%

\newcommand{\WQUADRIADJAR}[1]%
{\truex{800}%
\truey{1600}%
\X=\value{x}%
\multiply\X by 3%
\begin{picture}(#1,\X)%
\put(#1,0){\line(-1,0){#1}}%
\put(0,0){\whead}%
\put(0,\value{x}){\line(1,0){#1}}%
\put(#1,\value{x}){\ehead}%
\put(#1,\value{y}){\line(-1,0){#1}}%
\put(0,\value{y}){\whead}%
\put(0,\X){\line(1,0){#1}}%
\put(#1,\X){\ehead}%
\end{picture}}%

\newcommand{\WQUINTIAR}[1]%
{\truex{600}%
\truey{1200}%
\truez{1800}%
\X=\value{x}%
\multiply\X by 4%
\begin{picture}(#1,\X)%
\put(#1,0){\line(-1,0){#1}}%
\put(0,0){\whead}%
\put(#1,\value{x}){\line(-1,0){#1}}%
\put(0,\value{x}){\whead}%
\put(#1,\value{y}){\line(-1,0){#1}}%
\put(0,\value{y}){\whead}%
\put(#1,\value{z}){\line(-1,0){#1}}%
\put(0,\value{z}){\whead}%
\put(#1,\X){\line(-1,0){#1}}%
\put(0,\X){\whead}%
\end{picture}}%

\newcommand{\WQUINTIADJAR}[1]%
{\truex{600}%
\truey{1200}%
\truez{1800}%
\X=\value{x}%
\multiply\X by 4%
\begin{picture}(#1,\X)%
\put(#1,0){\line(-1,0){#1}}%
\put(0,0){\whead}%
\put(0,\value{x}){\line(1,0){#1}}%
\put(#1,\value{x}){\ehead}%
\put(#1,\value{y}){\line(-1,0){#1}}%
\put(0,\value{y}){\whead}%
\put(0,\value{z}){\line(1,0){#1}}%
\put(#1,\value{z}){\ehead}%
\put(#1,\X){\line(-1,0){#1}}%
\put(0,\X){\whead}%
\end{picture}}%

\def\basicwar[#1]{%
\Z=#1%
\multiply \Z by 100%
\hcase{\WAR}{\Z}}%

\newcommand{\war}{\@ifnextchar[{\basicwar}%
{\hspace{\SOURCE\unitlength}\basicwar[\ARROWLENGTH]}}%

\def\basicWar[#1]#2{%
\Z=#1%
\multiply \Z by 100%
\Hcase{\WAR}{#2}{\Z}}%

\newcommand{\War}{\@ifnextchar[{\basicWar}%
{\hspace{\SOURCE\unitlength}\basicWar[\ARROWLENGTH]}}%

\def\basicwaR[#1]#2{%
\Z=#1%
\multiply \Z by 100%
\hcasE{\WAR}{#2}{\Z}}%

\newcommand{\waR}{\@ifnextchar[{\basicwaR}%
{\hspace{\SOURCE\unitlength}\basicwaR[\ARROWLENGTH]}}%

\def\basicwdist[#1]{%
\Z=#1%
\multiply \Z by 100%
\hcase{\WDIST}{\Z}}%

\newcommand{\wdist}{\@ifnextchar[{\basicwdist}%
{\hspace{\SOURCE\unitlength}\basicwdist[\ARROWLENGTH]}}%

\def\basicWdist[#1]#2{%
\Z=#1%
\multiply \Z by 100%
\Hcase{\WDIST}{\DUP{#2}}{\Z}}%

\newcommand{\Wdist}{\@ifnextchar[{\basicWdist}%
{\hspace{\SOURCE\unitlength}\basicWdist[\ARROWLENGTH]}}%

\def\basicwdisT[#1]#2{%
\Z=#1%
\multiply \Z by 100%
\hcasE{\WDIST}{\DDOWN{#2}}{\Z}}%

\newcommand{\wdisT}{\@ifnextchar[{\basicwdisT}%
{\hspace{\SOURCE\unitlength}\basicwdisT[\ARROWLENGTH]}}%

\def\basicwdotar[#1]{%
\Z=#1%
\multiply \Z by 100%
\hcase{\WDOTAR}{\Z}}%

\newcommand{\wdotar}{\@ifnextchar[{\basicwdotar}%
{\hspace{\SOURCE\unitlength}\basicwdotar[\ARROWLENGTH]}}%

\def\basicWdotar[#1]#2{%
\Z=#1%
\multiply \Z by 100%
\Hcase{\WDOTAR}{#2}{\Z}}%

\newcommand{\Wdotar}{\@ifnextchar[{\basicWdotar}%
{\hspace{\SOURCE\unitlength}\basicWdotar[\ARROWLENGTH]}}%

\def\basicwdotaR[#1]#2{%
\Z=#1%
\multiply \Z by 100%
\hcasE{\WDOTAR}{#2}{\Z}}%

\newcommand{\wdotaR}{\@ifnextchar[{\basicwdotaR}%
{\hspace{\SOURCE\unitlength}\basicwdotaR[\ARROWLENGTH]}}%

\def\basicwmono[#1]{%
\Z=#1%
\multiply \Z by 100%
\hcase{\WMONO}{\Z}}%

\newcommand{\wmono}{\@ifnextchar[{\basicwmono}%
{\hspace{\SOURCE\unitlength}\basicwmono[\ARROWLENGTH]}}%

\def\basicWmono[#1]#2{%
\Z=#1%
\multiply \Z by 100%
\Hcase{\WMONO}{#2}{\Z}}%

\newcommand{\Wmono}{\@ifnextchar[{\basicWmono}%
{\hspace{\SOURCE\unitlength}\basicWmono[\ARROWLENGTH]}}%

\def\basicwmonO[#1]#2{%
\Z=#1%
\multiply \Z by 100%
\hcasE{\WMONO}{#2}{\Z}}%

\newcommand{\wmonO}{\@ifnextchar[{\basicwmonO}%
{\hspace{\SOURCE\unitlength}\basicwmonO[\ARROWLENGTH]}}%

\def\basicwepi[#1]{%
\Z=#1%
\multiply \Z by 100%
\hcase{\WEPI}{\Z}}%

\newcommand{\wepi}{\@ifnextchar[{\basicwepi}%
{\hspace{\SOURCE\unitlength}\basicwepi[\ARROWLENGTH]}}%

\def\basicWepi[#1]#2{%
\Z=#1%
\multiply \Z by 100%
\Hcase{\WEPI}{#2}{\Z}}%

\newcommand{\Wepi}{\@ifnextchar[{\basicWepi}%
{\hspace{\SOURCE\unitlength}\basicWepi[\ARROWLENGTH]}}%

\def\basicwepI[#1]#2{%
\Z=#1%
\multiply \Z by 100%
\hcasE{\WEPI}{#2}{\Z}}%

\newcommand{\wepI}{\@ifnextchar[{\basicwepI}%
{\hspace{\SOURCE\unitlength}\basicwepI[\ARROWLENGTH]}}%

\def\basicwbimo[#1]{%
\Z=#1%
\multiply \Z by 100%
\hcase{\WBIMO}{\Z}}%

\newcommand{\wbimo}{\@ifnextchar[{\basicwbimo}%
{\hspace{\SOURCE\unitlength}\basicwbimo[\ARROWLENGTH]}}%

\def\basicWbimo[#1]#2{%
\Z=#1%
\multiply \Z by 100%
\Hcase{\WBIMO}{#2}{\Z}}%

\newcommand{\Wbimo}{\@ifnextchar[{\basicWbimo}%
{\hspace{\SOURCE\unitlength}\basicWbimo[\ARROWLENGTH]}}%

\def\basicwbimO[#1]#2{%
\Z=#1%
\multiply \Z by 100%
\hcasE{\WBIMO}{#2}{\Z}}%

\newcommand{\wbimO}{\@ifnextchar[{\basicwbimO}%
{\hspace{\SOURCE\unitlength}\basicwbimO[\ARROWLENGTH]}}%

\def\basicwiso[#1]{%
\Z=#1%
\multiply \Z by 100%
\Hisocase{\WAR}{\isosign}{}{\Z}}%

\newcommand{\wiso}{\@ifnextchar[{\basicwiso}%
{\hspace{\SOURCE\unitlength}\basicwiso[\ARROWLENGTH]}}%

\def\basicWiso[#1]#2{%
\Z=#1%
\multiply \Z by 100%
\Hisocase{\WAR}{#2}{\isosign}{\Z}}%

\newcommand{\Wiso}{\@ifnextchar[{\basicWiso}%
{\hspace{\SOURCE\unitlength}\basicWiso[\ARROWLENGTH]}}%

\def\basicwisO[#1]#2{%
\Z=#1%
\multiply \Z by 100%
\Hisocase{\WAR}{\isosign}{#2}{\Z}}%

\newcommand{\wisO}{\@ifnextchar[{\basicwisO}%
{\hspace{\SOURCE\unitlength}\basicwisO[\ARROWLENGTH]}}%

\def\basicwbiar[#1]{%
\Z=#1%
\multiply \Z by 100%
\hbicase{\WBIAR}{\Z}}%

\newcommand{\wbiar}{\@ifnextchar[{\basicwbiar}%
{\hspace{\SOURCE\unitlength}\basicwbiar[\ARROWLENGTH]}}%

\def\basicWbiar[#1]#2#3{%
\Z=#1%
\multiply \Z by 100%
\Hbicase{\WBIAR}{#2}{#3}{\Z}}%

\newcommand{\Wbiar}{\@ifnextchar[{\basicWbiar}%
{\hspace{\SOURCE\unitlength}\basicWbiar[\ARROWLENGTH]}}%

\def\basicwbidist[#1]{%
\Z=#1%
\multiply \Z by 100%
\hbicase{\WBIDIST}{\Z}}%

\newcommand{\wbidist}{\@ifnextchar[{\basicwbidist}%
{\hspace{\SOURCE\unitlength}\basicwbidist[\ARROWLENGTH]}}%

\def\basicWbidist[#1]#2#3{%
\Z=#1%
\multiply \Z by 100%
\Hbicase{\WBIDIST}{\DUP{#2}}{\DDOWN{#3}}{\Z}}%

\newcommand{\Wbidist}{\@ifnextchar[{\basicWbidist}%
{\hspace{\SOURCE\unitlength}\basicWbidist[\ARROWLENGTH]}}%

\def\basicwadjar[#1]{%
\Z=#1%
\multiply \Z by 100%
\hbicase{\WADJAR}{\Z}}%

\newcommand{\wadjar}{\@ifnextchar[{\basicwadjar}%
{\hspace{\SOURCE\unitlength}\basicwadjar[\ARROWLENGTH]}}%

\def\basicWadjar[#1]#2#3{%
\Z=#1%
\multiply \Z by 100%
\Hbicase{\WADJAR}{#2}{#3}{\Z}}%

\newcommand{\Wadjar}{\@ifnextchar[{\basicWadjar}%
{\hspace{\SOURCE\unitlength}\basicWadjar[\ARROWLENGTH]}}%

\def\basicwadjdist[#1]{%
\Z=#1%
\multiply \Z by 100%
\hbicase{\WADJDIST}{\Z}}%

\newcommand{\wadjdist}{\@ifnextchar[{\basicwadjdist}%
{\hspace{\SOURCE\unitlength}\basicwadjdist[\ARROWLENGTH]}}%

\def\basicWadjdist[#1]#2#3{%
\Z=#1%
\multiply \Z by 100%
\Hbicase{\WADJDIST}{\DUP{#2}}{\DDOWN{#3}}{\Z}}%

\newcommand{\Wadjdist}{\@ifnextchar[{\basicWadjdist}%
{\hspace{\SOURCE\unitlength}\basicWadjdist[\ARROWLENGTH]}}%

\def\basicwtriar[#1]{%
\Z=#1%
\multiply \Z by 100%
\htricase{\WTRIAR}{\Z}}%

\newcommand{\wtriar}{\@ifnextchar[{\basicwtriar}%
{\hspace{\SOURCE\unitlength}\basicwtriar[\ARROWLENGTH]}}%

\def\basicWtriar[#1]#2#3#4{%
\Z=#1%
\multiply \Z by 100%
\Htricase{\WTRIAR}{#2}{#3}{#4}{}{\Z}}%

\newcommand{\Wtriar}{\@ifnextchar[{\basicWtriar}%
{\hspace{\SOURCE\unitlength}\basicWtriar[\ARROWLENGTH]}}%

\def\basicwtriaR[#1]#2#3#4{%
\Z=#1%
\multiply \Z by 100%
\Htricase{\WTRIAR}{}{#2}{#3}{#4}{\Z}}%

\newcommand{\wtriaR}{\@ifnextchar[{\basicwtriaR}%
{\hspace{\SOURCE\unitlength}\basicwtriaR[\ARROWLENGTH]}}%

\def\basicwtriadjar[#1]{%
\Z=#1%
\multiply \Z by 100%
\htricase{\WTRIADJAR}{\Z}}%

\newcommand{\wtriadjar}{\@ifnextchar[{\basicwtriadjar}%
{\hspace{\SOURCE\unitlength}\basicwtriadjar[\ARROWLENGTH]}}%

\def\basicWtriadjar[#1]#2#3#4{%
\Z=#1%
\multiply \Z by 100%
\Htricase{\WTRIADJAR}{#2}{#3}{#4}{}{\Z}}%

\newcommand{\Wtriadjar}{\@ifnextchar[{\basicWtriadjar}%
{\hspace{\SOURCE\unitlength}\basicWtriadjar[\ARROWLENGTH]}}%

\def\basicwtriadjaR[#1]#2#3#4{%
\Z=#1%
\multiply \Z by 100%
\Htricase{\WTRIADJAR}{}{#2}{#3}{#4}{\Z}}%

\newcommand{\wtriadjaR}{\@ifnextchar[{\basicwtriadjaR}%
{\hspace{\SOURCE\unitlength}\basicwtriadjaR[\ARROWLENGTH]}}%

\def\basicwquadriar[#1]{%
\Z=#1%
\multiply \Z by 100%
\hmulticase{\WQUADRIAR}{\Z}}%

\newcommand{\wquadriar}{\@ifnextchar[{\basicwquadriar}%
{\hspace{\SOURCE\unitlength}\basicwquadriar[\ARROWLENGTH]}}%

\def\basicwquadriadjar[#1]{%
\Z=#1%
\multiply \Z by 100%
\hmulticase{\WQUADRIADJAR}{\Z}}%

\newcommand{\wquadriadjar}{\@ifnextchar[{\basicwquadriadjar}%
{\hspace{\SOURCE\unitlength}\basicwquadriadjar[\ARROWLENGTH]}}%

\def\basicwquintiar[#1]{%
\Z=#1%
\multiply \Z by 100%
\hmulticase{\WQUINTIAR}{\Z}}%

\newcommand{\wquintiar}{\@ifnextchar[{\basicwquintiar}%
{\hspace{\SOURCE\unitlength}\basicwquintiar[\ARROWLENGTH]}}%

\def\basicwquintiadjar[#1]{%
\Z=#1%
\multiply \Z by 100%
\hmulticase{\WQUINTIADJAR}{\Z}}%

\newcommand{\wquintiadjar}{\@ifnextchar[{\basicwquintiadjar}%
{\hspace{\SOURCE\unitlength}\basicwquintiadjar[\ARROWLENGTH]}}%

\newcommand{\vcase}[2]{#1{#2}}%

\newcommand{\Vcase}[3]{\makebox[0pt]%
{\makebox[0pt][r]{\raisebox{0pt}[0pt][0pt]{${#2}\hspace{2pt}$}}}#1{#3}}%

\newcommand{\vcasE}[3]{\makebox[0pt]%
{#1{#3}\makebox[0pt][l]{\raisebox{0pt}[0pt][0pt]{\hspace{2pt}$#2$}}}}%

\newcommand{\Visocase}[4]{\makebox[0pt]%
{\makebox[0pt][r]{\raisebox{0pt}[0pt][0pt]{$#2$\hspace{2pt}}}#1{#4}%
\makebox[0pt][l]{\raisebox{0pt}[0pt][0pt]{\hspace{2pt}$#3$}}}}%

\newcommand{\vbicase}[2]{\makebox[0pt]{{#1{#2}}}}%

\newcommand{\Vbicase}[4]{\makebox[0pt]%
{\makebox[0pt][r]{\raisebox{0pt}[0pt][0pt]{$#2$\hspace{5.5pt}}}#1{#4}%
\makebox[0pt][l]{\raisebox{0pt}[0pt][0pt]{\hspace{6.5pt}$#3$}}}}%

\newcommand{\vtricase}[2]%
{\makebox[0pt]{#1{}{}{}{}{#2}}}%

\newcommand{\Vtricase}[6]%
{\makebox[0pt]{#1{#2}{#3}{#4}{#5}{#6}}}%

\newcommand{\vmulticase}[2]{\makebox[0pt]{#1{#2}}}%

\newcommand{\SAR}[1]%
{\begin{picture}(0,0)%
\put(0,0){\makebox(0,0)%
{\begin{picture}(0,#1)%
\put(0,#1){\line(0,-1){#1}}%
\put(0,0){\shead}%
\end{picture}}}\end{picture}}%

\newcommand{\SDIST}[1]%
{\begin{picture}(0,0)%
\put(0,0){\makebox(0,0)%
{\begin{picture}(0,#1)%
\put(0,#1){\line(0,-1){#1}}%
\put(0,0){\shead}%
\end{picture}}}%
\put(0,0){\distsign}%
\end{picture}}%

\newcommand{\SDOTAR}[1]%
{\truex{100}\truey{300}%
\NUMBEROFDOTS=#1%
\divide\NUMBEROFDOTS by \value{y}%
\advance\NUMBEROFDOTS by 1%
\begin{picture}(0,0)%
\put(0,0){\makebox(0,0)%
{\begin{picture}(0,#1)%
\multiput(0,#1)(0,-\value{y}){\NUMBEROFDOTS}%
{\circle*{\value{x}}}%
\put(0,0){\shead}%
\end{picture}}}\end{picture}}%

\newcommand{\SMONO}[1]%
{\truetail%
\monolength=#1%
\advance\monolength by -\truemonotail%
\begin{picture}(0,0)%
\put(0,0){\makebox(0,0)%
{\begin{picture}(0,#1)%
\put(0,\monolength){\line(0,-1){\monolength}}%
\put(0,\monolength){\shead}%
\put(0,0){\shead}%
\end{picture}}}\end{picture}}%

\newcommand{\SEPI}[1]%
{\truehead%
\epilength=#1%
\advance\epilength by -\trueepihead%
\begin{picture}(0,0)%
\put(0,0){\makebox(0,0)%
{\begin{picture}(0,#1)%
\put(0,#1){\line(0,-1){\epilength}}%
\put(0,\trueepihead){\shead}%
\put(0,0){\shead}%
\end{picture}}}\end{picture}}%

\newcommand{\SBIMO}[1]%
{\truehead\truetail%
\monolength=#1%
\advance\monolength by -\truemonotail%
\epilength=\monolength%
\advance\epilength by -\trueepihead%
\begin{picture}(0,0)%
\put(0,0){\makebox(0,0)%
{\begin{picture}(0,#1)%
\put(0,\monolength){\line(0,-1){\epilength}}%
\put(0,\monolength){\shead}%
\put(0,\trueepihead){\shead}%
\put(0,0){\shead}%
\end{picture}}}\end{picture}}%

\newcommand{\SBIAR}[1]%
{\begin{picture}(0,0)%
\truex{350}%
\put(0,0){\makebox(0,0)%
{\begin{picture}(0,#1)%
\put(-\value{x},#1){\line(0,-1){#1}}%
\put(-\value{x},0){\shead}%
\put(\value{x},#1){\line(0,-1){#1}}%
\put(\value{x},0){\shead}%
\end{picture}}}\end{picture}}%

\newcommand{\SBIDIST}[1]%
{\begin{picture}(0,0)%
\truex{350}%
\put(0,0){\makebox(0,0)%
{\begin{picture}(0,#1)%
\put(-\value{x},#1){\line(0,-1){#1}}%
\put(-\value{x},0){\shead}%
\put(\value{x},#1){\line(0,-1){#1}}%
\put(\value{x},0){\shead}%
\end{picture}}}%
\put(-\value{x},0){\distsign}%
\put(\value{x},0){\distsign}%
\end{picture}}%

\newcommand{\SEQL}[1]%
{\begin{picture}(0,0)%
\truex{100}%
\put(0,0){\makebox(0,0)%
{\begin{picture}(0,#1)\put(-\value{x},#1){\line(0,-1){#1}}%
\put(\value{x},#1){\line(0,-1){#1}}%
\end{picture}}}\end{picture}}%

\newcommand{\SLINE}[1]%
{\begin{picture}(0,0)%
\put(0,0){\makebox(0,0)%
{\begin{picture}(0,#1)\put(0,#1){\line(0,-1){#1}}%
\end{picture}}}\end{picture}}%

\newcommand{\SADJAR}[1]{\begin{picture}(0,0)%
\truex{350}%
\put(0,0){\makebox(0,0)%
{\begin{picture}(0,#1)%
\put(-\value{x},#1){\line(0,-1){#1}}%
\put(-\value{x},0){\shead}%
\put(\value{x},0){\line(0,1){#1}}%
\put(\value{x},#1){\nhead}%
\end{picture}}}\end{picture}}%

\newcommand{\SADJDIST}[1]{\begin{picture}(0,0)%
\truex{350}%
\put(0,0){\makebox(0,0)%
{\begin{picture}(0,#1)%
\put(-\value{x},#1){\line(0,-1){#1}}%
\put(-\value{x},0){\shead}%
\put(\value{x},0){\line(0,1){#1}}%
\put(\value{x},#1){\nhead}%
\end{picture}}}%
\put(-\value{x},0){\distsign}%
\put(\value{x},0){\distsign}%
\end{picture}}%

\newcommand{\STRIAR}[5]%
{\truex{1200}%
\truey{1800}%
\truez{600}%
\begin{picture}(0,0)%
\put(0,0){\makebox(0,0)%
{\begin{picture}(0,#5)%
\X=#5\divide\X by 2%
\put(-\value{x},#5){\line(0,-1){#5}}%
\put(-\value{x},0){\shead}%
\put(0,#5){\line(0,-1){#5}}%
\put(0,0){\shead}%
\put(\value{x},#5){\line(0,-1){#5}}%
\put(\value{x},0){\shead}%
\put(-\value{y},\X){\makebox(0,0){$#1$}}%
\put(\value{y},\X){\makebox(0,0){$#4$}}%
\put(-\value{z},\X){\makebox(0,0){$#2$}}%
\put(\value{z},\X){\makebox(0,0){$#3$}}%
\end{picture}}}\end{picture}}%

\newcommand{\STRIADJAR}[5]%
{\truex{1200}%
\truey{1800}%
\truez{600}%
\begin{picture}(0,0)%
\put(0,0){\makebox(0,0)%
{\begin{picture}(0,#5)%
\X=#5\divide\X by 2%
\put(-\value{x},0){\line(0,1){#5}}%
\put(-\value{x},#5){\nhead}%
\put(0,#5){\line(0,-1){#5}}%
\put(0,0){\shead}%
\put(\value{x},0){\line(0,1){#5}}%
\put(\value{x},#5){\nhead}%
\put(-\value{y},\X){\makebox(0,0){$#1$}}%
\put(\value{y},\X){\makebox(0,0){$#4$}}%
\put(-\value{z},\X){\makebox(0,0){$#2$}}%
\put(\value{z},\X){\makebox(0,0){$#3$}}%
\end{picture}}}\end{picture}}%

\newcommand{\SQUADRIAR}[1]%
{\truex{400}%
\truey{1200}%
\begin{picture}(0,0)%
\put(0,0){\makebox(0,0)%
{\begin{picture}(0,#1)%
\put(-\value{y},#1){\line(0,-1){#1}}%
\put(-\value{y},0){\shead}%
\put(-\value{x},#1){\line(0,-1){#1}}%
\put(-\value{x},0){\shead}%
\put(\value{x},#1){\line(0,-1){#1}}%
\put(\value{x},0){\shead}%
\put(\value{y},#1){\line(0,-1){#1}}%
\put(\value{y},0){\shead}%
\end{picture}}}\end{picture}}%

\newcommand{\SQUADRIADJAR}[1]%
{\truex{400}%
\truey{1200}%
\begin{picture}(0,0)%
\put(0,0){\makebox(0,0)%
{\begin{picture}(0,#1)%
\put(-\value{y},#1){\line(0,-1){#1}}%
\put(-\value{y},0){\shead}%
\put(-\value{x},0){\line(0,1){#1}}%
\put(-\value{x},#1){\nhead}%
\put(\value{x},#1){\line(0,-1){#1}}%
\put(\value{x},0){\shead}%
\put(\value{y},0){\line(0,1){#1}}%
\put(\value{y},#1){\nhead}%
\end{picture}}}\end{picture}}%

\newcommand{\SQUINTIAR}[1]%
{\truex{600}%
\truey{1200}%
\begin{picture}(0,0)%
\put(0,0){\makebox(0,0)%
{\begin{picture}(0,#1)%
\put(-\value{y},#1){\line(0,-1){#1}}%
\put(-\value{y},0){\shead}%
\put(-\value{x},#1){\line(0,-1){#1}}%
\put(-\value{x},0){\shead}%
\put(0,#1){\line(0,-1){#1}}%
\put(0,0){\shead}%
\put(\value{x},#1){\line(0,-1){#1}}%
\put(\value{x},0){\shead}%
\put(\value{y},#1){\line(0,-1){#1}}%
\put(\value{y},0){\shead}%
\end{picture}}}\end{picture}}%

\newcommand{\SQUINTIADJAR}[1]%
{\truex{600}%
\truey{1200}%
\begin{picture}(0,0)%
\put(0,0){\makebox(0,0)%
{\begin{picture}(0,#1)%
\put(-\value{y},#1){\line(0,-1){#1}}%
\put(-\value{y},0){\shead}%
\put(-\value{x},0){\line(0,1){#1}}%
\put(-\value{x},#1){\nhead}%
\put(0,#1){\line(0,-1){#1}}%
\put(0,0){\shead}%
\put(\value{x},0){\line(0,1){#1}}%
\put(\value{x},#1){\nhead}%
\put(\value{y},#1){\line(0,-1){#1}}%
\put(\value{y},0){\shead}%
\end{picture}}}\end{picture}}%

\def\basicsar[#1]{\vcase{\SAR}{#100}}%

\newcommand{\sar}{\@ifnextchar[{\basicsar}{\basicsar[50]}}%

\def\basicSar[#1]#2{\Vcase{\SAR}{#2}{#100}}%

\newcommand{\Sar}{\@ifnextchar[{\basicSar}{\basicSar[50]}}%

\def\basicsaR[#1]#2{\vcasE{\SAR}{#2}{#100}}%

\newcommand{\saR}{\@ifnextchar[{\basicsaR}{\basicsaR[50]}}%

\def\basicsdist[#1]{\vcase{\SDIST}{#100}}%

\newcommand{\sdist}{\@ifnextchar[{\basicsdist}{\basicsdist[50]}}%

\def\basicSdist[#1]#2{\Vcase{\SDIST}{#2\hspace*{2pt}}{#100}}%

\newcommand{\Sdist}{\@ifnextchar[{\basicSdist}{\basicSdist[50]}}%

\def\basicsdisT[#1]#2{\vcasE{\SDIST}{\hspace*{2pt}#2}{#100}}%

\newcommand{\sdisT}{\@ifnextchar[{\basicsdisT}{\basicsdisT[50]}}%

\def\basicsdotar[#1]{\vcase{\SDOTAR}{#100}}%

\newcommand{\sdotar}{\@ifnextchar[{\basicsdotar}{\basicsdotar[50]}}%

\def\basicSdotar[#1]#2{\Vcase{\SDOTAR}{#2}{#100}}%

\newcommand{\Sdotar}{\@ifnextchar[{\basicSdotar}{\basicSdotar[50]}}%

\def\basicsdotaR[#1]#2{\vcasE{\SDOTAR}{#2}{#100}}%

\newcommand{\sdotaR}{\@ifnextchar[{\basicsdotaR}{\basicsdotaR[50]}}%

\def\basicsmono[#1]{\vcase{\SMONO}{#100}}%

\newcommand{\smono}{\@ifnextchar[{\basicsmono}{\basicsmono[50]}}%

\def\basicSmono[#1]#2{\Vcase{\SMONO}{#2}{#100}}%

\newcommand{\Smono}{\@ifnextchar[{\basicSmono}{\basicSmono[50]}}%

\def\basicsmonO[#1]#2{\vcasE{\SMONO}{#2}{#100}}%

\newcommand{\smonO}{\@ifnextchar[{\basicsmonO}{\basicsmonO[50]}}%

\def\basicsepi[#1]{\vcase{\SEPI}{#100}}%

\newcommand{\sepi}{\@ifnextchar[{\basicsepi}{\basicsepi[50]}}%

\def\basicSepi[#1]#2{\Vcase{\SEPI}{#2}{#100}}%

\newcommand{\Sepi}{\@ifnextchar[{\basicSepi}{\basicSepi[50]}}%

\def\basicsepI[#1]#2{\vcasE{\SEPI}{#2}{#100}}%

\newcommand{\sepI}{\@ifnextchar[{\basicsepI}{\basicsepI[50]}}%

\def\basicsbimo[#1]{\vcase{\SBIMO}{#100}}%

\newcommand{\sbimo}{\@ifnextchar[{\basicsbimo}{\basicsbimo[50]}}%

\def\basicSbimo[#1]#2{\Vcase{\SBIMO}{#2}{#100}}%

\newcommand{\Sbimo}{\@ifnextchar[{\basicSbimo}{\basicSbimo[50]}}%

\def\basicsbimO[#1]#2{\vcasE{\SBIMO}{#2}{#100}}%

\newcommand{\sbimO}{\@ifnextchar[{\basicsbimO}{\basicsbimO[50]}}%

\def\basicsiso[#1]{\Visocase{\SAR}{\isosign}{}{#100}}%

\newcommand{\siso}{\@ifnextchar[{\basicsiso}{\basicsiso[50]}}%

\def\basicSiso[#1]#2{\Visocase{\SAR}{#2}{\isosign}{#100}}%

\newcommand{\Siso}{\@ifnextchar[{\basicSiso}{\basicSiso[50]}}%

\def\basicsisO[#1]#2{\Visocase{\SAR}{\isosign}{#2}{#100}}%

\newcommand{\sisO}{\@ifnextchar[{\basicsisO}{\basicsisO[50]}}%

\def\basicseql[#1]{\vcase{\SEQL}{#100}}%

\newcommand{\seql}{\@ifnextchar[{\basicseql}{\basicseql[50]}}%

\def\basicSeql[#1]#2{\Vcase{\SEQL}{#2\hspace*{2pt}}{#100}}%

\newcommand{\Seql}{\@ifnextchar[{\basicSeql}{\basicSeql[50]}}%

\def\basicseqL[#1]#2{\vcasE{\SEQL}{\hspace*{2pt}#2}{#100}}%

\newcommand{\seqL}{\@ifnextchar[{\basicseqL}{\basicseqL[50]}}%

\def\basicsline[#1]{\vcase{\SLINE}{#100}}%

\newcommand{\sline}{\@ifnextchar[{\basicsline}{\basicsline[50]}}%

\def\basicSline[#1]#2{\Vcase{\SLINE}{#2\hspace*{2pt}}{#100}}%

\newcommand{\Sline}{\@ifnextchar[{\basicSline}{\basicSline[50]}}%

\def\basicslinE[#1]#2{\vcasE{\SLINE}{\hspace*{2pt}#2}{#100}}%

\newcommand{\slinE}{\@ifnextchar[{\basicslinE}{\basicslinE[50]}}%

\def\basicsbiar[#1]{\vbicase{\SBIAR}{#100}}%

\newcommand{\sbiar}{\@ifnextchar[{\basicsbiar}{\basicsbiar[50]}}%

\def\basicSbiar[#1]#2#3{\Vbicase{\SBIAR}{#2}{#3}{#100}}%

\newcommand{\Sbiar}{\@ifnextchar[{\basicSbiar}{\basicSbiar[50]}}%

\def\basicsbidist[#1]{\vbicase{\SBIDIST}{#100}}%

\newcommand{\sbidist}{\@ifnextchar[{\basicsbidist}{\basicsbidist[50]}}%

\def\basicSbidist[#1]#2#3%
{\Vbicase{\SBIDIST}{#2\hspace*{2pt}}{\hspace*{2pt}#3}{#100}}%

\newcommand{\Sbidist}{\@ifnextchar[{\basicSbidist}{\basicSbidist[50]}}%

\def\basicsadjar[#1]{\vbicase{\SADJAR}{#100}}%

\newcommand{\sadjar}{\@ifnextchar[{\basicsadjar}{\basicsadjar[50]}}%

\def\basicSadjar[#1]#2#3{\Vbicase{\SADJAR}{#2}{#3}{#100}}%

\newcommand{\Sadjar}{\@ifnextchar[{\basicSadjar}{\basicSadjar[50]}}%

\def\basicsadjdist[#1]{\vbicase{\SADJDIST}{#100}}%

\newcommand{\sadjdist}{\@ifnextchar[{\basicsadjdist}{\basicsadjdist[50]}}%

\def\basicSadjdist[#1]#2#3%
{\Vbicase{\SADJDIST}{#2\hspace*{2pt}}{\hspace*{2pt}#3}{#100}}%

\newcommand{\Sadjdist}{\@ifnextchar[{\basicSadjdist}{\basicSadjdist[50]}}%

\def\basicstriar[#1]{\vtricase{\STRIAR}{#100}}%

\newcommand{\striar}{\@ifnextchar[{\basicstriar}{\basicstriar[50]}}%

\def\basicStriar[#1]#2#3#4{\Vtricase{\STRIAR}{#2}{#3}{#4}{}{#100}}%

\newcommand{\Striar}{\@ifnextchar[{\basicStriar}{\basicStriar[50]}}%

\def\basicstriaR[#1]#2#3#4{\Vtricase{\STRIAR}{}{#2}{#3}{#4}{#100}}%

\newcommand{\striaR}{\@ifnextchar[{\basicstriaR}{\basicstriaR[50]}}%

\def\basicstriadjar[#1]{\vtricase{\STRIADJAR}{#100}}%

\newcommand{\striadjar}{\@ifnextchar[{\basicstriadjar}{\basicstriadjar[50]}}%

\def\basicStriadjar[#1]#2#3#4{\Vtricase{\STRIADJAR}{#2}{#3}{#4}{}{#100}}%

\newcommand{\Striadjar}{\@ifnextchar[{\basicStriadjar}{\basicStriadjar[50]}}%

\def\basicstriadjaR[#1]#2#3#4{\Vtricase{\STRIADJAR}{}{#2}{#3}{#4}{#100}}%

\newcommand{\striadjaR}{\@ifnextchar[{\basicstriadjaR}{\basicstriadjaR[50]}}%

\def\basicsquadriar[#1]{%
\Z=#1%
\multiply \Z by 100%
\vmulticase{\SQUADRIAR}{\Z}}%

\newcommand{\squadriar}{\@ifnextchar[{\basicsquadriar}%
{\hspace{\SOURCE\unitlength}\basicsquadriar[\ARROWLENGTH]}}%

\def\basicsquadriadjar[#1]{%
\Z=#1%
\multiply \Z by 100%
\vmulticase{\SQUADRIADJAR}{\Z}}%

\newcommand{\squadriadjar}{\@ifnextchar[{\basicsquadriadjar}%
{\hspace{\SOURCE\unitlength}\basicsquadriadjar[\ARROWLENGTH]}}%

\def\basicsquintiar[#1]{%
\Z=#1%
\multiply \Z by 100%
\vmulticase{\SQUINTIAR}{\Z}}%

\newcommand{\squintiar}{\@ifnextchar[{\basicsquintiar}%
{\hspace{\SOURCE\unitlength}\basicsquintiar[\ARROWLENGTH]}}%

\def\basicsquintiadjar[#1]{%
\Z=#1%
\multiply \Z by 100%
\vmulticase{\SQUINTIADJAR}{\Z}}%

\newcommand{\squintiadjar}{\@ifnextchar[{\basicsquintiadjar}%
{\hspace{\SOURCE\unitlength}\basicsquintiadjar[\ARROWLENGTH]}}%

\newcommand{\NAR}[1]%
{\begin{picture}(0,0)%
\put(0,0){\makebox(0,0)%
{\begin{picture}(0,#1)%
\put(0,0){\line(0,1){#1}}%
\put(0,#1){\nhead}%
\end{picture}}}\end{picture}}%

\newcommand{\NDIST}[1]%
{\begin{picture}(0,0)%
\put(0,0){\makebox(0,0)%
{\begin{picture}(0,#1)%
\put(0,0){\line(0,1){#1}}%
\put(0,#1){\nhead}%
\end{picture}}}
\put(0,0){\distsign}%
\end{picture}}%

\newcommand{\NDOTAR}[1]%
{\truex{100}\truey{300}%
\NUMBEROFDOTS=#1%
\divide\NUMBEROFDOTS by \value{y}%
\advance\NUMBEROFDOTS by 1%
\begin{picture}(0,0)%
\put(0,0){\makebox(0,0)%
{\begin{picture}(0,#1)%
\multiput(0,0)(0,\value{y}){\NUMBEROFDOTS}%
{\circle*{\value{x}}}%
\put(0,#1){\nhead}%
\end{picture}}}\end{picture}}%

\newcommand{\NMONO}[1]%
{\truetail%
\monolength=#1%
\advance\monolength by -\truemonotail%
\begin{picture}(0,0)%
\put(0,0){\makebox(0,0)%
{\begin{picture}(0,#1)%
\put(0,\truemonotail){\line(0,1){\monolength}}%
\put(0,#1){\nhead}%
\put(0,\truemonotail){\nhead}%
\end{picture}}}\end{picture}}%

\newcommand{\NEPI}[1]%
{\truehead%
\epilength=#1%
\advance\epilength by -\trueepihead%
\begin{picture}(0,0)%
\put(0,0){\makebox(0,0)%
{\begin{picture}(0,#1)%
\put(0,0){\line(0,1){\epilength}}%
\put(0,#1){\nhead}%
\put(0,\epilength){\nhead}%
\end{picture}}}\end{picture}}%

\newcommand{\NBIMO}[1]%
{\truehead\truetail%
\epilength=#1%
\advance\epilength by -\trueepihead%
\monolength=\epilength%
\advance\monolength by -\truemonotail%
\begin{picture}(0,0)%
\put(0,0){\makebox(0,0)%
{\begin{picture}(0,#1)%
\put(0,\truemonotail){\line(0,1){\monolength}}%
\put(0,#1){\nhead}%
\put(0,\truemonotail){\nhead}%
\put(0,\epilength){\nhead}%
\end{picture}}}\end{picture}}%

\newcommand{\NBIAR}[1]%
{\begin{picture}(0,0)%
\truex{350}%
\put(0,0){\makebox(0,0)%
{\begin{picture}(0,#1)%
\put(-\value{x},0){\line(0,1){#1}}%
\put(-\value{x},#1){\nhead}%
\put(\value{x},0){\line(0,1){#1}}%
\put(\value{x},#1){\nhead}%
\end{picture}}}\end{picture}}%

\newcommand{\NBIDIST}[1]%
{\begin{picture}(0,0)%
\truex{350}%
\put(0,0){\makebox(0,0)%
{\begin{picture}(0,#1)%
\put(-\value{x},0){\line(0,1){#1}}%
\put(-\value{x},#1){\nhead}%
\put(\value{x},0){\line(0,1){#1}}%
\put(\value{x},#1){\nhead}%
\end{picture}}}
\put(-\value{x},0){\distsign}%
\put(\value{x},0){\distsign}%
\end{picture}}%

\newcommand{\NADJAR}[1]{\begin{picture}(0,0)%
\truex{350}%
\put(0,0){\makebox(0,0)%
{\begin{picture}(0,#1)%
\put(\value{x},#1){\line(0,-1){#1}}%
\put(\value{x},0){\shead}%
\put(-\value{x},0){\line(0,1){#1}}%
\put(-\value{x},#1){\nhead}%
\end{picture}}}\end{picture}}%

\newcommand{\NADJDIST}[1]{\begin{picture}(0,0)%
\truex{350}%
\put(0,0){\makebox(0,0)%
{\begin{picture}(0,#1)%
\put(\value{x},#1){\line(0,-1){#1}}%
\put(\value{x},0){\shead}%
\put(-\value{x},0){\line(0,1){#1}}%
\put(-\value{x},#1){\nhead}%
\end{picture}}}
\put(-\value{x},0){\distsign}%
\put(\value{x},0){\distsign}%
\end{picture}}%

\newcommand{\NTRIAR}[5]%
{\truex{1200}%
\truey{1800}%
\truez{600}%
\begin{picture}(0,0)%
\put(0,0){\makebox(0,0)%
{\begin{picture}(0,#5)%
\X=#5\divide\X by 2%
\put(-\value{x},0){\line(0,1){#5}}%
\put(-\value{x},#5){\nhead}%
\put(0,0){\line(0,1){#5}}%
\put(0,#5){\nhead}%
\put(\value{x},0){\line(0,1){#5}}%
\put(\value{x},#5){\nhead}%
\put(-\value{y},\X){\makebox(0,0){$#1$}}%
\put(\value{y},\X){\makebox(0,0){$#4$}}%
\put(-\value{z},\X){\makebox(0,0){$#2$}}%
\put(\value{z},\X){\makebox(0,0){$#3$}}%
\end{picture}}}\end{picture}}%

\newcommand{\NTRIADJAR}[5]%
{\truex{1200}%
\truey{1800}%
\truez{600}%
\begin{picture}(0,0)%
\put(0,0){\makebox(0,0)%
{\begin{picture}(0,#5)%
\X=#5\divide\X by 2%
\put(-\value{x},#5){\line(0,-1){#5}}%
\put(-\value{x},0){\shead}%
\put(0,0){\line(0,1){#5}}%
\put(0,#5){\nhead}%
\put(\value{x},#5){\line(0,-1){#5}}%
\put(\value{x},0){\shead}%
\put(-\value{y},\X){\makebox(0,0){$#1$}}%
\put(\value{y},\X){\makebox(0,0){$#4$}}%
\put(-\value{z},\X){\makebox(0,0){$#2$}}%
\put(\value{z},\X){\makebox(0,0){$#3$}}%
\end{picture}}}\end{picture}}%

\newcommand{\NQUADRIAR}[1]%
{\truex{400}%
\truey{1200}%
\begin{picture}(0,0)%
\put(0,0){\makebox(0,0)%
{\begin{picture}(0,#1)%
\put(-\value{y},0){\line(0,1){#1}}%
\put(-\value{y},#1){\nhead}%
\put(-\value{x},0){\line(0,1){#1}}%
\put(-\value{x},#1){\nhead}%
\put(\value{x},0){\line(0,1){#1}}%
\put(\value{x},#1){\nhead}%
\put(\value{y},0){\line(0,1){#1}}%
\put(\value{y},#1){\nhead}%
\end{picture}}}\end{picture}}%

\newcommand{\NQUADRIADJAR}[1]%
{\truex{400}%
\truey{1200}%
\begin{picture}(0,0)%
\put(0,0){\makebox(0,0)%
{\begin{picture}(0,#1)%
\put(-\value{y},0){\line(0,1){#1}}%
\put(-\value{y},#1){\nhead}%
\put(-\value{x},#1){\line(0,-1){#1}}%
\put(-\value{x},0){\shead}%
\put(\value{x},0){\line(0,1){#1}}%
\put(\value{x},#1){\nhead}%
\put(\value{y},#1){\line(0,-1){#1}}%
\put(\value{y},0){\shead}%
\end{picture}}}\end{picture}}%

\newcommand{\NQUINTIAR}[1]%
{\truex{600}%
\truey{1200}%
\begin{picture}(0,0)%
\put(0,0){\makebox(0,0)%
{\begin{picture}(0,#1)%
\put(-\value{y},0){\line(0,1){#1}}%
\put(-\value{y},#1){\nhead}%
\put(-\value{x},0){\line(0,1){#1}}%
\put(-\value{x},#1){\nhead}%
\put(0,0){\line(0,1){#1}}%
\put(0,#1){\nhead}%
\put(\value{x},0){\line(0,1){#1}}%
\put(\value{x},#1){\nhead}%
\put(\value{y},0){\line(0,1){#1}}%
\put(\value{y},#1){\nhead}%
\end{picture}}}\end{picture}}%

\newcommand{\NQUINTIADJAR}[1]%
{\truex{600}%
\truey{1200}%
\begin{picture}(0,0)%
\put(0,0){\makebox(0,0)%
{\begin{picture}(0,#1)%
\put(-\value{y},0){\line(0,1){#1}}%
\put(-\value{y},#1){\nhead}%
\put(-\value{x},#1){\line(0,-1){#1}}%
\put(-\value{x},0){\shead}%
\put(0,0){\line(0,1){#1}}%
\put(0,#1){\nhead}%
\put(\value{x},#1){\line(0,-1){#1}}%
\put(\value{x},0){\shead}%
\put(\value{y},0){\line(0,1){#1}}%
\put(\value{y},#1){\nhead}%
\end{picture}}}\end{picture}}%

\def\basicnar[#1]{\vcase{\NAR}{#100}}%

\newcommand{\nar}{\@ifnextchar[{\basicnar}{\basicnar[50]}}%

\def\basicNar[#1]#2{\Vcase{\NAR}{#2}{#100}}%

\newcommand{\Nar}{\@ifnextchar[{\basicNar}{\basicNar[50]}}%

\def\basicnaR[#1]#2{\vcasE{\NAR}{#2}{#100}}%

\newcommand{\naR}{\@ifnextchar[{\basicnaR}{\basicnaR[50]}}%

\def\basicndist[#1]{\vcase{\NDIST}{#100}}%

\newcommand{\ndist}{\@ifnextchar[{\basicndist}{\basicndist[50]}}%

\def\basicNdist[#1]#2{\Vcase{\NDIST}{#2\hspace*{2pt}}{#100}}%

\newcommand{\Ndist}{\@ifnextchar[{\basicNdist}{\basicNdist[50]}}%

\def\basicndisT[#1]#2{\vcasE{\NDIST}{\hspace*{2pt}#2}{#100}}%

\newcommand{\ndisT}{\@ifnextchar[{\basicndisT}{\basicndisT[50]}}%

\def\basicndotar[#1]{\vcase{\NDOTAR}{#100}}%

\newcommand{\ndotar}{\@ifnextchar[{\basicndotar}{\basicndotar[50]}}%

\def\basicNdotar[#1]#2{\Vcase{\NDOTAR}{#2}{#100}}%

\newcommand{\Ndotar}{\@ifnextchar[{\basicNdotar}{\basicNdotar[50]}}%

\def\basicndotaR[#1]#2{\vcasE{\NDOTAR}{#2}{#100}}%

\newcommand{\ndotaR}{\@ifnextchar[{\basicndotaR}{\basicndotaR[50]}}%

\def\basicnmono[#1]{\vcase{\NMONO}{#100}}%

\newcommand{\nmono}{\@ifnextchar[{\basicnmono}%
{\basicnmono[50]}}%

\def\basicNmono[#1]#2{\Vcase{\NMONO}{#2}{#100}}%

\newcommand{\Nmono}{\@ifnextchar[{\basicNmono}{\basicNmono[50]}}%

\def\basicnmonO[#1]#2{\vcasE{\NMONO}{#2}{#100}}%

\newcommand{\nmonO}{\@ifnextchar[{\basicnmonO}{\basicnmonO[50]}}%

\def\basicnepi[#1]{\vcase{\NEPI}{#100}}%

\newcommand{\nepi}{\@ifnextchar[{\basicnepi}{\basicnepi[50]}}%

\def\basicNepi[#1]#2{\Vcase{\NEPI}{#2}{#100}}%

\newcommand{\Nepi}{\@ifnextchar[{\basicNepi}{\basicNepi[50]}}%

\def\basicnepI[#1]#2{\vcasE{\NEPI}{#2}{#100}}%

\newcommand{\nepI}{\@ifnextchar[{\basicnepI}{\basicnepI[50]}}%

\def\basicnbimo[#1]{\vcase{\NBIMO}{#100}}%

\newcommand{\nbimo}{\@ifnextchar[{\basicnbimo}{\basicnbimo[50]}}%

\def\basicNbimo[#1]#2{\Vcase{\NBIMO}{#2}{#100}}%

\newcommand{\Nbimo}{\@ifnextchar[{\basicNbimo}{\basicNbimo[50]}}%

\def\basicnbimO[#1]#2{\vcasE{\NBIMO}{#2}{#100}}%

\newcommand{\nbimO}{\@ifnextchar[{\basicnbimO}{\basicnbimO[50]}}%

\def\basicniso[#1]{\Visocase{\NAR}{\isosign}{}{#100}}%

\newcommand{\niso}{\@ifnextchar[{\basicniso}{\basicniso[50]}}%

\def\basicNiso[#1]#2{\Visocase{\NAR}{#2}{\isosign}{#100}}%

\newcommand{\Niso}{\@ifnextchar[{\basicNiso}{\basicNiso[50]}}%

\def\basicnisO[#1]#2{\Visocase{\NAR}{\isosign}{#2}{#100}}%

\newcommand{\nisO}{\@ifnextchar[{\basicnisO}{\basicnisO[50]}}%

\def\basicnbiar[#1]{\vbicase{\NBIAR}{#100}}%

\newcommand{\nbiar}{\@ifnextchar[{\basicnbiar}{\basicnbiar[50]}}%

\def\basicNbiar[#1]#2#3{\Vbicase{\NBIAR}{#2}{#3}{#100}}%

\newcommand{\Nbiar}{\@ifnextchar[{\basicNbiar}{\basicNbiar[50]}}%

\def\basicnbidist[#1]{\vbicase{\NBIDIST}{#100}}%

\newcommand{\nbidist}{\@ifnextchar[{\basicnbidist}{\basicnbidist[50]}}%

\def\basicNbidist[#1]#2#3%
{\Vbicase{\NBIDIST}{#2\hspace*{2pt}}{\hspace*{2pt}#3}{#100}}%

\newcommand{\Nbidist}{\@ifnextchar[{\basicNbidist}{\basicNbidist[50]}}%

\def\basicnadjar[#1]{\vbicase{\NADJAR}{#100}}%

\newcommand{\nadjar}{\@ifnextchar[{\basicnadjar}{\basicnadjar[50]}}%

\def\basicNadjar[#1]#2#3{\Vbicase{\NADJAR}{#2}{#3}{#100}}%

\newcommand{\Nadjar}{\@ifnextchar[{\basicNadjar}{\basicNadjar[50]}}%

\def\basicnadjdist[#1]{\vbicase{\NADJDIST}{#100}}%

\newcommand{\nadjdist}{\@ifnextchar[{\basicnadjdist}{\basicnadjdist[50]}}%

\def\basicNadjdist[#1]#2#3%
{\Vbicase{\NADJDIST}{#2\hspace*{2pt}}{\hspace*{2pt}#3}{#100}}%

\newcommand{\Nadjdist}{\@ifnextchar[{\basicNadjdist}{\basicNadjdist[50]}}%

\def\basicntriar[#1]{\vtricase{\NTRIAR}{#100}}%

\newcommand{\ntriar}{\@ifnextchar[{\basicntriar}{\basicntriar[50]}}%

\def\basicNtriar[#1]#2#3#4{\Vtricase{\NTRIAR}{#2}{#3}{#4}{}{#100}}%

\newcommand{\Ntriar}{\@ifnextchar[{\basicNtriar}{\basicNtriar[50]}}%

\def\basicntriaR[#1]#2#3#4{\Vtricase{\NTRIAR}{}{#2}{#3}{#4}{#100}}%

\newcommand{\ntriaR}{\@ifnextchar[{\basicntriaR}{\basicntriaR[50]}}%

\def\basicntriadjar[#1]{\vtricase{\NTRIADJAR}{#100}}%

\newcommand{\ntriadjar}{\@ifnextchar[{\basicntriadjar}{\basicntriadjar[50]}}%

\def\basicNtriadjar[#1]#2#3#4{\Vtricase{\NTRIADJAR}{#2}{#3}{#4}{}{#100}}%

\newcommand{\Ntriadjar}{\@ifnextchar[{\basicNtriadjar}{\basicNtriadjar[50]}}%

\def\basicntriadjaR[#1]#2#3#4{\Vtricase{\NTRIADJAR}{}{#2}{#3}{#4}{#100}}%

\newcommand{\ntriadjaR}{\@ifnextchar[{\basicntriadjaR}{\basicntriadjaR[50]}}%

\def\basicnquadriar[#1]{%
\Z=#1%
\multiply \Z by 100%
\vmulticase{\NQUADRIAR}{\Z}}%

\newcommand{\nquadriar}{\@ifnextchar[{\basicnquadriar}%
{\hspace{\SOURCE\unitlength}\basicnquadriar[\ARROWLENGTH]}}%

\def\basicnquadriadjar[#1]{%
\Z=#1%
\multiply \Z by 100%
\vmulticase{\NQUADRIADJAR}{\Z}}%

\newcommand{\nquadriadjar}{\@ifnextchar[{\basicnquadriadjar}%
{\hspace{\SOURCE\unitlength}\basicnquadriadjar[\ARROWLENGTH]}}%

\def\basicnquintiar[#1]{%
\Z=#1%
\multiply \Z by 100%
\vmulticase{\NQUINTIAR}{\Z}}%

\newcommand{\nquintiar}{\@ifnextchar[{\basicnquintiar}%
{\hspace{\SOURCE\unitlength}\basicnquintiar[\ARROWLENGTH]}}%

\def\basicnquintiadjar[#1]{%
\Z=#1%
\multiply \Z by 100%
\vmulticase{\NQUINTIADJAR}{\Z}}%

\newcommand{\nquintiadjar}{\@ifnextchar[{\basicnquintiadjar}%
{\hspace{\SOURCE\unitlength}\basicnquintiadjar[\ARROWLENGTH]}}%

\newcommand{\fdcase}[4]{\begin{picture}(0,0)%
\put(0,0){#1{#4}}%
\truex{200}\truey{600}\truez{600}%
\put(-\value{x},-\value{x}){\makebox(0,\value{z})[r]{${#2}$}}%
\put(\value{x},-\value{y}){\makebox(0,\value{z})[l]{${#3}$}}%
\end{picture}}%

\newcommand{\fdbicase}[4]{\begin{picture}(0,0)%
\put(0,0){#1{#4}}%
\truex{600}\truey{150}\truez{600}%
\put(-\value{x},\value{y}){\makebox(0,\value{z})[r]{${#2}$}}%
\truex{300}\truey{1050}\truez{600}%
\put(\value{x},-\value{y}){\makebox(0,\value{z})[l]{${#3}$}}%
\end{picture}}%

\newcommand{\NEAR}[1]{%
\Y=#1%
\divide\Y by 2%
\begin{picture}(0,0)%
\put(-\Y,-\Y){\line(1,1){#1}}%
\put(\Y,\Y){\nehead}%
\end{picture}}%

\newcommand{\NEDIST}[1]{%
\Y=#1%
\divide\Y by 2%
\begin{picture}(0,0)%
\put(-\Y,-\Y){\line(1,1){#1}}%
\put(\Y,\Y){\nehead}%
\put(0,0){\distsign}%
\end{picture}}%

\newcommand{\NEDOTAR}[1]%
{\truex{100}\truey{212}%
\Y=#1%
\divide\Y by 2%
\NUMBEROFDOTS=#1%
\divide\NUMBEROFDOTS by \value{y}%
\advance\NUMBEROFDOTS by 1%
\begin{picture}(0,0)%
\multiput(-\Y,-\Y)(\value{y},\value{y}){\NUMBEROFDOTS}%
{\circle*{\value{x}}}%
\put(\Y,\Y){\nehead}%
\end{picture}}%

\newcommand{\NEMONO}[1]{%
\Y=#1%
\divide \Y by 2%
\Truetail%
\bimolength=#1%
\advance\bimolength by -\Truemonotail%
\monolength=\bimolength%
\advance\monolength by -\Y%
\begin{picture}(0,0)%
\put(-\monolength,-\monolength){\line(1,1){\bimolength}}%
\put(-\monolength,-\monolength){\nehead}%
\put(\Y,\Y){\nehead}%
\end{picture}}%

\newcommand{\NEEPI}[1]{%
\Y=#1%
\divide\Y by 2%
\Truehead%
\bimolength=#1%
\advance\bimolength by -\Trueepihead%
\epilength=\bimolength%
\advance\epilength by -\Y%
\begin{picture}(0,0)%
\put(-\Y,-\Y){\line(1,1){\bimolength}}%
\put(\epilength,\epilength){\nehead}%
\put(\Y,\Y){\nehead}%
\end{picture}}%

\newcommand{\NEBIMO}[1]{%
\Y=#1%
\divide\Y by 2%
\Truetail\Truehead%
\bimolength=#1%
\advance\bimolength by -\Truemonotail%
\monolength=\bimolength%
\advance\monolength by -\Y%
\advance\bimolength by -\Trueepihead%
\epilength=\bimolength%
\advance\epilength by -\monolength%
\begin{picture}(0,0)%
\put(-\monolength,-\monolength){\line(1,1){\bimolength}}%
\put(-\monolength,-\monolength){\nehead}%
\put(\epilength,\epilength){\nehead}%
\put(\Y,\Y){\nehead}%
\end{picture}}%

\newcommand{\NEBIAR}[1]{%
\Y=#1%
\divide\Y by 2%
\begin{picture}(0,0)%
\put(-\Y,-\Y){\begin{picture}(0,0)%
\truex{247}%
\put(-\value{x},\value{x}){\line(1,1){#1}}%
\put(\value{x},-\value{x}){\line(1,1){#1}}%
\monolength=#1%
\advance\monolength by -\value{x}%
\epilength=#1%
\advance\epilength by \value{x}%
\put(\monolength,\epilength){\nehead}%
\put(\epilength,\monolength){\nehead}%
\end{picture}}\end{picture}}%

\newcommand{\NEBIDIST}[1]{%
\Y=#1%
\divide\Y by 2%
\begin{picture}(0,0)%
\put(-\Y,-\Y){\begin{picture}(0,0)%
\truex{247}%
\monolength=#1%
\advance\monolength by -\value{x}%
\epilength=#1%
\advance\epilength by \value{x}%
\put(\value{x},-\value{x}){\line(1,1){#1}}%
\put(\epilength,\monolength){\nehead}%
\end{picture}}%
\put(-\Y,-\Y){\begin{picture}(0,0)%
\truex{247}%
\monolength=#1%
\advance\monolength by \value{x}%
\epilength=#1%
\advance\epilength by -\value{x}%
\put(-\value{x},\value{x}){\line(1,1){#1}}%
\put(\epilength,\monolength){\nehead}%
\end{picture}}%
\put(-\value{x},\value{x}){\distsign}%
\put(\value{x},-\value{x}){\distsign}%
\end{picture}}%

\newcommand{\NEEQL}[1]{%
\Y=#1%
\divide\Y by 2%
\begin{picture}(0,0)%
\put(-\Y,-\Y){\begin{picture}(0,0)%
\truex{70}%
\put(-\value{x},\value{x}){\line(1,1){#1}}%
\put(\value{x},-\value{x}){\line(1,1){#1}}%
\end{picture}}\end{picture}}%

\newcommand{\NELINE}[1]{%
\Y=#1%
\divide\Y by 2%
\begin{picture}(0,0)%
\put(-\Y,-\Y){\begin{picture}(0,0)%
\put(0,0){\line(1,1){#1}}%
\end{picture}}\end{picture}}%

\newcommand{\NEADJAR}[1]{%
\Y=#1%
\divide\Y by 2%
\begin{picture}(0,0)%
\put(-\Y,-\Y){\begin{picture}(0,0)%
\truex{247}%
\monolength=#1%
\advance\monolength by -\value{x}%
\epilength=#1%
\advance\epilength by \value{x}%
\put(\value{x},-\value{x}){\line(1,1){#1}}%
\put(\epilength,\monolength){\nehead}%
\end{picture}}%
\put(\Y,\Y){\begin{picture}(0,0)%
\truex{247}%
\monolength=#1%
\advance\monolength by -\value{x}%
\epilength=#1%
\advance\epilength by \value{x}%
\put(-\value{x},\value{x}){\line(-1,-1){#1}}%
\put(-\epilength,-\monolength){\swhead}%
\end{picture}}\end{picture}}%

\newcommand{\NEADJDIST}[1]{%
\Y=#1%
\divide\Y by 2%
\begin{picture}(0,0)%
\put(-\Y,-\Y){\begin{picture}(0,0)%
\truex{247}%
\monolength=#1%
\advance\monolength by -\value{x}%
\epilength=#1%
\advance\epilength by \value{x}%
\put(\value{x},-\value{x}){\line(1,1){#1}}%
\put(\epilength,\monolength){\nehead}%
\end{picture}}%
\put(\Y,\Y){\begin{picture}(0,0)%
\truex{247}%
\monolength=#1%
\advance\monolength by -\value{x}%
\epilength=#1%
\advance\epilength by \value{x}%
\put(-\value{x},\value{x}){\line(-1,-1){#1}}%
\put(-\epilength,-\monolength){\swhead}%
\end{picture}}%
\put(-\value{x},\value{x}){\distsign}%
\put(\value{x},-\value{x}){\distsign}%
\end{picture}}%

\def\basicnear[#1]{\fdcase{\NEAR}{}{}{#100}}%

\newcommand{\near}{\@ifnextchar[{\basicnear}{\basicnear[59]}}%

\def\basicNear[#1]#2{\fdcase{\NEAR}{#2}{}{#100}}%

\newcommand{\Near}{\@ifnextchar[{\basicNear}{\basicNear[59]}}%

\def\basicneaR[#1]#2{\fdcase{\NEAR}{}{#2}{#100}}%

\newcommand{\neaR}{\@ifnextchar[{\basicneaR}{\basicneaR[59]}}%

\def\basicnedist[#1]{\fdcase{\NEDIST}{}{}{#100}}%

\newcommand{\nedist}{\@ifnextchar[{\basicnedist}{\basicnedist[59]}}%

\def\basicNedist[#1]#2{\fdcase{\NEDIST}{#2}{}{#100}}%

\newcommand{\Nedist}{\@ifnextchar[{\basicNedist}{\basicNedist[59]}}%

\def\basicnedisT[#1]#2{\fdcase{\NEDIST}{}{#2}{#100}}%

\newcommand{\nedisT}{\@ifnextchar[{\basicnedisT}{\basicnedisT[59]}}%

\def\basicnedotar[#1]{\fdcase{\NEDOTAR}{}{}{#100}}%

\newcommand{\nedotar}{\@ifnextchar[{\basicnedotar}{\basicnedotar[59]}}%

\def\basicNedotar[#1]#2{\fdcase{\NEDOTAR}{#2}{}{#100}}%

\newcommand{\Nedotar}{\@ifnextchar[{\basicNedotar}{\basicNedotar[59]}}%

\def\basicnedotaR[#1]#2{\fdcase{\NEDOTAR}{}{#2}{#100}}%

\newcommand{\nedotaR}{\@ifnextchar[{\basicnedotaR}{\basicnedotaR[59]}}%

\def\basicnemono[#1]{\fdcase{\NEMONO}{}{}{#100}}%

\newcommand{\nemono}{\@ifnextchar[{\basicnemono}{\basicnemono[59]}}%

\def\basicNemono[#1]#2{\fdcase{\NEMONO}{#2}{}{#100}}%

\newcommand{\Nemono}{\@ifnextchar[{\basicNemono}{\basicNemono[59]}}%

\def\basicnemonO[#1]#2{\fdcase{\NEMONO}{}{#2}{#100}}%

\newcommand{\nemonO}{\@ifnextchar[{\basicnemonO}{\basicnemonO[59]}}%

\def\basicneepi[#1]{\fdcase{\NEEPI}{}{}{#100}}%

\newcommand{\neepi}{\@ifnextchar[{\basicneepi}{\basicneepi[59]}}%

\def\basicNeepi[#1]#2{\fdcase{\NEEPI}{#2}{}{#100}}%

\newcommand{\Neepi}{\@ifnextchar[{\basicNeepi}{\basicNeepi[59]}}%

\def\basicneepI[#1]#2{\fdcase{\NEEPI}{}{#2}{#100}}%

\newcommand{\neepI}{\@ifnextchar[{\basicneepI}{\basicneepI[59]}}%

\def\basicnebimo[#1]{\fdcase{\NEBIMO}{}{}{#100}}%

\newcommand{\nebimo}{\@ifnextchar[{\basicnebimo}{\basicnebimo[59]}}%

\def\basicNebimo[#1]#2{\fdcase{\NEBIMO}{#2}{}{#100}}%

\newcommand{\Nebimo}{\@ifnextchar[{\basicNebimo}{\basicNebimo[59]}}%

\def\basicnebimO[#1]#2{\fdcase{\NEBIMO}{}{#2}{#100}}%

\newcommand{\nebimO}{\@ifnextchar[{\basicnebimO}{\basicnebimO[59]}}%

\def\basicneiso[#1]{\fdcase{\NEAR}{\hspace{-2pt}\isosign}{}{#100}}%

\newcommand{\neiso}{\@ifnextchar[{\basicneiso}{\basicneiso[59]}}%

\def\basicNeiso[#1]#2{\fdcase{\NEAR}{#2}{\isosign}{#100}}%

\newcommand{\Neiso}{\@ifnextchar[{\basicNeiso}{\basicNeiso[59]}}%

\def\basicneisO[#1]#2{\fdcase{\NEAR}{\hspace{-2pt}\isosign}{#2}{#100}}%

\newcommand{\neisO}{\@ifnextchar[{\basicneisO}{\basicneisO[59]}}%

\def\basicneeql[#1]{\fdcase{\NEEQL}{}{}{#100}}%

\newcommand{\neeql}{\@ifnextchar[{\basicneeql}{\basicneeql[59]}}%

\def\basicNeeql[#1]#2{\fdcase{\NEEQL}{#2}{}{#100}}%

\newcommand{\Neeql}{\@ifnextchar[{\basicNeeql}{\basicNeeql[59]}}%

\def\basicneeqL[#1]#2{\fdcase{\NEEQL}{}{#2}{#100}}%

\newcommand{\neeqL}{\@ifnextchar[{\basicneeqL}{\basicneeqL[59]}}%

\def\basicneline[#1]{\fdcase{\NELINE}{}{}{#100}}%

\newcommand{\neline}{\@ifnextchar[{\basicneline}{\basicneline[59]}}%

\def\basicNeline[#1]#2{\fdcase{\NELINE}{#2}{}{#100}}%

\newcommand{\Neline}{\@ifnextchar[{\basicNeline}{\basicNeline[59]}}%

\def\basicnelinE[#1]#2{\fdcase{\NELINE}{}{#2}{#100}}%

\newcommand{\nelinE}{\@ifnextchar[{\basicnelinE}{\basicneeqL[59]}}%

\def\basicnebiar[#1]{\fdbicase{\NEBIAR}{}{}{#100}}%

\newcommand{\nebiar}{\@ifnextchar[{\basicnebiar}{\basicnebiar[59]}}%

\def\basicNebiar[#1]#2#3{\fdbicase{\NEBIAR}{#2}{#3}{#100}}%

\newcommand{\Nebiar}{\@ifnextchar[{\basicNebiar}{\basicNebiar[59]}}%

\def\basicneadjar[#1]{\fdbicase{\NEADJAR}{}{}{#100}}%

\newcommand{\neadjar}{\@ifnextchar[{\basicneadjar}{\basicneadjar[59]}}%

\def\basicNeadjar[#1]#2#3{\fdbicase{\NEADJAR}{#2}{#3}{#100}}%

\newcommand{\Neadjar}{\@ifnextchar[{\basicNeadjar}{\basicNeadjar[59]}}%

\def\basicnebidist[#1]{\fdbicase{\NEBIDIST}{}{}{#100}}%

\newcommand{\nebidist}{\@ifnextchar[{\basicnebidist}{\basicnebidist[59]}}%

\def\basicNebidist[#1]#2#3{\fdbicase{\NEBIDIST}{#2}{#3}{#100}}%

\newcommand{\Nebidist}{\@ifnextchar[{\basicNebidist}{\basicNebidist[59]}}%

\def\basicneadjdist[#1]{\fdbicase{\NEADJDIST}{}{}{#100}}%

\newcommand{\neadjdist}{\@ifnextchar[{\basicneadjdist}{\basicneadjdist[59]}}%

\def\basicNeadjdist[#1]#2#3{\fdbicase{\NEADJDIST}{#2}{#3}{#100}}%

\newcommand{\Neadjdist}{\@ifnextchar[{\basicNeadjdist}{\basicNeadjdist[59]}}%

\newcommand{\SWAR}[1]{%
\Y=#1%
\divide\Y by 2%
\begin{picture}(0,0)%
\put(\Y,\Y){\line(-1,-1){#1}}%
\put(-\Y,-\Y){\swhead}%
\end{picture}}%

\newcommand{\SWDIST}[1]{%
\Y=#1%
\divide\Y by 2%
\begin{picture}(0,0)%
\put(\Y,\Y){\line(-1,-1){#1}}%
\put(-\Y,-\Y){\swhead}%
\put(0,0){\distsign}%
\end{picture}}%

\newcommand{\SWDOTAR}[1]%
{\truex{100}\truey{212}%
\Y=#1%
\divide\Y by 2%
\NUMBEROFDOTS=#1%
\divide\NUMBEROFDOTS by \value{y}%
\advance\NUMBEROFDOTS by 1%
\begin{picture}(0,0)%
\multiput(\Y,\Y)(-\value{y},-\value{y}){\NUMBEROFDOTS}%
{\circle*{\value{x}}}%
\put(-\Y,-\Y){\swhead}%
\end{picture}}%

\newcommand{\SWMONO}[1]{%
\Y=#1%
\divide \Y by 2%
\Truetail%
\bimolength=#1%
\advance\bimolength by -\Truemonotail%
\monolength=\bimolength%
\advance\monolength by -\Y%
\begin{picture}(0,0)%
\put(\monolength,\monolength){\line(-1,-1){\bimolength}}%
\put(\monolength,\monolength){\swhead}%
\put(-\Y,-\Y){\swhead}%
\end{picture}}%

\newcommand{\SWEPI}[1]{%
\Y=#1%
\divide\Y by 2%
\Truehead%
\bimolength=#1%
\advance\bimolength by -\Trueepihead%
\epilength=\bimolength%
\advance\epilength by -\Y%
\begin{picture}(0,0)%
\put(\Y,\Y){\line(-1,-1){\bimolength}}%
\put(-\epilength,-\epilength){\swhead}%
\put(-\Y,-\Y){\swhead}%
\end{picture}}%

\newcommand{\SWBIMO}[1]{%
\Y=#1%
\divide\Y by 2%
\Truetail\Truehead%
\bimolength=#1%
\advance\bimolength by -\Truemonotail%
\monolength=\bimolength%
\advance\monolength by -\Y%
\advance\bimolength by -\Trueepihead%
\epilength=\bimolength%
\advance\epilength by -\monolength%
\begin{picture}(0,0)%
\put(\monolength,\monolength){\line(-1,-1){\bimolength}}%
\put(\monolength,\monolength){\swhead}%
\put(-\epilength,-\epilength){\swhead}%
\put(-\Y,-\Y){\swhead}%
\end{picture}}%

\newcommand{\SWBIAR}[1]{%
\Y=#1%
\divide\Y by 2%
\begin{picture}(0,0)%
\put(\Y,\Y){\begin{picture}(0,0)%
\truex{247}%
\put(\value{x},-\value{x}){\line(-1,-1){#1}}%
\put(-\value{x},\value{x}){\line(-1,-1){#1}}%
\monolength=#1%
\advance\monolength by -\value{x}%
\epilength=#1%
\advance\epilength by \value{x}%
\put(-\monolength,-\epilength){\swhead}%
\put(-\epilength,-\monolength){\swhead}%
\end{picture}}\end{picture}}%

\newcommand{\SWBIDIST}[1]{%
\Y=#1%
\divide\Y by 2%
\begin{picture}(0,0)%
\put(\Y,\Y){\begin{picture}(0,0)%
\truex{247}%
\monolength=#1%
\advance\monolength by -\value{x}%
\epilength=#1%
\advance\epilength by \value{x}%
\put(-\value{x},\value{x}){\line(-1,-1){#1}}%
\put(-\epilength,-\monolength){\swhead}%
\end{picture}}%
\put(\Y,\Y){\begin{picture}(0,0)%
\truex{247}%
\monolength=#1%
\advance\monolength by \value{x}%
\epilength=#1%
\advance\epilength by -\value{x}%
\put(\value{x},-\value{x}){\line(-1,-1){#1}}%
\put(-\epilength,-\monolength){\swhead}%
\end{picture}}%
\put(\value{x},-\value{x}){\distsign}%
\put(-\value{x},\value{x}){\distsign}%
\end{picture}}%

\newcommand{\SWADJAR}[1]{%
\Y=#1%
\divide\Y by 2%
\begin{picture}(0,0)%
\put(\Y,\Y){\begin{picture}(0,0)%
\truex{247}%
\monolength=#1%
\advance\monolength by -\value{x}%
\epilength=#1%
\advance\epilength by \value{x}%
\put(\value{x},-\value{x}){\line(-1,-1){#1}}%
\put(-\monolength,-\epilength){\swhead}%
\end{picture}}%
\put(-\Y,-\Y){\begin{picture}(0,0)%
\truex{247}%
\monolength=#1%
\advance\monolength by -\value{x}%
\epilength=#1%
\advance\epilength by \value{x}%
\put(-\value{x},\value{x}){\line(1,1){#1}}%
\put(\monolength,\epilength){\nehead}%
\end{picture}}\end{picture}}%

\newcommand{\SWADJDIST}[1]{%
\Y=#1%
\divide\Y by 2%
\begin{picture}(0,0)%
\put(\Y,\Y){\begin{picture}(0,0)%
\truex{247}%
\monolength=#1%
\advance\monolength by -\value{x}%
\epilength=#1%
\advance\epilength by \value{x}%
\put(\value{x},-\value{x}){\line(-1,-1){#1}}%
\put(-\monolength,-\epilength){\swhead}%
\end{picture}}%
\put(-\Y,-\Y){\begin{picture}(0,0)%
\truex{247}%
\monolength=#1%
\advance\monolength by -\value{x}%
\epilength=#1%
\advance\epilength by \value{x}%
\put(-\value{x},\value{x}){\line(1,1){#1}}%
\put(\monolength,\epilength){\nehead}%
\end{picture}}%
\put(-\value{x},\value{x}){\distsign}%
\put(\value{x},-\value{x}){\distsign}%
\end{picture}}%

\def\basicswar[#1]{\fdcase{\SWAR}{}{}{#100}}%

\newcommand{\swar}{\@ifnextchar[{\basicswar}{\basicswar[59]}}%

\def\basicSwar[#1]#2{\fdcase{\SWAR}{#2}{}{#100}}%

\newcommand{\Swar}{\@ifnextchar[{\basicSwar}{\basicSwar[59]}}%

\def\basicswaR[#1]#2{\fdcase{\SWAR}{}{#2}{#100}}%

\newcommand{\swaR}{\@ifnextchar[{\basicswaR}{\basicswaR[59]}}%

\def\basicswdist[#1]{\fdcase{\SWDIST}{}{}{#100}}%

\newcommand{\swdist}{\@ifnextchar[{\basicswdist}{\basicswdist[59]}}%

\def\basicSwdist[#1]#2{\fdcase{\SWDIST}{#2}{}{#100}}%

\newcommand{\Swdist}{\@ifnextchar[{\basicSwdist}{\basicSwdist[59]}}%

\def\basicswdisT[#1]#2{\fdcase{\SWDIST}{}{#2}{#100}}%

\newcommand{\swdisT}{\@ifnextchar[{\basicswdisT}{\basicswdisT[59]}}%

\def\basicswdotar[#1]{\fdcase{\SWDOTAR}{}{}{#100}}%

\newcommand{\swdotar}{\@ifnextchar[{\basicswdotar}{\basicswdotar[59]}}%

\def\basicSwdotar[#1]#2{\fdcase{\SWDOTAR}{#2}{}{#100}}%

\newcommand{\Swdotar}{\@ifnextchar[{\basicSwdotar}{\basicSwdotar[59]}}%

\def\basicswdotaR[#1]#2{\fdcase{\SWDOTAR}{}{#2}{#100}}%

\newcommand{\swdotaR}{\@ifnextchar[{\basicswdotaR}{\basicswdotaR[59]}}%

\def\basicswmono[#1]{\fdcase{\SWMONO}{}{}{#100}}%

\newcommand{\swmono}{\@ifnextchar[{\basicswmono}{\basicswmono[59]}}%

\def\basicSwmono[#1]#2{\fdcase{\SWMONO}{#2}{}{#100}}%

\newcommand{\Swmono}{\@ifnextchar[{\basicSwmono}{\basicSwmono[59]}}%

\def\basicswmonO[#1]#2{\fdcase{\SWMONO}{}{#2}{#100}}%

\newcommand{\swmonO}{\@ifnextchar[{\basicswmonO}{\basicswmonO[59]}}%

\def\basicswepi[#1]{\fdcase{\SWEPI}{}{}{#100}}%

\newcommand{\swepi}{\@ifnextchar[{\basicswepi}{\basicswepi[59]}}%

\def\basicSwepi[#1]#2{\fdcase{\SWEPI}{#2}{}{#100}}%

\newcommand{\Swepi}{\@ifnextchar[{\basicSwepi}{\basicSwepi[59]}}%

\def\basicswepI[#1]#2{\fdcase{\SWEPI}{}{#2}{#100}}%

\newcommand{\swepI}{\@ifnextchar[{\basicswepI}{\basicswepI[59]}}%

\def\basicswbimo[#1]{\fdcase{\SWBIMO}{}{}{#100}}%

\newcommand{\swbimo}{\@ifnextchar[{\basicswbimo}{\basicswbimo[59]}}%

\def\basicSwbimo[#1]#2{\fdcase{\SWBIMO}{#2}{}{#100}}%

\newcommand{\Swbimo}{\@ifnextchar[{\basicSwbimo}{\basicSwbimo[59]}}%

\def\basicswbimO[#1]#2{\fdcase{\SWBIMO}{}{#2}{#100}}%

\newcommand{\swbimO}{\@ifnextchar[{\basicswbimO}{\basicswbimO[59]}}%

\def\basicswiso[#1]{\fdcase{\SWAR}{\hspace{-2pt}\isosign}{}{#100}}%

\newcommand{\swiso}{\@ifnextchar[{\basicswiso}{\basicswiso[59]}}%

\def\basicSwiso[#1]#2{\fdcase{\SWAR}{#2}{\isosign}{#100}}%

\newcommand{\Swiso}{\@ifnextchar[{\basicSwiso}{\basicSwiso[59]}}%

\def\basicswisO[#1]#2{\fdcase{\SWAR}{\hspace{-2pt}\isosign}{#2}{#100}}%

\newcommand{\swisO}{\@ifnextchar[{\basicswisO}{\basicswisO[59]}}%

\def\basicswbiar[#1]{\fdbicase{\SWBIAR}{}{}{#100}}%

\newcommand{\swbiar}{\@ifnextchar[{\basicswbiar}{\basicswbiar[59]}}%

\def\basicSwbiar[#1]#2#3{\fdbicase{\SWBIAR}{#2}{#3}{#100}}%

\newcommand{\Swbiar}{\@ifnextchar[{\basicSwbiar}{\basicSwbiar[59]}}%

\def\basicswadjar[#1]{\fdbicase{\SWADJAR}{}{}{#100}}%

\newcommand{\swadjar}{\@ifnextchar[{\basicswadjar}{\basicswadjar[59]}}%

\def\basicSwadjar[#1]#2#3{\fdbicase{\SWADJAR}{#2}{#3}{#100}}%

\newcommand{\Swadjar}{\@ifnextchar[{\basicSwadjar}{\basicSwadjar[59]}}%

\def\basicswbidist[#1]{\fdbicase{\SWBIDIST}{}{}{#100}}%

\newcommand{\swbidist}{\@ifnextchar[{\basicswbidist}{\basicswbidist[59]}}%

\def\basicSwbidist[#1]#2#3{\fdbicase{\SWBIDIST}{#2}{#3}{#100}}%

\newcommand{\Swbidist}{\@ifnextchar[{\basicSwbidist}{\basicSwbidist[59]}}%

\def\basicswadjdist[#1]{\fdbicase{\SWADJDIST}{}{}{#100}}%

\newcommand{\swadjdist}{\@ifnextchar[{\basicswadjdist}{\basicswadjdist[59]}}%

\def\basicSwadjdist[#1]#2#3{\fdbicase{\SWADJDIST}{#2}{#3}{#100}}%

\newcommand{\Swadjdist}{\@ifnextchar[{\basicSwadjdist}{\basicSwadjdist[59]}}%

\newcommand{\sdcase}[4]{\begin{picture}(0,0)%
\put(0,0){#1{#4}}%
\truex{100}\truez{600}%
\put(\value{x},\value{x}){\makebox(0,\value{z})[l]{${#2}$}}%
\truex{300}\truey{800}%
\put(-\value{x},-\value{y}){\makebox(0,\value{z})[r]{${#3}$}}%
\end{picture}}%

\newcommand{\sdbicase}[4]{\begin{picture}(0,0)%
\put(0,0){#1{#4}}%
\truex{350}\truey{600}\truez{950}%
\put(\value{x},\value{x}){\makebox(0,\value{y})[l]{${#2}$}}%
\truex{450}\truey{600}\truez{1050}%
\put(-\value{x},-\value{z}){\makebox(0,\value{y})[r]{${#3}$}}%
\end{picture}}%

\newcommand{\SEAR}[1]{%
\Y=#1%
\divide\Y by 2%
\begin{picture}(0,0)%
\put(-\Y,\Y){\line(1,-1){#1}}%
\put(\Y,-\Y){\sehead}%
\end{picture}}%

\newcommand{\SEDIST}[1]{%
\Y=#1%
\divide\Y by 2%
\begin{picture}(0,0)%
\put(-\Y,\Y){\line(1,-1){#1}}%
\put(\Y,-\Y){\sehead}%
\put(0,0){\distsign}%
\end{picture}}%

\newcommand{\SEDOTAR}[1]%
{\truex{100}\truey{212}%
\Y=#1%
\divide\Y by 2%
\NUMBEROFDOTS=#1%
\divide\NUMBEROFDOTS by \value{y}%
\advance\NUMBEROFDOTS by 1%
\begin{picture}(0,0)%
\multiput(-\Y,\Y)(\value{y},-\value{y}){\NUMBEROFDOTS}%
{\circle*{\value{x}}}%
\put(\Y,-\Y){\sehead}%
\end{picture}}%

\newcommand{\SEMONO}[1]{%
\Y=#1%
\divide \Y by 2%
\Truetail%
\bimolength=#1%
\advance\bimolength by -\Truemonotail%
\monolength=\bimolength%
\advance\monolength by -\Y%
\begin{picture}(0,0)%
\put(-\monolength,\monolength){\line(1,-1){\bimolength}}%
\put(-\monolength,\monolength){\sehead}%
\put(\Y,-\Y){\sehead}%
\end{picture}}%

\newcommand{\SEEPI}[1]{%
\Y=#1%
\divide\Y by 2%
\Truehead%
\bimolength=#1%
\advance\bimolength by -\Trueepihead%
\epilength=\bimolength%
\advance\epilength by -\Y%
\begin{picture}(0,0)%
\put(-\Y,\Y){\line(1,-1){\bimolength}}%
\put(\epilength,-\epilength){\sehead}%
\put(\Y,-\Y){\sehead}%
\end{picture}}%

\newcommand{\SEBIMO}[1]{%
\Y=#1%
\divide\Y by 2%
\Truetail\Truehead%
\bimolength=#1%
\advance\bimolength by -\Truemonotail%
\monolength=\bimolength%
\advance\monolength by -\Y%
\advance\bimolength by -\Trueepihead%
\epilength=\bimolength%
\advance\epilength by -\monolength%
\begin{picture}(0,0)%
\put(-\monolength,\monolength){\line(1,-1){\bimolength}}%
\put(-\monolength,\monolength){\sehead}%
\put(\epilength,-\epilength){\sehead}%
\put(\Y,-\Y){\sehead}%
\end{picture}}%

\newcommand{\SEBIAR}[1]{%
\Y=#1%
\divide\Y by 2%
\begin{picture}(0,0)%
\put(-\Y,\Y){\begin{picture}(0,0)%
\truex{247}%
\put(-\value{x},-\value{x}){\line(1,-1){#1}}%
\put(\value{x},\value{x}){\line(1,-1){#1}}%
\monolength=#1%
\advance\monolength by -\value{x}%
\epilength=#1%
\advance\epilength by \value{x}%
\put(\monolength,-\epilength){\sehead}%
\put(\epilength,-\monolength){\sehead}%
\end{picture}}\end{picture}}%

\newcommand{\SEBIDIST}[1]{%
\Y=#1%
\divide\Y by 2%
\begin{picture}(0,0)%
\put(-\Y,\Y){\begin{picture}(0,0)%
\truex{247}%
\monolength=#1%
\advance\monolength by -\value{x}%
\epilength=#1%
\advance\epilength by \value{x}%
\put(\value{x},\value{x}){\line(1,-1){#1}}%
\put(\epilength,-\monolength){\sehead}%
\end{picture}}%
\put(-\Y,\Y){\begin{picture}(0,0)%
\truex{247}%
\monolength=#1%
\advance\monolength by \value{x}%
\epilength=#1%
\advance\epilength by -\value{x}%
\put(-\value{x},-\value{x}){\line(1,-1){#1}}%
\put(\epilength,-\monolength){\sehead}%
\end{picture}}%
\put(-\value{x},-\value{x}){\distsign}%
\put(\value{x},\value{x}){\distsign}%
\end{picture}}%

\newcommand{\SEEQL}[1]{%
\Y=#1%
\divide\Y by 2%
\begin{picture}(0,0)%
\put(-\Y,\Y){\begin{picture}(0,0)%
\truex{70}%
\put(-\value{x},-\value{x}){\line(1,-1){#1}}%
\put(\value{x},\value{x}){\line(1,-1){#1}}%
\end{picture}}\end{picture}}%

\newcommand{\SELINE}[1]{%
\Y=#1%
\divide\Y by 2%
\begin{picture}(0,0)%
\put(-\Y,\Y){\begin{picture}(0,0)%
\put(0,0){\line(1,-1){#1}}%
\end{picture}}\end{picture}}%

\newcommand{\SEADJAR}[1]{%
\Y=#1%
\divide\Y by 2%
\begin{picture}(0,0)%
\put(-\Y,\Y){\begin{picture}(0,0)%
\truex{247}%
\monolength=#1%
\advance\monolength by -\value{x}%
\epilength=#1%
\advance\epilength by \value{x}%
\put(-\value{x},-\value{x}){\line(1,-1){#1}}%
\put(\monolength,-\epilength){\sehead}%
\end{picture}}%
\put(\Y,-\Y){\begin{picture}(0,0)%
\truex{247}%
\monolength=#1%
\advance\monolength by -\value{x}%
\epilength=#1%
\advance\epilength by \value{x}%
\put(\value{x},\value{x}){\line(-1,1){#1}}%
\put(-\monolength,\epilength){\nwhead}%
\end{picture}}\end{picture}}%

\newcommand{\SEADJDIST}[1]{%
\Y=#1%
\divide\Y by 2%
\begin{picture}(0,0)%
\put(-\Y,\Y){\begin{picture}(0,0)%
\truex{247}%
\monolength=#1%
\advance\monolength by -\value{x}%
\epilength=#1%
\advance\epilength by \value{x}%
\put(-\value{x},-\value{x}){\line(1,-1){#1}}%
\put(\monolength,-\epilength){\sehead}%
\end{picture}}%
\put(\Y,-\Y){\begin{picture}(0,0)%
\truex{247}%
\monolength=#1%
\advance\monolength by -\value{x}%
\epilength=#1%
\advance\epilength by \value{x}%
\put(\value{x},\value{x}){\line(-1,1){#1}}%
\put(-\monolength,\epilength){\nwhead}%
\end{picture}}%
\put(-\value{x},-\value{x}){\distsign}%
\put(\value{x},\value{x}){\distsign}%
\end{picture}}%

\def\basicsear[#1]{\sdcase{\SEAR}{}{}{#100}}%

\newcommand{\sear}{\@ifnextchar[{\basicsear}{\basicsear[59]}}%

\def\basicSear[#1]#2{\sdcase{\SEAR}{#2}{}{#100}}%

\newcommand{\Sear}{\@ifnextchar[{\basicSear}{\basicSear[59]}}%

\def\basicseaR[#1]#2{\sdcase{\SEAR}{}{#2}{#100}}%

\newcommand{\seaR}{\@ifnextchar[{\basicseaR}{\basicseaR[59]}}%

\def\basicsedist[#1]{\sdcase{\SEDIST}{}{}{#100}}%

\newcommand{\sedist}{\@ifnextchar[{\basicsedist}{\basicsedist[59]}}%

\def\basicSedist[#1]#2{\sdcase{\SEDIST}{#2}{}{#100}}%

\newcommand{\Sedist}{\@ifnextchar[{\basicSedist}{\basicSedist[59]}}%

\def\basicsedisT[#1]#2{\sdcase{\SEDIST}{}{#2}{#100}}%

\newcommand{\sedisT}{\@ifnextchar[{\basicsedisT}{\basicsedisT[59]}}%

\def\basicsedotar[#1]{\sdcase{\SEDOTAR}{}{}{#100}}%

\newcommand{\sedotar}{\@ifnextchar[{\basicsedotar}{\basicsedotar[59]}}%

\def\basicSedotar[#1]#2{\sdcase{\SEDOTAR}{#2}{}{#100}}%

\newcommand{\Sedotar}{\@ifnextchar[{\basicSedotar}{\basicSedotar[59]}}%

\def\basicsedotaR[#1]#2{\sdcase{\SEDOTAR}{}{#2}{#100}}%

\newcommand{\sedotaR}{\@ifnextchar[{\basicsedotaR}{\basicsedotaR[59]}}%

\def\basicsemono[#1]{\sdcase{\SEMONO}{}{}{#100}}%

\newcommand{\semono}{\@ifnextchar[{\basicsemono}{\basicsemono[59]}}%

\def\basicSemono[#1]#2{\sdcase{\SEMONO}{#2}{}{#100}}%

\newcommand{\Semono}{\@ifnextchar[{\basicSemono}{\basicSemono[59]}}%

\def\basicsemonO[#1]#2{\sdcase{\SEMONO}{}{#2}{#100}}%

\newcommand{\semonO}{\@ifnextchar[{\basicsemonO}{\basicsemonO[59]}}%

\def\basicseepi[#1]{\sdcase{\SEEPI}{}{}{#100}}%

\newcommand{\seepi}{\@ifnextchar[{\basicseepi}{\basicseepi[59]}}%

\def\basicSeepi[#1]#2{\sdcase{\SEEPI}{#2}{}{#100}}%

\newcommand{\Seepi}{\@ifnextchar[{\basicSeepi}{\basicSeepi[59]}}%

\def\basicseepI[#1]#2{\sdcase{\SEEPI}{}{#2}{#100}}%

\newcommand{\seepI}{\@ifnextchar[{\basicseepI}{\basicseepI[59]}}%

\def\basicsebimo[#1]{\sdcase{\SEBIMO}{}{}{#100}}%

\newcommand{\sebimo}{\@ifnextchar[{\basicsebimo}{\basicsebimo[59]}}%

\def\basicSebimo[#1]#2{\sdcase{\SEBIMO}{#2}{}{#100}}%

\newcommand{\Sebimo}{\@ifnextchar[{\basicSebimo}{\basicSebimo[59]}}%

\def\basicsebimO[#1]#2{\sdcase{\SEBIMO}{}{#2}{#100}}%

\newcommand{\sebimO}{\@ifnextchar[{\basicsebimO}{\basicsebimO[59]}}%

\def\basicseiso[#1]{\sdcase{\SEAR}{\hspace{-2pt}\isosign}{}{#100}}%

\newcommand{\seiso}{\@ifnextchar[{\basicseiso}{\basicseiso[59]}}%

\def\basicSeiso[#1]#2{\sdcase{\SEAR}{#2}{\isosign}{#100}}%

\newcommand{\Seiso}{\@ifnextchar[{\basicSeiso}{\basicSeiso[59]}}%

\def\basicseisO[#1]#2{\sdcase{\SEAR}{\hspace{-2pt}\isosign}{#2}{#100}}%

\newcommand{\seisO}{\@ifnextchar[{\basicseisO}{\basicseisO[59]}}%

\def\basicseeql[#1]{\sdcase{\SEEQL}{}{}{#100}}%

\newcommand{\seeql}{\@ifnextchar[{\basicseeql}{\basicseeql[59]}}%

\def\basicSeeql[#1]#2{\sdcase{\SEEQL}{#2}{}{#100}}%

\newcommand{\Seeql}{\@ifnextchar[{\basicSeeql}{\basicSeeql[59]}}%

\def\basicseeqL[#1]#2{\sdcase{\SEEQL}{}{#2}{#100}}%

\newcommand{\seeqL}{\@ifnextchar[{\basicseeqL}{\basicseeqL[59]}}%

\def\basicseline[#1]{\sdcase{\SELINE}{}{}{#100}}%

\newcommand{\seline}{\@ifnextchar[{\basicseline}{\basicseline[59]}}%

\def\basicSeline[#1]#2{\sdcase{\SELINE}{#2}{}{#100}}%

\newcommand{\Seline}{\@ifnextchar[{\basicSeline}{\basicSeline[59]}}%

\def\basicselinE[#1]#2{\sdcase{\SELINE}{}{#2}{#100}}%

\newcommand{\selinE}{\@ifnextchar[{\basicselinE}{\basicselinE[59]}}%

\def\basicsebiar[#1]{\sdbicase{\SEBIAR}{}{}{#100}}%

\newcommand{\sebiar}{\@ifnextchar[{\basicsebiar}{\basicsebiar[59]}}%

\def\basicSebiar[#1]#2#3{\sdbicase{\SEBIAR}{#2}{#3}{#100}}%

\newcommand{\Sebiar}{\@ifnextchar[{\basicSebiar}{\basicSebiar[59]}}%

\def\basicseadjar[#1]{\sdbicase{\SEADJAR}{}{}{#100}}%

\newcommand{\seadjar}{\@ifnextchar[{\basicseadjar}{\basicseadjar[59]}}%

\def\basicSeadjar[#1]#2#3{\sdbicase{\SEADJAR}{#2}{#3}{#100}}%

\newcommand{\Seadjar}{\@ifnextchar[{\basicSeadjar}{\basicSeadjar[59]}}%

\def\basicsebidist[#1]{\sdbicase{\SEBIDIST}{}{}{#100}}%

\newcommand{\sebidist}{\@ifnextchar[{\basicsebidist}{\basicsebidist[59]}}%

\def\basicSebidist[#1]#2#3{\sdbicase{\SEBIDIST}{#2}{#3}{#100}}%

\newcommand{\Sebidist}{\@ifnextchar[{\basicSebidist}{\basicSebidist[59]}}%

\def\basicseadjdist[#1]{\sdbicase{\SEADJDIST}{}{}{#100}}%

\newcommand{\seadjdist}{\@ifnextchar[{\basicseadjdist}{\basicseadjdist[59]}}%

\def\basicSeadjdist[#1]#2#3{\sdbicase{\SEADJDIST}{#2}{#3}{#100}}%

\newcommand{\Seadjdist}{\@ifnextchar[{\basicSeadjdist}{\basicSeadjdist[59]}}%

\newcommand{\NWAR}[1]{%
\Y=#1%
\divide\Y by 2%
\begin{picture}(0,0)%
\put(\Y,-\Y){\line(-1,1){#1}}%
\put(-\Y,\Y){\nwhead}%
\end{picture}}%

\newcommand{\NWDIST}[1]{%
\Y=#1%
\divide\Y by 2%
\begin{picture}(0,0)%
\put(\Y,-\Y){\line(-1,1){#1}}%
\put(-\Y,\Y){\nwhead}%
\put(0,0){\distsign}%
\end{picture}}%

\newcommand{\NWDOTAR}[1]%
{\truex{100}\truey{212}%
\Y=#1%
\divide\Y by 2%
\NUMBEROFDOTS=#1%
\divide\NUMBEROFDOTS by \value{y}%
\advance\NUMBEROFDOTS by 1%
\begin{picture}(0,0)%
\multiput(\Y,-\Y)(-\value{y},\value{y}){\NUMBEROFDOTS}%
{\circle*{\value{x}}}%
\put(-\Y,\Y){\nwhead}%
\end{picture}}%

\newcommand{\NWMONO}[1]{%
\Y=#1%
\divide \Y by 2%
\Truetail%
\bimolength=#1%
\advance\bimolength by -\Truemonotail%
\monolength=\bimolength%
\advance\monolength by -\Y%
\begin{picture}(0,0)%
\put(\monolength,-\monolength){\line(-1,1){\bimolength}}%
\put(\monolength,-\monolength){\nwhead}%
\put(-\Y,\Y){\nwhead}%
\end{picture}}%

\newcommand{\NWEPI}[1]{%
\Y=#1%
\divide\Y by 2%
\Truehead%
\bimolength=#1%
\advance\bimolength by -\Trueepihead%
\epilength=\bimolength%
\advance\epilength by -\Y%
\begin{picture}(0,0)%
\put(\Y,-\Y){\line(-1,1){\bimolength}}%
\put(-\epilength,\epilength){\nwhead}%
\put(-\Y,\Y){\nwhead}%
\end{picture}}%

\newcommand{\NWBIMO}[1]{%
\Y=#1%
\divide\Y by 2%
\Truetail\Truehead%
\bimolength=#1%
\advance\bimolength by -\Truemonotail%
\monolength=\bimolength%
\advance\monolength by -\Y%
\advance\bimolength by -\Trueepihead%
\epilength=\bimolength%
\advance\epilength by -\monolength%
\begin{picture}(0,0)%
\put(\monolength,-\monolength){\line(-1,1){\bimolength}}%
\put(\monolength,-\monolength){\nwhead}%
\put(-\epilength,\epilength){\nwhead}%
\put(-\Y,\Y){\nwhead}%
\end{picture}}%

\newcommand{\NWBIAR}[1]{%
\Y=#1%
\divide\Y by 2%
\begin{picture}(0,0)%
\put(\Y,-\Y){\begin{picture}(0,0)%
\truex{247}%
\put(-\value{x},-\value{x}){\line(-1,1){#1}}%
\put(\value{x},\value{x}){\line(-1,1){#1}}%
\monolength=#1%
\advance\monolength by -\value{x}%
\epilength=#1%
\advance\epilength by \value{x}%
\put(-\monolength,\epilength){\nwhead}%
\put(-\epilength,\monolength){\nwhead}%
\end{picture}}\end{picture}}%

\newcommand{\NWBIDIST}[1]{%
\Y=#1%
\divide\Y by 2%
\begin{picture}(0,0)%
\put(\Y,-\Y){\begin{picture}(0,0)%
\truex{247}%
\monolength=#1%
\advance\monolength by -\value{x}%
\epilength=#1%
\advance\epilength by \value{x}%
\put(-\value{x},-\value{x}){\line(-1,1){#1}}%
\put(-\epilength,\monolength){\nwhead}%
\end{picture}}%
\put(\Y,-\Y){\begin{picture}(0,0)%
\truex{247}%
\monolength=#1%
\advance\monolength by \value{x}%
\epilength=#1%
\advance\epilength by -\value{x}%
\put(\value{x},\value{x}){\line(-1,1){#1}}%
\put(-\epilength,\monolength){\nwhead}%
\end{picture}}%
\put(-\value{x},-\value{x}){\distsign}%
\put(\value{x},\value{x}){\distsign}%
\end{picture}}%

\newcommand{\NWADJAR}[1]{%
\Y=#1%
\divide\Y by 2%
\begin{picture}(0,0)%
\put(\Y,-\Y){\begin{picture}(0,0)%
\truex{247}%
\monolength=#1%
\advance\monolength by -\value{x}%
\epilength=#1%
\advance\epilength by \value{x}%
\put(-\value{x},-\value{x}){\line(-1,1){#1}}%
\put(-\epilength,\monolength){\nwhead}%
\end{picture}}%
\put(-\Y,\Y){\begin{picture}(0,0)%
\truex{247}%
\monolength=#1%
\advance\monolength by -\value{x}%
\epilength=#1%
\advance\epilength by \value{x}%
\put(\value{x},\value{x}){\line(1,-1){#1}}%
\put(\epilength,-\monolength){\sehead}%
\end{picture}}\end{picture}}%

\newcommand{\NWADJDIST}[1]{%
\Y=#1%
\divide\Y by 2%
\begin{picture}(0,0)%
\put(\Y,-\Y){\begin{picture}(0,0)%
\truex{247}%
\monolength=#1%
\advance\monolength by -\value{x}%
\epilength=#1%
\advance\epilength by \value{x}%
\put(-\value{x},-\value{x}){\line(-1,1){#1}}%
\put(-\epilength,\monolength){\nwhead}%
\end{picture}}%
\put(-\Y,\Y){\begin{picture}(0,0)%
\truex{247}%
\monolength=#1%
\advance\monolength by -\value{x}%
\epilength=#1%
\advance\epilength by \value{x}%
\put(\value{x},\value{x}){\line(1,-1){#1}}%
\put(\epilength,-\monolength){\sehead}%
\end{picture}}%
\put(-\value{x},-\value{x}){\distsign}%
\put(\value{x},\value{x}){\distsign}%
\end{picture}}%

\def\basicnwar[#1]{\sdcase{\NWAR}{}{}{#100}}%

\newcommand{\nwar}{\@ifnextchar[{\basicnwar}{\basicnwar[59]}}%

\def\basicNwar[#1]#2{\sdcase{\NWAR}{#2}{}{#100}}%

\newcommand{\Nwar}{\@ifnextchar[{\basicNwar}{\basicNwar[59]}}%

\def\basicnwaR[#1]#2{\sdcase{\NWAR}{}{#2}{#100}}%

\newcommand{\nwaR}{\@ifnextchar[{\basicnwaR}{\basicnwaR[59]}}%

\def\basicnwdist[#1]{\sdcase{\NWDIST}{}{}{#100}}%

\newcommand{\nwdist}{\@ifnextchar[{\basicnwdist}{\basicnwdist[59]}}%

\def\basicNwdist[#1]#2{\sdcase{\NWDIST}{#2}{}{#100}}%

\newcommand{\Nwdist}{\@ifnextchar[{\basicNwdist}{\basicNwdist[59]}}%

\def\basicnwdisT[#1]#2{\sdcase{\NWDIST}{}{#2}{#100}}%

\newcommand{\nwdisT}{\@ifnextchar[{\basicnwdisT}{\basicnwdisT[59]}}%

\def\basicnwdotar[#1]{\sdcase{\NWDOTAR}{}{}{#100}}%

\newcommand{\nwdotar}{\@ifnextchar[{\basicnwdotar}{\basicnwdotar[59]}}%

\def\basicNwdotar[#1]#2{\sdcase{\NWDOTAR}{#2}{}{#100}}%

\newcommand{\Nwdotar}{\@ifnextchar[{\basicNwdotar}{\basicNwdotar[59]}}%

\def\basicnwdotaR[#1]#2{\sdcase{\NWDOTAR}{}{#2}{#100}}%

\newcommand{\nwdotaR}{\@ifnextchar[{\basicnwdotaR}{\basicnwdotaR[59]}}%

\def\basicnwmono[#1]{\sdcase{\NWMONO}{}{}{#100}}%

\newcommand{\nwmono}{\@ifnextchar[{\basicnwmono}{\basicnwmono[59]}}%

\def\basicNwmono[#1]#2{\sdcase{\NWMONO}{#2}{}{#100}}%

\newcommand{\Nwmono}{\@ifnextchar[{\basicNwmono}{\basicNwmono[59]}}%

\def\basicnwmonO[#1]#2{\sdcase{\NWMONO}{}{#2}{#100}}%

\newcommand{\nwmonO}{\@ifnextchar[{\basicnwmonO}{\basicnwmonO[59]}}%

\def\basicnwepi[#1]{\sdcase{\NWEPI}{}{}{#100}}%

\newcommand{\nwepi}{\@ifnextchar[{\basicnwepi}{\basicnwepi[59]}}%

\def\basicNwepi[#1]#2{\sdcase{\NWEPI}{#2}{}{#100}}%

\newcommand{\Nwepi}{\@ifnextchar[{\basicNwepi}{\basicNwepi[59]}}%

\def\basicnwepI[#1]#2{\sdcase{\NWEPI}{}{#2}{#100}}%

\newcommand{\nwepI}{\@ifnextchar[{\basicnwepI}{\basicnwepI[59]}}%

\def\basicnwbimo[#1]{\sdcase{\NWBIMO}{}{}{#100}}%

\newcommand{\nwbimo}{\@ifnextchar[{\basicnwbimo}{\basicnwbimo[59]}}%

\def\basicNwbimo[#1]#2{\sdcase{\NWBIMO}{#2}{}{#100}}%

\newcommand{\Nwbimo}{\@ifnextchar[{\basicNwbimo}{\basicNwbimo[59]}}%

\def\basicnwbimO[#1]#2{\sdcase{\NWBIMO}{}{#2}{#100}}%

\newcommand{\nwbimO}{\@ifnextchar[{\basicnwbimO}{\basicnwbimO[59]}}%

\def\basicnwiso[#1]{\sdcase{\NWAR}{\hspace{-2pt}\isosign}{}{#100}}%

\newcommand{\nwiso}{\@ifnextchar[{\basicnwiso}{\basicnwiso[59]}}%

\def\basicNwiso[#1]#2{\sdcase{\NWAR}{#2}{\isosign}{#100}}%

\newcommand{\Nwiso}{\@ifnextchar[{\basicNwiso}{\basicNwiso[59]}}%

\def\basicnwisO[#1]#2{\sdcase{\NWAR}{\hspace{-2pt}\isosign}{#2}{#100}}%

\newcommand{\nwisO}{\@ifnextchar[{\basicnwisO}{\basicnwisO[59]}}%

\def\basicnwbiar[#1]{\sdbicase{\NWBIAR}{}{}{#100}}%

\newcommand{\nwbiar}{\@ifnextchar[{\basicnwbiar}{\basicnwbiar[59]}}%

\def\basicNwbiar[#1]#2#3{\sdbicase{\NWBIAR}{#2}{#3}{#100}}%

\newcommand{\Nwbiar}{\@ifnextchar[{\basicNwbiar}{\basicNwbiar[59]}}%

\def\basicnwadjar[#1]{\sdbicase{\NWADJAR}{}{}{#100}}%

\newcommand{\nwadjar}{\@ifnextchar[{\basicnwadjar}{\basicnwadjar[59]}}%

\def\basicNwadjar[#1]#2#3{\sdbicase{\NWADJAR}{#2}{#3}{#100}}%

\newcommand{\Nwadjar}{\@ifnextchar[{\basicNwadjar}{\basicNwadjar[59]}}%

\def\basicnwbidist[#1]{\sdbicase{\NWBIDIST}{}{}{#100}}%

\newcommand{\nwbidist}{\@ifnextchar[{\basicnwbidist}{\basicnwbidist[59]}}%

\def\basicNwbidist[#1]#2#3{\sdbicase{\NWBIDIST}{#2}{#3}{#100}}%

\newcommand{\Nwbidist}{\@ifnextchar[{\basicNwbidist}{\basicNwbidist[59]}}%

\def\basicnwadjdist[#1]{\sdbicase{\NWADJDIST}{}{}{#100}}%

\newcommand{\nwadjdist}{\@ifnextchar[{\basicnwadjdist}{\basicnwadjdist[59]}}%

\def\basicNwadjdist[#1]#2#3{\sdbicase{\NWADJDIST}{#2}{#3}{#100}}%

\newcommand{\Nwadjdist}{\@ifnextchar[{\basicNwadjdist}{\basicNwadjdist[59]}}%

\newcommand{\ENEAR}[3]{%
\Y=#3%
\divide\Y by 2%
\Z=\Y%
\divide\Z by 2%
\begin{picture}(0,0)%
\put(-\Y,-\Z){\line(2,1){#3}}%
\put(\Y,\Z){\enehead}%
\truex{200}\truey{800}\truez{600}%
\put(-\value{x},\value{x}){\makebox(0,\value{z})[r]{${#1}$}}%
\put(\value{x},-\value{y}){\makebox(0,\value{z})[l]{${#2}$}}%
\end{picture}}%

\newcommand{\ENEDIST}[3]{%
\Y=#3%
\divide\Y by 2%
\Z=\Y%
\divide\Z by 2%
\begin{picture}(0,0)%
\put(-\Y,-\Z){\line(2,1){#3}}%
\put(\Y,\Z){\enehead}%
\put(0,0){\distsign}%
\truex{200}\truey{800}\truez{600}%
\put(-\value{x},\value{x}){\makebox(0,\value{z})[r]{${#1}$}}%
\put(\value{x},-\value{y}){\makebox(0,\value{z})[l]{${#2}$}}%
\end{picture}}%

\newcommand{\ENEDOTAR}[3]{%
\truex{100}\truey{268}\truez{134}%
\Y=#3%
\divide\Y by 2%
\Z=\Y%
\divide\Z by 2%
\NUMBEROFDOTS=#3%
\divide\NUMBEROFDOTS by \value{y}%
\advance\NUMBEROFDOTS by 1%
\begin{picture}(0,0)%
\multiput(-\Y,-\Z)(\value{y},\value{z}){\NUMBEROFDOTS}%
{\circle*{\value{x}}}%
\put(\Y,\Z){\enehead}%
\truex{200}\truey{800}\truez{600}%
\put(-\value{x},\value{x}){\makebox(0,\value{z})[r]{${#1}$}}%
\put(\value{x},-\value{y}){\makebox(0,\value{z})[l]{${#2}$}}%
\end{picture}}%

\newcommand{\ENEMONO}[3]{%
\Y=#3%
\divide\Y by 2%
\Z=\Y%
\divide\Z by 2%
\TrueTail%
\bimolength=#3%
\advance\bimolength by -\TrueMonoTail%
\monolength=\bimolength%
\advance\monolength by -\Y%
\secondmonolength=\monolength%
\divide\secondmonolength by 2%
\begin{picture}(0,0)%
\put(-\monolength,-\secondmonolength){\line(2,1){\bimolength}}%
\put(-\monolength,-\secondmonolength){\enehead}%
\put(\Y,\Z){\enehead}%
\truex{200}\truey{800}\truez{600}%
\put(-\value{x},\value{x}){\makebox(0,\value{z})[r]{${#1}$}}%
\put(\value{x},-\value{y}){\makebox(0,\value{z})[l]{${#2}$}}%
\end{picture}}%

\newcommand{\ENEEPI}[3]{%
\Y=#3%
\divide\Y by 2%
\Z=\Y%
\divide\Z by 2%
\TrueHead%
\bimolength=#3%
\advance\bimolength by -\TrueEpiHead%
\epilength=\bimolength%
\advance\epilength by -\Y%
\secondepilength=\epilength%
\divide\secondepilength by 2%
\begin{picture}(0,0)%
\put(-\Y,-\Z){\line(2,1){\bimolength}}%
\put(\epilength,\secondepilength){\enehead}%
\put(\Y,\Z){\enehead}%
\truex{200}\truey{800}\truez{600}%
\put(-\value{x},\value{x}){\makebox(0,\value{z})[r]{${#1}$}}%
\put(\value{x},-\value{y}){\makebox(0,\value{z})[l]{${#2}$}}%
\end{picture}}%

\newcommand{\ENEBIMO}[3]{%
\Y=#3%
\divide\Y by 2%
\Z=\Y%
\divide\Z by 2%
\TrueTail\TrueHead%
\bimolength=#3%
\advance\bimolength by -\TrueMonoTail%
\monolength=\bimolength%
\advance\monolength by -\Y%
\advance\bimolength by -\TrueEpiHead%
\epilength=\bimolength%
\advance\epilength by -\monolength%
\secondmonolength=\monolength%
\divide\secondmonolength by 2%
\secondepilength=\epilength%
\divide\secondepilength by 2%
\begin{picture}(0,0)%
\put(-\monolength,-\secondmonolength){\line(2,1){\bimolength}}%
\put(-\monolength,-\secondmonolength){\enehead}%
\put(\epilength,\secondepilength){\enehead}%
\put(\Y,\Z){\enehead}%
\truex{200}\truey{800}\truez{600}%
\put(-\value{x},\value{x}){\makebox(0,\value{z})[r]{${#1}$}}%
\put(\value{x},-\value{y}){\makebox(0,\value{z})[l]{${#2}$}}%
\end{picture}}%

\newcommand{\ENEEQL}[3]{%
\Y=#3%
\divide\Y by 2%
\Z=\Y%
\divide\Z by 2%
\begin{picture}(0,0)%
\put(-\Y,-\Z){\begin{picture}(0,0)%
\truex{44}\truey{89}%
\put(-\value{x},\value{y}){\line(2,1){#3}}%
\put(\value{x},-\value{y}){\line(2,1){#3}}%
\end{picture}}%
\truex{200}\truey{800}\truez{600}%
\put(-\value{x},\value{x}){\makebox(0,\value{z})[r]{${#1}$}}%
\put(\value{x},-\value{y}){\makebox(0,\value{z})[l]{${#2}$}}%
\end{picture}}%

\newcommand{\ENEBIAR}[3]{%
\Y=#3%
\divide\Y by 2%
\Z=\Y%
\divide\Z by 2%
\begin{picture}(0,0)%
\put(-\Y,-\Z){\begin{picture}(0,0)%
\truex{156}\truey{313}%
\put(-\value{x},\value{y}){\line(2,1){#3}}%
\put(\value{x},-\value{y}){\line(2,1){#3}}%
\monolength=#3%
\advance\monolength by -\value{x}%
\epilength=#3%
\advance\epilength by \value{x}%
\secondmonolength=\Y%
\advance\secondmonolength by -\value{y}%
\secondepilength=\Y%
\advance\secondepilength by \value{y}%
\put(\monolength,\secondepilength){\enehead}%
\put(\epilength,\secondmonolength){\enehead}%
\end{picture}}
\truex{300}\truey{1000}\truez{600}%
\put(-\value{x},\value{x}){\makebox(0,\value{z})[r]{${#1}$}}%
\put(\value{x},-\value{y}){\makebox(0,\value{z})[l]{${#2}$}}%
\end{picture}}%

\newcommand{\ENEBIDIST}[3]{%
\Y=#3%
\divide\Y by 2%
\Z=\Y%
\divide\Z by 2%
\begin{picture}(0,0)%
\truex{156}\truey{313}\truez{400}%
\put(-\Y,-\Z){\begin{picture}(0,0)%
\put(-\value{x},\value{y}){\line(2,1){#3}}%
\put(\value{x},-\value{y}){\line(2,1){#3}}%
\monolength=#3%
\advance\monolength by -\value{x}%
\epilength=#3%
\advance\epilength by \value{x}%
\secondmonolength=\Y%
\advance\secondmonolength by -\value{y}%
\secondepilength=\Y%
\advance\secondepilength by \value{y}%
\put(\monolength,\secondepilength){\enehead}%
\put(\epilength,\secondmonolength){\enehead}%
\end{picture}}
\put(-\value{x},\value{y}){\distsign}%
\put(\value{x},-\value{y}){\distsign}%
\truex{300}\truey{1000}\truez{600}%
\put(-\value{x},\value{x}){\makebox(0,\value{z})[r]{${#1}$}}%
\put(\value{x},-\value{y}){\makebox(0,\value{z})[l]{${#2}$}}%
\end{picture}}%

\newcommand{\ENEADJAR}[3]{%
\Y=#3%
\divide\Y by 2%
\Z=\Y%
\divide\Z by 2%
\begin{picture}(0,0)%
\put(-\Y,-\Z){\begin{picture}(0,0)%
\truex{156}\truey{313}%
\monolength=#3%
\advance\monolength by -\value{x}%
\epilength=#3%
\advance\epilength by \value{x}%
\secondmonolength=\Y%
\advance\secondmonolength by -\value{y}%
\secondepilength=\Y%
\advance\secondepilength by \value{y}%
\put(\value{x},-\value{y}){\line(2,1){#3}}%
\put(\epilength,\secondmonolength){\enehead}%
\put(\monolength,\secondepilength){\line(-2,-1){#3}}%
\put(-\value{x},\value{y}){\wswhead}%
\end{picture}}
\truex{300}\truey{1000}\truez{600}%
\put(-\value{x},\value{x}){\makebox(0,\value{z})[r]{${#1}$}}%
\put(\value{x},-\value{y}){\makebox(0,\value{z})[l]{${#2}$}}%
\end{picture}}%

\newcommand{\ENEADJDIST}[3]{%
\Y=#3%
\divide\Y by 2%
\Z=\Y%
\divide\Z by 2%
\begin{picture}(0,0)%
\truex{156}\truey{313}\truez{400}%
\put(-\Y,-\Z){\begin{picture}(0,0)%
\monolength=#3%
\advance\monolength by -\value{x}%
\epilength=#3%
\advance\epilength by \value{x}%
\secondmonolength=\Y%
\advance\secondmonolength by -\value{y}%
\secondepilength=\Y%
\advance\secondepilength by \value{y}%
\put(\value{x},-\value{y}){\line(2,1){#3}}%
\put(\epilength,\secondmonolength){\enehead}%
\put(\monolength,\secondepilength){\line(-2,-1){#3}}%
\put(-\value{x},\value{y}){\wswhead}%
\end{picture}}
\put(-\value{x},\value{y}){\distsign}%
\put(\value{x},-\value{y}){\distsign}%
\truex{300}\truey{1000}\truez{600}%
\put(-\value{x},\value{x}){\makebox(0,\value{z})[r]{${#1}$}}%
\put(\value{x},-\value{y}){\makebox(0,\value{z})[l]{${#2}$}}%
\end{picture}}%

\def\basicenear[#1]{\ENEAR{}{}{#100}}%

\newcommand{\enear}{\@ifnextchar[{\basicenear}{\basicenear[133]}}%

\def\basicEnear[#1]#2{\ENEAR{#2}{}{#100}}%

\newcommand{\Enear}{\@ifnextchar[{\basicEnear}{\basicEnear[133]}}%

\def\basiceneaR[#1]#2{\ENEAR{}{#2}{#100}}%

\newcommand{\eneaR}{\@ifnextchar[{\basiceneaR}{\basiceneaR[133]}}%

\def\basicenedist[#1]{\ENEDIST{}{}{#100}}%

\newcommand{\enedist}{\@ifnextchar[{\basicenedist}{\basicenedist[133]}}%

\def\basicEnedist[#1]#2{\ENEDIST{#2}{}{#100}}%

\newcommand{\Enedist}{\@ifnextchar[{\basicEnedist}{\basicEnedist[133]}}%

\def\basicenedisT[#1]#2{\ENEDIST{}{#2}{#100}}%

\newcommand{\enedisT}{\@ifnextchar[{\basicenedisT}{\basicenedisT[133]}}%

\def\basicenedotar[#1]{\ENEDOTAR{}{}{#100}}%

\newcommand{\enedotar}{\@ifnextchar[{\basicenedotar}{\basicenedotar[133]}}%

\def\basicEnedotar[#1]#2{\ENEDOTAR{#2}{}{#100}}%

\newcommand{\Enedotar}{\@ifnextchar[{\basicEnedotar}{\basicEnedotar[133]}}%

\def\basicenedotaR[#1]#2{\ENEDOTAR{}{#2}{#100}}%

\newcommand{\enedotaR}{\@ifnextchar[{\basicenedotaR}{\basicenedotaR[133]}}%

\def\basicenemono[#1]{\ENEMONO{}{}{#100}}%

\newcommand{\enemono}{\@ifnextchar[{\basicenemono}{\basicenemono[133]}}%

\def\basicEnemono[#1]#2{\ENEMONO{#2}{}{#100}}%

\newcommand{\Enemono}{\@ifnextchar[{\basicEnemono}{\basicEnemono[133]}}%

\def\basicenemonO[#1]#2{\ENEMONO{}{#2}{#100}}%

\newcommand{\enemonO}{\@ifnextchar[{\basicenemonO}{\basicenemonO[133]}}%

\def\basiceneepi[#1]{\ENEEPI{}{}{#100}}%

\newcommand{\eneepi}{\@ifnextchar[{\basiceneepi}{\basiceneepi[133]}}%

\def\basicEneepi[#1]#2{\ENEEPI{#2}{}{#100}}%

\newcommand{\Eneepi}{\@ifnextchar[{\basicEneepi}{\basicEneepi[133]}}%

\def\basiceneepI[#1]#2{\ENEEPI{}{#2}{#100}}%

\newcommand{\eneepI}{\@ifnextchar[{\basiceneepI}{\basiceneepI[133]}}%

\def\basicenebimo[#1]{\ENEBIMO{}{}{#100}}%

\newcommand{\enebimo}{\@ifnextchar[{\basicenebimo}{\basicenebimo[133]}}%

\def\basicEnebimo[#1]#2{\ENEBIMO{#2}{}{#100}}%

\newcommand{\Enebimo}{\@ifnextchar[{\basicEnebimo}{\basicEnebimo[133]}}%

\def\basicenebimO[#1]#2{\ENEBIMO{}{#2}{#100}}%

\newcommand{\enebimO}{\@ifnextchar[{\basicenebimO}{\basicenebimO[133]}}%

\def\basiceneiso[#1]{\ENEAR{\isosign}{}{#100}}%

\newcommand{\eneiso}{\@ifnextchar[{\basiceneiso}{\basiceneiso[133]}}%

\def\basicEneiso[#1]#2{\ENEAR{#2}{\isosign}{#100}}%

\newcommand{\Eneiso}{\@ifnextchar[{\basicEneiso}{\basicEneiso[133]}}%

\def\basiceneisO[#1]#2{\ENEAR{\isosign}{#2}{#100}}%

\newcommand{\eneisO}{\@ifnextchar[{\basiceneisO}{\basiceneisO[133]}}%

\def\basiceneeql[#1]{\ENEEQL{}{}{#100}}%

\newcommand{\eneeql}{\@ifnextchar[{\basiceneeql}{\basiceneeql[133]}}%

\def\basicEneeql[#1]#2{\ENEEQL{#2}{}{#100}}%

\newcommand{\Eneeql}{\@ifnextchar[{\basicEneeql}{\basicEneeql[133]}}%

\def\basiceneeqL[#1]#2{\ENEEQL{}{#2}{#100}}%

\newcommand{\eneeqL}{\@ifnextchar[{\basiceneeqL}{\basiceneeqL[133]}}%

\def\basicenebiar[#1]{\ENEBIAR{}{}{#100}}%

\newcommand{\enebiar}{\@ifnextchar[{\basicenebiar}{\basicenebiar[133]}}%

\def\basicEnebiar[#1]#2#3{\ENEBIAR{#2}{#3}{#100}}%

\newcommand{\Enebiar}{\@ifnextchar[{\basicEnebiar}{\basicEnebiar[133]}}%

\def\basicenebidist[#1]{\ENEBIDIST{}{}{#100}}%

\newcommand{\enebidist}{\@ifnextchar[{\basicenebidist}{\basicenebidist[133]}}%

\def\basicEnebidist[#1]#2#3{\ENEBIDIST{#2}{#3}{#100}}%

\newcommand{\Enebidist}{\@ifnextchar[{\basicEnebidist}{\basicEnebidist[133]}}%

\def\basiceneadjar[#1]{\ENEADJAR{}{}{#100}}%

\newcommand{\eneadjar}{\@ifnextchar[{\basiceneadjar}{\basiceneadjar[133]}}%

\def\basicEneadjar[#1]#2#3{\ENEADJAR{#2}{#3}{#100}}%

\newcommand{\Eneadjar}{\@ifnextchar[{\basicEneadjar}{\basicEneadjar[133]}}%

\def\basiceneadjdist[#1]{\ENEADJDIST{}{}{#100}}%

\newcommand{\eneadjdist}{\@ifnextchar[{\basiceneadjdist}{\basiceneadjdist[13
3]}}%

\def\basicEneadjdist[#1]#2#3{\ENEADJDIST{#2}{#3}{#100}}%

\newcommand{\Eneadjdist}{\@ifnextchar[{\basicEneadjdist}{\basicEneadjdist[13
3]}}%

\newcommand{\ESEAR}[3]{%
\Y=#3%
\divide\Y by 2%
\Z=\Y%
\divide\Z by 2%
\begin{picture}(0,0)%
\put(-\Y,\Z){\line(2,-1){#3}}%
\put(\Y,-\Z){\esehead}%
\truex{200}\truey{800}\truez{600}%
\put(\value{x},\value{x}){\makebox(0,\value{z})[l]{${#1}$}}%
\put(-\value{x},-\value{y}){\makebox(0,\value{z})[r]{${#2}$}}%
\end{picture}}%

\newcommand{\ESEDIST}[3]{%
\Y=#3%
\divide\Y by 2%
\Z=\Y%
\divide\Z by 2%
\begin{picture}(0,0)%
\put(-\Y,\Z){\line(2,-1){#3}}%
\put(\Y,-\Z){\esehead}%
\put(0,0){\distsign}%
\truex{200}\truey{800}\truez{600}%
\put(\value{x},\value{x}){\makebox(0,\value{z})[l]{${#1}$}}%
\put(-\value{x},-\value{y}){\makebox(0,\value{z})[r]{${#2}$}}%
\end{picture}}%

\newcommand{\ESEDOTAR}[3]{%
\truex{100}\truey{268}\truez{134}%
\Y=#3%
\divide\Y by 2%
\Z=\Y%
\divide\Z by 2%
\NUMBEROFDOTS=#3%
\divide\NUMBEROFDOTS by \value{y}%
\advance\NUMBEROFDOTS by 1%
\begin{picture}(0,0)%
\multiput(-\Y,\Z)(\value{y},-\value{z}){\NUMBEROFDOTS}%
{\circle*{\value{x}}}%
\put(\Y,-\Z){\esehead}%
\truex{200}\truey{800}\truez{600}%
\put(\value{x},\value{x}){\makebox(0,\value{z})[l]{${#1}$}}%
\put(-\value{x},-\value{y}){\makebox(0,\value{z})[r]{${#2}$}}%
\end{picture}}%

\newcommand{\ESEMONO}[3]{%
\Y=#3%
\divide\Y by 2%
\Z=\Y%
\divide\Z by 2%
\TrueTail%
\bimolength=#3%
\advance\bimolength by -\TrueMonoTail%
\monolength=\bimolength%
\advance\monolength by -\Y%
\secondmonolength=\monolength%
\divide\secondmonolength by 2%
\begin{picture}(0,0)%
\put(-\monolength,\secondmonolength){\line(2,-1){\bimolength}}%
\put(-\monolength,\secondmonolength){\esehead}%
\put(\Y,-\Z){\esehead}%
\truex{200}\truey{800}\truez{600}%
\put(\value{x},\value{x}){\makebox(0,\value{z})[l]{${#1}$}}%
\put(-\value{x},-\value{y}){\makebox(0,\value{z})[r]{${#2}$}}%
\end{picture}}%

\newcommand{\ESEEPI}[3]{%
\Y=#3%
\divide\Y by 2%
\Z=\Y%
\divide\Z by 2%
\TrueHead%
\bimolength=#3%
\advance\bimolength by -\TrueEpiHead%
\epilength=\bimolength%
\advance\epilength by -\Y%
\secondepilength=\epilength%
\divide\secondepilength by 2%
\begin{picture}(0,0)%
\put(-\Y,\Z){\line(2,-1){\bimolength}}%
\put(\epilength,-\secondepilength){\esehead}%
\put(\Y,-\Z){\esehead}%
\truex{200}\truey{800}\truez{600}%
\put(\value{x},\value{x}){\makebox(0,\value{z})[l]{${#1}$}}%
\put(-\value{x},-\value{y}){\makebox(0,\value{z})[r]{${#2}$}}%
\end{picture}}%

\newcommand{\ESEBIMO}[3]{%
\Y=#3%
\divide\Y by 2%
\Z=\Y%
\divide\Z by 2%
\TrueTail\TrueHead%
\bimolength=#3%
\advance\bimolength by -\TrueMonoTail%
\monolength=\bimolength%
\advance\monolength by -\Y%
\advance\bimolength by -\TrueEpiHead%
\epilength=\bimolength%
\advance\epilength by -\monolength%
\secondmonolength=\monolength%
\divide\secondmonolength by 2%
\secondepilength=\epilength%
\divide\secondepilength by 2%
\begin{picture}(0,0)%
\put(-\monolength,\secondmonolength){\line(2,-1){\bimolength}}%
\put(-\monolength,\secondmonolength){\esehead}%
\put(\epilength,-\secondepilength){\esehead}%
\put(\Y,-\Z){\esehead}%
\truex{200}\truey{800}\truez{600}%
\put(\value{x},\value{x}){\makebox(0,\value{z})[l]{${#1}$}}%
\put(-\value{x},-\value{y}){\makebox(0,\value{z})[r]{${#2}$}}%
\end{picture}}%

\newcommand{\ESEEQL}[3]{%
\Y=#3%
\divide\Y by 2%
\Z=\Y%
\divide\Z by 2%
\begin{picture}(0,0)%
\put(-\Y,\Z){\begin{picture}(0,0)%
\truex{44}\truey{89}%
\put(-\value{x},-\value{y}){\line(2,-1){#3}}%
\put(\value{x},\value{y}){\line(2,-1){#3}}%
\end{picture}}%
\truex{200}\truey{800}\truez{600}%
\put(\value{x},\value{x}){\makebox(0,\value{z})[l]{${#1}$}}%
\put(-\value{x},-\value{y}){\makebox(0,\value{z})[r]{${#2}$}}%
\end{picture}}%

\newcommand{\ESEBIAR}[3]{%
\Y=#3%
\divide\Y by 2%
\Z=\Y%
\divide\Z by 2%
\begin{picture}(0,0)%
\put(-\Y,\Z){\begin{picture}(0,0)%
\truex{156}\truey{313}%
\put(-\value{x},-\value{y}){\line(2,-1){#3}}%
\put(\value{x},\value{y}){\line(2,-1){#3}}%
\monolength=#3%
\advance\monolength by -\value{x}%
\epilength=#3%
\advance\epilength by \value{x}%
\secondmonolength=\Y%
\advance\secondmonolength by -\value{y}%
\secondepilength=\Y%
\advance\secondepilength by \value{y}%
\put(\monolength,-\secondepilength){\esehead}%
\put(\epilength,-\secondmonolength){\esehead}%
\end{picture}}
\truex{400}\truey{1000}\truez{600}%
\put(\value{x},\value{x}){\makebox(0,\value{z})[l]{${#1}$}}%
\put(-\value{x},-\value{y}){\makebox(0,\value{z})[r]{${#2}$}}%
\end{picture}}%

\newcommand{\ESEBIDIST}[3]{%
\Y=#3%
\divide\Y by 2%
\Z=\Y%
\divide\Z by 2%
\begin{picture}(0,0)%
\truex{156}\truey{313}\truez{400}%
\put(-\Y,\Z){\begin{picture}(0,0)%
\put(-\value{x},-\value{y}){\line(2,-1){#3}}%
\put(\value{x},\value{y}){\line(2,-1){#3}}%
\monolength=#3%
\advance\monolength by -\value{x}%
\epilength=#3%
\advance\epilength by \value{x}%
\secondmonolength=\Y%
\advance\secondmonolength by -\value{y}%
\secondepilength=\Y%
\advance\secondepilength by \value{y}%
\put(\monolength,-\secondepilength){\esehead}%
\put(\epilength,-\secondmonolength){\esehead}%
\end{picture}}
\put(\value{x},\value{y}){\distsign}%
\put(-\value{x},-\value{y}){\distsign}%
\truex{400}\truey{1000}\truez{600}%
\put(\value{x},\value{x}){\makebox(0,\value{z})[l]{${#1}$}}%
\put(-\value{x},-\value{y}){\makebox(0,\value{z})[r]{${#2}$}}%
\end{picture}}%

\newcommand{\ESEADJAR}[3]{%
\Y=#3%
\divide\Y by 2%
\Z=\Y%
\divide\Z by 2%
\begin{picture}(0,0)%
\put(-\Y,\Z){\begin{picture}(0,0)%
\truex{156}\truey{313}%
\monolength=#3%
\advance\monolength by -\value{x}%
\epilength=#3%
\advance\epilength by \value{x}%
\secondmonolength=\Y%
\advance\secondmonolength by -\value{y}%
\secondepilength=\Y%
\advance\secondepilength by \value{y}%
\put(-\value{x},-\value{y}){\line(2,-1){#3}}%
\put(\monolength,-\secondepilength){\esehead}%
\put(\epilength,-\secondmonolength){\line(-2,1){#3}}%
\put(\value{x},\value{y}){\wnwhead}%
\end{picture}}
\truex{400}\truey{1000}\truez{600}%
\put(\value{x},\value{x}){\makebox(0,\value{z})[l]{${#1}$}}%
\put(-\value{x},-\value{y}){\makebox(0,\value{z})[r]{${#2}$}}%
\end{picture}}%

\newcommand{\ESEADJDIST}[3]{%
\Y=#3%
\divide\Y by 2%
\Z=\Y%
\divide\Z by 2%
\begin{picture}(0,0)%
\truex{156}\truey{313}\truez{400}%
\put(-\Y,\Z){\begin{picture}(0,0)%
\monolength=#3%
\advance\monolength by -\value{x}%
\epilength=#3%
\advance\epilength by \value{x}%
\secondmonolength=\Y%
\advance\secondmonolength by -\value{y}%
\secondepilength=\Y%
\advance\secondepilength by \value{y}%
\put(-\value{x},-\value{y}){\line(2,-1){#3}}%
\put(\monolength,-\secondepilength){\esehead}%
\put(\epilength,-\secondmonolength){\line(-2,1){#3}}%
\put(\value{x},\value{y}){\wnwhead}%
\end{picture}}
\put(\value{x},\value{y}){\distsign}%
\put(-\value{x},-\value{y}){\distsign}%
\truex{400}\truey{1000}\truez{600}%
\put(\value{x},\value{x}){\makebox(0,\value{z})[l]{${#1}$}}%
\put(-\value{x},-\value{y}){\makebox(0,\value{z})[r]{${#2}$}}%
\end{picture}}%

\def\basicesear[#1]{\ESEAR{}{}{#100}}%

\newcommand{\esear}{\@ifnextchar[{\basicesear}{\basicesear[133]}}%

\def\basicEsear[#1]#2{\ESEAR{#2}{}{#100}}%

\newcommand{\Esear}{\@ifnextchar[{\basicEsear}{\basicEsear[133]}}%

\def\basiceseaR[#1]#2{\ESEAR{}{#2}{#100}}%

\newcommand{\eseaR}{\@ifnextchar[{\basiceseaR}{\basiceseaR[133]}}%

\def\basicesedist[#1]{\ESEDIST{}{}{#100}}%

\newcommand{\esedist}{\@ifnextchar[{\basicesedist}{\basicesedist[133]}}%

\def\basicEsedist[#1]#2{\ESEDIST{#2}{}{#100}}%

\newcommand{\Esedist}{\@ifnextchar[{\basicEsedist}{\basicEsedist[133]}}%

\def\basicesedisT[#1]#2{\ESEDIST{}{#2}{#100}}%

\newcommand{\esedisT}{\@ifnextchar[{\basicesedisT}{\basicesedisT[133]}}%

\def\basicesedotar[#1]{\ESEDOTAR{}{}{#100}}%

\newcommand{\esedotar}{\@ifnextchar[{\basicesedotar}{\basicesedotar[133]}}%

\def\basicEsedotar[#1]#2{\ESEDOTAR{#2}{}{#100}}%

\newcommand{\Esedotar}{\@ifnextchar[{\basicEsedotar}{\basicEsedotar[133]}}%

\def\basicesedotaR[#1]#2{\ESEDOTAR{}{#2}{#100}}%

\newcommand{\esedotaR}{\@ifnextchar[{\basicesedotaR}{\basicesedotaR[133]}}%

\def\basicesemono[#1]{\ESEMONO{}{}{#100}}%

\newcommand{\esemono}{\@ifnextchar[{\basicesemono}{\basicesemono[133]}}%

\def\basicEsemono[#1]#2{\ESEMONO{#2}{}{#100}}%

\newcommand{\Esemono}{\@ifnextchar[{\basicEsemono}{\basicEsemono[133]}}%

\def\basicesemonO[#1]#2{\ESEMONO{}{#2}{#100}}%

\newcommand{\esemonO}{\@ifnextchar[{\basicesemonO}{\basicesemonO[133]}}%

\def\basiceseepi[#1]{\ESEEPI{}{}{#100}}%

\newcommand{\eseepi}{\@ifnextchar[{\basiceseepi}{\basiceseepi[133]}}%

\def\basicEseepi[#1]#2{\ESEEPI{#2}{}{#100}}%

\newcommand{\Eseepi}{\@ifnextchar[{\basicEseepi}{\basicEseepi[133]}}%

\def\basiceseepI[#1]#2{\ESEEPI{}{#2}{#100}}%

\newcommand{\eseepI}{\@ifnextchar[{\basiceseepI}{\basiceseepI[133]}}%

\def\basicesebimo[#1]{\ESEBIMO{}{}{#100}}%

\newcommand{\esebimo}{\@ifnextchar[{\basicesebimo}{\basicesebimo[133]}}%

\def\basicEsebimo[#1]#2{\ESEBIMO{#2}{}{#100}}%

\newcommand{\Esebimo}{\@ifnextchar[{\basicEsebimo}{\basicEsebimo[133]}}%

\def\basicesebimO[#1]#2{\ESEBIMO{}{#2}{#100}}%

\newcommand{\esebimO}{\@ifnextchar[{\basicesebimO}{\basicesebimO[133]}}%

\def\basiceseiso[#1]{\ESEAR{\isosign}{}{#100}}%

\newcommand{\eseiso}{\@ifnextchar[{\basiceseiso}{\basiceseiso[133]}}%

\def\basicEseiso[#1]#2{\ESEAR{#2}{\isosign}{#100}}%

\newcommand{\Eseiso}{\@ifnextchar[{\basicEseiso}{\basicEseiso[133]}}%

\def\basiceseisO[#1]#2{\ESEAR{\isosign}{#2}{#100}}%

\newcommand{\eseisO}{\@ifnextchar[{\basiceseisO}{\basiceseisO[133]}}%

\def\basiceseeql[#1]{\ESEEQL{}{}{#100}}%

\newcommand{\eseeql}{\@ifnextchar[{\basiceseeql}{\basiceseeql[133]}}%

\def\basicEseeql[#1]#2{\ESEEQL{#2}{}{#100}}%

\newcommand{\Eseeql}{\@ifnextchar[{\basicEseeql}{\basicEseeql[133]}}%

\def\basiceseeqL[#1]#2{\ESEEQL{}{#2}{#100}}%

\newcommand{\eseeqL}{\@ifnextchar[{\basiceseeqL}{\basiceseeqL[133]}}%

\def\basicesebiar[#1]{\ESEBIAR{}{}{#100}}%

\newcommand{\esebiar}{\@ifnextchar[{\basicesebiar}{\basicesebiar[133]}}%

\def\basicEsebiar[#1]#2#3{\ESEBIAR{#2}{#3}{#100}}%

\newcommand{\Esebiar}{\@ifnextchar[{\basicEsebiar}{\basicEsebiar[133]}}%

\def\basicesebidist[#1]{\ESEBIDIST{}{}{#100}}%

\newcommand{\esebidist}{\@ifnextchar[{\basicesebidist}{\basicesebidist[133]}}%

\def\basicEsebidist[#1]#2#3{\ESEBIDIST{#2}{#3}{#100}}%

\newcommand{\Esebidist}{\@ifnextchar[{\basicEsebidist}{\basicEsebidist[133]}}%

\def\basiceseadjar[#1]{\ESEADJAR{}{}{#100}}%

\newcommand{\eseadjar}{\@ifnextchar[{\basiceseadjar}{\basiceseadjar[133]}}%

\def\basicEseadjar[#1]#2#3{\ESEADJAR{#2}{#3}{#100}}%

\newcommand{\Eseadjar}{\@ifnextchar[{\basicEseadjar}{\basicEseadjar[133]}}%

\def\basiceseadjdist[#1]{\ESEADJDIST{}{}{#100}}%

\newcommand{\eseadjdist}{\@ifnextchar[{\basiceseadjdist}{\basiceseadjdist[13
3]}}%

\def\basicEseadjdist[#1]#2#3{\ESEADJDIST{#2}{#3}{#100}}%

\newcommand{\Eseadjdist}{\@ifnextchar[{\basicEseadjdist}{\basicEseadjdist[13
3]}}%

\newcommand{\WSWAR}[3]{%
\Y=#3%
\divide\Y by 2%
\Z=\Y%
\divide\Z by 2%
\begin{picture}(0,0)%
\put(\Y,\Z){\line(-2,-1){#3}}%
\put(-\Y,-\Z){\wswhead}%
\truex{200}\truey{800}\truez{600}%
\put(-\value{x},\value{x}){\makebox(0,\value{z})[r]{${#1}$}}%
\put(\value{x},-\value{y}){\makebox(0,\value{z})[l]{${#2}$}}%
\end{picture}}%

\newcommand{\WSWDIST}[3]{%
\Y=#3%
\divide\Y by 2%
\Z=\Y%
\divide\Z by 2%
\begin{picture}(0,0)%
\put(\Y,\Z){\line(-2,-1){#3}}%
\put(-\Y,-\Z){\wswhead}%
\truex{400}%
\put(0,0){\distsign}%
\truex{200}\truey{800}\truez{600}%
\put(-\value{x},\value{x}){\makebox(0,\value{z})[r]{${#1}$}}%
\put(\value{x},-\value{y}){\makebox(0,\value{z})[l]{${#2}$}}%
\end{picture}}%

\newcommand{\WSWDOTAR}[3]{%
\truex{100}\truey{268}\truez{134}%
\Y=#3%
\divide\Y by 2%
\Z=\Y%
\divide\Z by 2%
\NUMBEROFDOTS=#3%
\divide\NUMBEROFDOTS by \value{y}%
\advance\NUMBEROFDOTS by 1%
\begin{picture}(0,0)%
\multiput(\Y,\Z)(-\value{y},-\value{z}){\NUMBEROFDOTS}%
{\circle*{\value{x}}}%
\put(-\Y,-\Z){\wswhead}%
\truex{200}\truey{800}\truez{600}%
\put(-\value{x},\value{x}){\makebox(0,\value{z})[r]{${#1}$}}%
\put(\value{x},-\value{y}){\makebox(0,\value{z})[l]{${#2}$}}%
\end{picture}}%

\newcommand{\WSWMONO}[3]{%
\Y=#3%
\divide\Y by 2%
\Z=\Y%
\divide\Z by 2%
\TrueTail%
\bimolength=#3%
\advance\bimolength by -\TrueMonoTail%
\monolength=\bimolength%
\advance\monolength by -\Y%
\secondmonolength=\monolength%
\divide\secondmonolength by 2%
\begin{picture}(0,0)%
\put(\monolength,\secondmonolength){\line(-2,-1){\bimolength}}%
\put(\monolength,\secondmonolength){\wswhead}%
\put(-\Y,-\Z){\wswhead}%
\truex{200}\truey{800}\truez{600}%
\put(-\value{x},\value{x}){\makebox(0,\value{z})[r]{${#1}$}}%
\put(\value{x},-\value{y}){\makebox(0,\value{z})[l]{${#2}$}}%
\end{picture}}%

\newcommand{\WSWEPI}[3]{%
\Y=#3%
\divide\Y by 2%
\Z=\Y%
\divide\Z by 2%
\TrueHead%
\bimolength=#3%
\advance\bimolength by -\TrueEpiHead%
\epilength=\bimolength%
\advance\epilength by -\Y%
\secondepilength=\epilength%
\divide\secondepilength by 2%
\begin{picture}(0,0)%
\put(\Y,\Z){\line(-2,-1){\bimolength}}%
\put(-\epilength,-\secondepilength){\wswhead}%
\put(-\Y,-\Z){\wswhead}%
\truex{200}\truey{800}\truez{600}%
\put(-\value{x},\value{x}){\makebox(0,\value{z})[r]{${#1}$}}%
\put(\value{x},-\value{y}){\makebox(0,\value{z})[l]{${#2}$}}%
\end{picture}}%

\newcommand{\WSWBIMO}[3]{%
\Y=#3%
\divide\Y by 2%
\Z=\Y%
\divide\Z by 2%
\TrueTail\TrueHead%
\bimolength=#3%
\advance\bimolength by -\TrueMonoTail%
\monolength=\bimolength%
\advance\monolength by -\Y%
\advance\bimolength by -\TrueEpiHead%
\epilength=\bimolength%
\advance\epilength by -\monolength%
\secondmonolength=\monolength%
\divide\secondmonolength by 2%
\secondepilength=\epilength%
\divide\secondepilength by 2%
\begin{picture}(0,0)%
\put(\monolength,\secondmonolength){\line(-2,-1){\bimolength}}%
\put(\monolength,\secondmonolength){\wswhead}%
\put(-\epilength,-\secondepilength){\wswhead}%
\put(-\Y,-\Z){\wswhead}%
\truex{200}\truey{800}\truez{600}%
\put(-\value{x},\value{x}){\makebox(0,\value{z})[r]{${#1}$}}%
\put(\value{x},-\value{y}){\makebox(0,\value{z})[l]{${#2}$}}%
\end{picture}}%

\newcommand{\WSWBIAR}[3]{%
\Y=#3%
\divide\Y by 2%
\Z=\Y%
\divide\Z by 2%
\begin{picture}(0,0)%
\put(\Y,\Z){\begin{picture}(0,0)%
\truex{156}\truey{313}%
\put(-\value{x},\value{y}){\line(-2,-1){#3}}%
\put(\value{x},-\value{y}){\line(-2,-1){#3}}%
\monolength=#3%
\advance\monolength by -\value{x}%
\epilength=#3%
\advance\epilength by \value{x}%
\secondmonolength=\Y%
\advance\secondmonolength by -\value{y}%
\secondepilength=\Y%
\advance\secondepilength by \value{y}%
\put(-\monolength,-\secondepilength){\wswhead}%
\put(-\epilength,-\secondmonolength){\wswhead}%
\end{picture}}
\truex{300}\truey{1000}\truez{600}%
\put(-\value{x},\value{x}){\makebox(0,\value{z})[r]{${#1}$}}%
\put(\value{x},-\value{y}){\makebox(0,\value{z})[l]{${#2}$}}%
\end{picture}}%

\newcommand{\WSWBIDIST}[3]{%
\Y=#3%
\divide\Y by 2%
\Z=\Y%
\divide\Z by 2%
\begin{picture}(0,0)%
\truex{156}\truey{313}\truez{400}%
\put(\Y,\Z){\begin{picture}(0,0)%
\put(-\value{x},\value{y}){\line(-2,-1){#3}}%
\put(\value{x},-\value{y}){\line(-2,-1){#3}}%
\monolength=#3%
\advance\monolength by -\value{x}%
\epilength=#3%
\advance\epilength by \value{x}%
\secondmonolength=\Y%
\advance\secondmonolength by -\value{y}%
\secondepilength=\Y%
\advance\secondepilength by \value{y}%
\put(-\monolength,-\secondepilength){\wswhead}%
\put(-\epilength,-\secondmonolength){\wswhead}%
\end{picture}}
\put(-\value{x},\value{y}){\distsign}%
\put(\value{x},-\value{y}){\distsign}%
\truex{300}\truey{1000}\truez{600}%
\put(-\value{x},\value{x}){\makebox(0,\value{z})[r]{${#1}$}}%
\put(\value{x},-\value{y}){\makebox(0,\value{z})[l]{${#2}$}}%
\end{picture}}%

\newcommand{\WSWADJAR}[3]{%
\Y=#3%
\divide\Y by 2%
\Z=\Y%
\divide\Z by 2%
\begin{picture}(0,0)%
\put(\Y,\Z){\begin{picture}(0,0)%
\truex{156}\truey{313}%
\monolength=#3%
\advance\monolength by -\value{x}%
\epilength=#3%
\advance\epilength by \value{x}%
\secondmonolength=\Y%
\advance\secondmonolength by -\value{y}%
\secondepilength=\Y%
\advance\secondepilength by \value{y}%
\put(\value{x},-\value{y}){\line(-2,-1){#3}}%
\put(-\monolength,-\secondepilength){\wswhead}%
\put(-\epilength,-\secondmonolength){\line(2,1){#3}}%
\put(-\value{x},\value{y}){\enehead}%
\end{picture}}
\truex{300}\truey{1000}\truez{600}%
\put(-\value{x},\value{x}){\makebox(0,\value{z})[r]{${#1}$}}%
\put(\value{x},-\value{y}){\makebox(0,\value{z})[l]{${#2}$}}%
\end{picture}}%

\newcommand{\WSWADJDIST}[3]{%
\Y=#3%
\divide\Y by 2%
\Z=\Y%
\divide\Z by 2%
\begin{picture}(0,0)%
\truex{156}\truey{313}\truez{400}%
\put(\Y,\Z){\begin{picture}(0,0)%
\monolength=#3%
\advance\monolength by -\value{x}%
\epilength=#3%
\advance\epilength by \value{x}%
\secondmonolength=\Y%
\advance\secondmonolength by -\value{y}%
\secondepilength=\Y%
\advance\secondepilength by \value{y}%
\put(\value{x},-\value{y}){\line(-2,-1){#3}}%
\put(-\monolength,-\secondepilength){\wswhead}%
\put(-\epilength,-\secondmonolength){\line(2,1){#3}}%
\put(-\value{x},\value{y}){\enehead}%
\end{picture}}
\put(-\value{x},\value{y}){\distsign}%
\put(\value{x},-\value{y}){\distsign}%
\truex{300}\truey{1000}\truez{600}%
\put(-\value{x},\value{x}){\makebox(0,\value{z})[r]{${#1}$}}%
\put(\value{x},-\value{y}){\makebox(0,\value{z})[l]{${#2}$}}%
\end{picture}}%

\def\basicwswar[#1]{\WSWAR{}{}{#100}}%

\newcommand{\wswar}{\@ifnextchar[{\basicwswar}{\basicwswar[133]}}%

\def\basicWswar[#1]#2{\WSWAR{#2}{}{#100}}%

\newcommand{\Wswar}{\@ifnextchar[{\basicWswar}{\basicWswar[133]}}%

\def\basicwswaR[#1]#2{\WSWAR{}{#2}{#100}}%

\newcommand{\wswaR}{\@ifnextchar[{\basicwswaR}{\basicwswaR[133]}}%

\def\basicwswdist[#1]{\WSWDIST{}{}{#100}}%

\newcommand{\wswdist}{\@ifnextchar[{\basicwswdist}{\basicwswdist[133]}}%

\def\basicWswdist[#1]#2{\WSWDIST{#2}{}{#100}}%

\newcommand{\Wswdist}{\@ifnextchar[{\basicWswdist}{\basicWswdist[133]}}%

\def\basicwswdisT[#1]#2{\WSWDIST{}{#2}{#100}}%

\newcommand{\wswdisT}{\@ifnextchar[{\basicwswdisT}{\basicwswdisT[133]}}%

\def\basicwswdotar[#1]{\WSWDOTAR{}{}{#100}}%

\newcommand{\wswdotar}{\@ifnextchar[{\basicwswdotar}{\basicwswdotar[133]}}%

\def\basicWswdotar[#1]#2{\WSWDOTAR{#2}{}{#100}}%

\newcommand{\Wswdotar}{\@ifnextchar[{\basicWswdotar}{\basicWswdotar[133]}}%

\def\basicwswdotaR[#1]#2{\WSWDOTAR{}{#2}{#100}}%

\newcommand{\wswdotaR}{\@ifnextchar[{\basicwswdotaR}{\basicwswdotaR[133]}}%

\def\basicwswmono[#1]{\WSWMONO{}{}{#100}}%

\newcommand{\wswmono}{\@ifnextchar[{\basicwswmono}{\basicwswmono[133]}}%

\def\basicWswmono[#1]#2{\WSWMONO{#2}{}{#100}}%

\newcommand{\Wswmono}{\@ifnextchar[{\basicWswmono}{\basicWswmono[133]}}%

\def\basicwswmonO[#1]#2{\WSWMONO{}{#2}{#100}}%

\newcommand{\wswmonO}{\@ifnextchar[{\basicwswmonO}{\basicwswmonO[133]}}%

\def\basicwswepi[#1]{\WSWEPI{}{}{#100}}%

\newcommand{\wswepi}{\@ifnextchar[{\basicwswepi}{\basicwswepi[133]}}%

\def\basicWswepi[#1]#2{\WSWEPI{#2}{}{#100}}%

\newcommand{\Wswepi}{\@ifnextchar[{\basicWswepi}{\basicWswepi[133]}}%

\def\basicwswepI[#1]#2{\WSWEPI{}{#2}{#100}}%

\newcommand{\wswepI}{\@ifnextchar[{\basicwswepI}{\basicwswepI[133]}}%

\def\basicwswbimo[#1]{\WSWBIMO{}{}{#100}}%

\newcommand{\wswbimo}{\@ifnextchar[{\basicwswbimo}{\basicwswbimo[133]}}%

\def\basicWswbimo[#1]#2{\WSWBIMO{#2}{}{#100}}%

\newcommand{\Wswbimo}{\@ifnextchar[{\basicWswbimo}{\basicWswbimo[133]}}%

\def\basicwswbimO[#1]#2{\WSWBIMO{}{#2}{#100}}%

\newcommand{\wswbimO}{\@ifnextchar[{\basicwswbimO}{\basicwswbimO[133]}}%

\def\basicwswiso[#1]{\WSWAR{\isosign}{}{#100}}%

\newcommand{\wswiso}{\@ifnextchar[{\basicwswiso}{\basicwswiso[133]}}%

\def\basicWswiso[#1]#2{\WSWAR{#2}{\isosign}{#100}}%

\newcommand{\Wswiso}{\@ifnextchar[{\basicWswiso}{\basicWswiso[133]}}%

\def\basicwswisO[#1]#2{\WSWAR{\isosign}{#2}{#100}}%

\newcommand{\wswisO}{\@ifnextchar[{\basicwswisO}{\basicwswisO[133]}}%

\def\basicwswbiar[#1]{\WSWBIAR{}{}{#100}}%

\newcommand{\wswbiar}{\@ifnextchar[{\basicwswbiar}{\basicwswbiar[133]}}%

\def\basicWswbiar[#1]#2#3{\WSWBIAR{#2}{#3}{#100}}%

\newcommand{\Wswbiar}{\@ifnextchar[{\basicWswbiar}{\basicWswbiar[133]}}%

\def\basicwswbidist[#1]{\WSWBIDIST{}{}{#100}}%

\newcommand{\wswbidist}{\@ifnextchar[{\basicwswbidist}{\basicwswbidist[133]}}%

\def\basicWswbidist[#1]#2#3{\WSWBIDIST{#2}{#3}{#100}}%

\newcommand{\Wswbidist}{\@ifnextchar[{\basicWswbidist}{\basicWswbidist[133]}}%

\def\basicwswadjar[#1]{\WSWADJAR{}{}{#100}}%

\newcommand{\wswadjar}{\@ifnextchar[{\basicwswadjar}{\basicwswadjar[133]}}%

\def\basicWswadjar[#1]#2#3{\WSWADJAR{#2}{#3}{#100}}%

\newcommand{\Wswadjar}{\@ifnextchar[{\basicWswadjar}{\basicWswadjar[133]}}%

\def\basicwswadjdist[#1]{\WSWADJDIST{}{}{#100}}%

\newcommand{\wswadjdist}{\@ifnextchar[{\basicwswadjdist}{\basicwswadjdist[13
3]}}%

\def\basicWswadjdist[#1]#2#3{\WSWADJDIST{#2}{#3}{#100}}%

\newcommand{\Wswadjdist}{\@ifnextchar[{\basicWswadjdist}{\basicWswadjdist[13
3]}}%

\newcommand{\WNWAR}[3]{%
\Y=#3%
\divide\Y by 2%
\Z=\Y%
\divide\Z by 2%
\begin{picture}(0,0)%
\put(\Y,-\Z){\line(-2,1){#3}}%
\put(-\Y,\Z){\wnwhead}%
\truex{200}\truey{800}\truez{600}%
\put(\value{x},\value{x}){\makebox(0,\value{z})[l]{${#1}$}}%
\put(-\value{x},-\value{y}){\makebox(0,\value{z})[r]{${#2}$}}%
\end{picture}}%

\newcommand{\WNWDIST}[3]{%
\Y=#3%
\divide\Y by 2%
\Z=\Y%
\divide\Z by 2%
\begin{picture}(0,0)%
\put(\Y,-\Z){\line(-2,1){#3}}%
\put(-\Y,\Z){\wnwhead}%
\truex{400}%
\put(0,0){\distsign}%
\truex{200}\truey{800}\truez{600}%
\put(\value{x},\value{x}){\makebox(0,\value{z})[l]{${#1}$}}%
\put(-\value{x},-\value{y}){\makebox(0,\value{z})[r]{${#2}$}}%
\end{picture}}%

\newcommand{\WNWDOTAR}[3]{%
\truex{100}\truey{268}\truez{134}%
\Y=#3%
\divide\Y by 2%
\Z=\Y%
\divide\Z by 2%
\NUMBEROFDOTS=#3%
\divide\NUMBEROFDOTS by \value{y}%
\advance\NUMBEROFDOTS by 1%
\begin{picture}(0,0)%
\multiput(\Y,-\Z)(-\value{y},\value{z}){\NUMBEROFDOTS}%
{\circle*{\value{x}}}%
\put(-\Y,\Z){\wnwhead}%
\truex{200}\truey{800}\truez{600}%
\put(\value{x},\value{x}){\makebox(0,\value{z})[l]{${#1}$}}%
\put(-\value{x},-\value{y}){\makebox(0,\value{z})[r]{${#2}$}}%
\end{picture}}%

\newcommand{\WNWMONO}[3]{%
\Y=#3%
\divide\Y by 2%
\Z=\Y%
\divide\Z by 2%
\TrueTail%
\bimolength=#3%
\advance\bimolength by -\TrueMonoTail%
\monolength=\bimolength%
\advance\monolength by -\Y%
\secondmonolength=\monolength%
\divide\secondmonolength by 2%
\begin{picture}(0,0)%
\put(\monolength,-\secondmonolength){\line(-2,1){\bimolength}}%
\put(\monolength,-\secondmonolength){\wnwhead}%
\put(-\Y,\Z){\wnwhead}%
\truex{200}\truey{800}\truez{600}%
\put(\value{x},\value{x}){\makebox(0,\value{z})[l]{${#1}$}}%
\put(-\value{x},-\value{y}){\makebox(0,\value{z})[r]{${#2}$}}%
\end{picture}}%

\newcommand{\WNWEPI}[3]{%
\Y=#3%
\divide\Y by 2%
\Z=\Y%
\divide\Z by 2%
\TrueHead%
\bimolength=#3%
\advance\bimolength by -\TrueEpiHead%
\epilength=\bimolength%
\advance\epilength by -\Y%
\secondepilength=\epilength%
\divide\secondepilength by 2%
\begin{picture}(0,0)%
\put(\Y,-\Z){\line(-2,1){\bimolength}}%
\put(-\epilength,\secondepilength){\wnwhead}%
\put(-\Y,\Z){\wnwhead}%
\truex{200}\truey{800}\truez{600}%
\put(\value{x},\value{x}){\makebox(0,\value{z})[l]{${#1}$}}%
\put(-\value{x},-\value{y}){\makebox(0,\value{z})[r]{${#2}$}}%
\end{picture}}%

\newcommand{\WNWBIMO}[3]{%
\Y=#3%
\divide\Y by 2%
\Z=\Y%
\divide\Z by 2%
\TrueTail\TrueHead%
\bimolength=#3%
\advance\bimolength by -\TrueMonoTail%
\monolength=\bimolength%
\advance\monolength by -\Y%
\advance\bimolength by -\TrueEpiHead%
\epilength=\bimolength%
\advance\epilength by -\monolength%
\secondmonolength=\monolength%
\divide\secondmonolength by 2%
\secondepilength=\epilength%
\divide\secondepilength by 2%
\begin{picture}(0,0)%
\put(\monolength,-\secondmonolength){\line(-2,1){\bimolength}}%
\put(\monolength,-\secondmonolength){\wnwhead}%
\put(-\epilength,\secondepilength){\wnwhead}%
\put(-\Y,\Z){\wnwhead}%
\truex{200}\truey{800}\truez{600}%
\put(\value{x},\value{x}){\makebox(0,\value{z})[l]{${#1}$}}%
\put(-\value{x},-\value{y}){\makebox(0,\value{z})[r]{${#2}$}}%
\end{picture}}%

\newcommand{\WNWBIAR}[3]{%
\Y=#3%
\divide\Y by 2%
\Z=\Y%
\divide\Z by 2%
\begin{picture}(0,0)%
\put(\Y,-\Z){\begin{picture}(0,0)%
\truex{156}\truey{313}%
\put(-\value{x},-\value{y}){\line(-2,1){#3}}%
\put(\value{x},\value{y}){\line(-2,1){#3}}%
\monolength=#3%
\advance\monolength by -\value{x}%
\epilength=#3%
\advance\epilength by \value{x}%
\secondmonolength=\Y%
\advance\secondmonolength by -\value{y}%
\secondepilength=\Y%
\advance\secondepilength by \value{y}%
\put(-\monolength,\secondepilength){\wnwhead}%
\put(-\epilength,\secondmonolength){\wnwhead}%
\end{picture}}
\truex{400}\truey{1000}\truez{600}%
\put(\value{x},\value{x}){\makebox(0,\value{z})[l]{${#1}$}}%
\put(-\value{x},-\value{y}){\makebox(0,\value{z})[r]{${#2}$}}%
\end{picture}}%

\newcommand{\WNWBIDIST}[3]{%
\Y=#3%
\divide\Y by 2%
\Z=\Y%
\divide\Z by 2%
\begin{picture}(0,0)%
\truex{156}\truey{313}\truez{400}%
\put(\Y,-\Z){\begin{picture}(0,0)%
\put(-\value{x},-\value{y}){\line(-2,1){#3}}%
\put(\value{x},\value{y}){\line(-2,1){#3}}%
\monolength=#3%
\advance\monolength by -\value{x}%
\epilength=#3%
\advance\epilength by \value{x}%
\secondmonolength=\Y%
\advance\secondmonolength by -\value{y}%
\secondepilength=\Y%
\advance\secondepilength by \value{y}%
\put(-\monolength,\secondepilength){\wnwhead}%
\put(-\epilength,\secondmonolength){\wnwhead}%
\end{picture}}
\put(\value{x},\value{y}){\distsign}%
\put(-\value{x},-\value{y}){\distsign}%
\truex{400}\truey{1000}\truez{600}%
\put(\value{x},\value{x}){\makebox(0,\value{z})[l]{${#1}$}}%
\put(-\value{x},-\value{y}){\makebox(0,\value{z})[r]{${#2}$}}%
\end{picture}}%

\newcommand{\WNWADJAR}[3]{%
\Y=#3%
\divide\Y by 2%
\Z=\Y%
\divide\Z by 2%
\begin{picture}(0,0)%
\put(\Y,-\Z){\begin{picture}(0,0)%
\truex{156}\truey{313}%
\monolength=#3%
\advance\monolength by -\value{x}%
\epilength=#3%
\advance\epilength by \value{x}%
\secondmonolength=\Y%
\advance\secondmonolength by -\value{y}%
\secondepilength=\Y%
\advance\secondepilength by \value{y}%
\put(-\value{x},-\value{y}){\line(-2,1){#3}}%
\put(-\epilength,\secondmonolength){\wnwhead}%
\put(-\monolength,\secondepilength){\line(2,-1){#3}}%
\put(\value{x},\value{y}){\esehead}%
\end{picture}}
\truex{400}\truey{1000}\truez{600}%
\put(\value{x},\value{x}){\makebox(0,\value{z})[l]{${#1}$}}%
\put(-\value{x},-\value{y}){\makebox(0,\value{z})[r]{${#2}$}}%
\end{picture}}%

\newcommand{\WNWADJDIST}[3]{%
\Y=#3%
\divide\Y by 2%
\Z=\Y%
\divide\Z by 2%
\begin{picture}(0,0)%
\truex{156}\truey{313}\truez{400}%
\put(\Y,-\Z){\begin{picture}(0,0)%
\monolength=#3%
\advance\monolength by -\value{x}%
\epilength=#3%
\advance\epilength by \value{x}%
\secondmonolength=\Y%
\advance\secondmonolength by -\value{y}%
\secondepilength=\Y%
\advance\secondepilength by \value{y}%
\put(-\value{x},-\value{y}){\line(-2,1){#3}}%
\put(-\epilength,\secondmonolength){\wnwhead}%
\put(-\monolength,\secondepilength){\line(2,-1){#3}}%
\put(\value{x},\value{y}){\esehead}%
\end{picture}}
\put(\value{x},\value{y}){\distsign}%
\put(-\value{x},-\value{y}){\distsign}%
\truex{400}\truey{1000}\truez{600}%
\put(\value{x},\value{x}){\makebox(0,\value{z})[l]{${#1}$}}%
\put(-\value{x},-\value{y}){\makebox(0,\value{z})[r]{${#2}$}}%
\end{picture}}%

\def\basicwnwar[#1]{\WNWAR{}{}{#100}}%

\newcommand{\wnwar}{\@ifnextchar[{\basicwnwar}{\basicwnwar[133]}}%

\def\basicWnwar[#1]#2{\WNWAR{#2}{}{#100}}%

\newcommand{\Wnwar}{\@ifnextchar[{\basicWnwar}{\basicWnwar[133]}}%

\def\basicwnwaR[#1]#2{\WNWAR{}{#2}{#100}}%

\newcommand{\wnwaR}{\@ifnextchar[{\basicwnwaR}{\basicwnwaR[133]}}%

\def\basicwnwdist[#1]{\WNWDIST{}{}{#100}}%

\newcommand{\wnwdist}{\@ifnextchar[{\basicwnwdist}{\basicwnwdist[133]}}%

\def\basicWnwdist[#1]#2{\WNWDIST{#2}{}{#100}}%

\newcommand{\Wnwdist}{\@ifnextchar[{\basicWnwdist}{\basicWnwdist[133]}}%

\def\basicwnwdisT[#1]#2{\WNWDIST{}{#2}{#100}}%

\newcommand{\wnwdisT}{\@ifnextchar[{\basicwnwdisT}{\basicwnwdisT[133]}}%

\def\basicwnwdotar[#1]{\WNWDOTAR{}{}{#100}}%

\newcommand{\wnwdotar}{\@ifnextchar[{\basicwnwdotar}{\basicwnwdotar[133]}}%

\def\basicWnwdotar[#1]#2{\WNWDOTAR{#2}{}{#100}}%

\newcommand{\Wnwdotar}{\@ifnextchar[{\basicWnwdotar}{\basicWnwdotar[133]}}%

\def\basicwnwdotaR[#1]#2{\WNWDOTAR{}{#2}{#100}}%

\newcommand{\wnwdotaR}{\@ifnextchar[{\basicwnwdotaR}{\basicwnwdotaR[133]}}%

\def\basicwnwmono[#1]{\WNWMONO{}{}{#100}}%

\newcommand{\wnwmono}{\@ifnextchar[{\basicwnwmono}{\basicwnwmono[133]}}%

\def\basicWnwmono[#1]#2{\WNWMONO{#2}{}{#100}}%

\newcommand{\Wnwmono}{\@ifnextchar[{\basicWnwmono}{\basicWnwmono[133]}}%

\def\basicwnwmonO[#1]#2{\WNWMONO{}{#2}{#100}}%

\newcommand{\wnwmonO}{\@ifnextchar[{\basicwnwmonO}{\basicwnwmonO[133]}}%

\def\basicwnwepi[#1]{\WNWEPI{}{}{#100}}%

\newcommand{\wnwepi}{\@ifnextchar[{\basicwnwepi}{\basicwnwepi[133]}}%

\def\basicWnwepi[#1]#2{\WNWEPI{#2}{}{#100}}%

\newcommand{\Wnwepi}{\@ifnextchar[{\basicWnwepi}{\basicWnwepi[133]}}%

\def\basicwnwepI[#1]#2{\WNWEPI{}{#2}{#100}}%

\newcommand{\wnwepI}{\@ifnextchar[{\basicwnwepI}{\basicwnwepI[133]}}%

\def\basicwnwbimo[#1]{\WNWBIMO{}{}{#100}}%

\newcommand{\wnwbimo}{\@ifnextchar[{\basicwnwbimo}{\basicwnwbimo[133]}}%

\def\basicWnwbimo[#1]#2{\WNWBIMO{#2}{}{#100}}%

\newcommand{\Wnwbimo}{\@ifnextchar[{\basicWnwbimo}{\basicWnwbimo[133]}}%

\def\basicwnwbimO[#1]#2{\WNWBIMO{}{#2}{#100}}%

\newcommand{\wnwbimO}{\@ifnextchar[{\basicwnwbimO}{\basicwnwbimO[133]}}%

\def\basicwnwiso[#1]{\WNWAR{\isosign}{}{#100}}%

\newcommand{\wnwiso}{\@ifnextchar[{\basicwnwiso}{\basicwnwiso[133]}}%

\def\basicWnwiso[#1]#2{\WNWAR{#2}{\isosign}{#100}}%

\newcommand{\Wnwiso}{\@ifnextchar[{\basicWnwiso}{\basicWnwiso[133]}}%

\def\basicwnwisO[#1]#2{\WNWAR{\isosign}{#2}{#100}}%

\newcommand{\wnwisO}{\@ifnextchar[{\basicwnwisO}{\basicwnwisO[133]}}%

\def\basicwnwbiar[#1]{\WNWBIAR{}{}{#100}}%

\newcommand{\wnwbiar}{\@ifnextchar[{\basicwnwbiar}{\basicwnwbiar[133]}}%

\def\basicWnwbiar[#1]#2#3{\WNWBIAR{#2}{#3}{#100}}%

\newcommand{\Wnwbiar}{\@ifnextchar[{\basicWnwbiar}{\basicWnwbiar[133]}}%

\def\basicwnwbidist[#1]{\WNWBIDIST{}{}{#100}}%

\newcommand{\wnwbidist}{\@ifnextchar[{\basicwnwbidist}{\basicwnwbidist[133]}}%

\def\basicWnwbidist[#1]#2#3{\WNWBIDIST{#2}{#3}{#100}}%

\newcommand{\Wnwbidist}{\@ifnextchar[{\basicWnwbidist}{\basicWnwbidist[133]}}%

\def\basicwnwadjar[#1]{\WNWADJAR{}{}{#100}}%

\newcommand{\wnwadjar}{\@ifnextchar[{\basicwnwadjar}{\basicwnwadjar[133]}}%

\def\basicWnwadjar[#1]#2#3{\WNWADJAR{#2}{#3}{#100}}%

\newcommand{\Wnwadjar}{\@ifnextchar[{\basicWnwadjar}{\basicWnwadjar[133]}}%

\def\basicwnwadjdist[#1]{\WNWADJDIST{}{}{#100}}%

\newcommand{\wnwadjdist}{\@ifnextchar[{\basicwnwadjdist}{\basicwnwadjdist[13
3]}}%

\def\basicWnwadjdist[#1]#2#3{\WNWADJDIST{#2}{#3}{#100}}%

\newcommand{\Wnwadjdist}{\@ifnextchar[{\basicWnwadjdist}{\basicWnwadjdist[13
3]}}%

\newcommand{\NNEAR}[3]{%
\Z=#3%
\divide\Z by 2%
\begin{picture}(0,0)%
\put(-\Z,-#3){\line(1,2){#3}}%
\put(\Z,#3){\nnehead}%
\truex{200}\truey{800}\truez{600}%
\put(-\value{x},\value{x}){\makebox(0,\value{z})[r]{${#1}$}}%
\put(\value{x},-\value{y}){\makebox(0,\value{z})[l]{${#2}$}}%
\end{picture}}%

\newcommand{\NNEDIST}[3]{%
\Z=#3%
\divide\Z by 2%
\begin{picture}(0,0)%
\put(-\Z,-#3){\line(1,2){#3}}%
\put(\Z,#3){\nnehead}%
\truex{400}%
\put(0,0){\distsign}%
\truex{200}\truey{800}\truez{600}%
\put(-\value{x},\value{x}){\makebox(0,\value{z})[r]{${#1}$}}%
\put(\value{x},-\value{y}){\makebox(0,\value{z})[l]{${#2}$}}%
\end{picture}}%

\newcommand{\NNEDOTAR}[3]{%
\truex{100}\truey{268}\truez{134}%
\Z=#3%
\divide\Z by 2%
\NUMBEROFDOTS=#3%
\divide\NUMBEROFDOTS by \value{z}%
\advance\NUMBEROFDOTS by 1%
\begin{picture}(0,0)%
\multiput(-\Z,-#3)(\value{z},\value{y}){\NUMBEROFDOTS}%
{\circle*{\value{x}}}%
\put(\Z,#3){\nnehead}%
\truex{200}\truey{800}\truez{600}%
\put(-\value{x},\value{x}){\makebox(0,\value{z})[r]{${#1}$}}%
\put(\value{x},-\value{y}){\makebox(0,\value{z})[l]{${#2}$}}%
\end{picture}}%

\newcommand{\NNEMONO}[3]{%
\Z=#3%
\divide\Z by 2%
\truetaiL%
\bimolength=#3%
\advance\bimolength by -\truemonotaiL%
\monolength=\bimolength%
\advance\monolength by -\Z%
\secondmonolength=\monolength%
\multiply\secondmonolength by 2%
\begin{picture}(0,0)%
\put(-\monolength,-\secondmonolength){\line(1,2){\bimolength}}%
\put(-\monolength,-\secondmonolength){\nnehead}%
\put(\Z,#3){\nnehead}%
\truex{200}\truey{800}\truez{600}%
\put(-\value{x},\value{x}){\makebox(0,\value{z})[r]{${#1}$}}%
\put(\value{x},-\value{y}){\makebox(0,\value{z})[l]{${#2}$}}%
\end{picture}}%

\newcommand{\NNEEPI}[3]{%
\Z=#3%
\divide\Z by 2%
\trueheaD%
\bimolength=#3%
\advance\bimolength by -\trueepiheaD%
\epilength=\bimolength%
\advance\epilength by -\Z%
\secondepilength=\epilength%
\multiply\secondepilength by 2%
\begin{picture}(0,0)%
\put(-\Z,-#3){\line(1,2){\bimolength}}%
\put(\epilength,\secondepilength){\nnehead}%
\put(\Z,#3){\nnehead}%
\truex{200}\truey{800}\truez{600}%
\put(-\value{x},\value{x}){\makebox(0,\value{z})[r]{${#1}$}}%
\put(\value{x},-\value{y}){\makebox(0,\value{z})[l]{${#2}$}}%
\end{picture}}%

\newcommand{\NNEBIMO}[3]{%
\Z=#3%
\divide\Z by 2%
\truetaiL\trueheaD%
\bimolength=#3%
\advance\bimolength by -\truemonotaiL%
\monolength=\bimolength%
\advance\monolength by -\Z%
\advance\bimolength by -\trueepiheaD%
\epilength=\bimolength%
\advance\epilength by -\monolength%
\secondmonolength=\monolength%
\multiply\secondmonolength by 2%
\secondepilength=\epilength%
\multiply\secondepilength by 2%
\begin{picture}(0,0)%
\put(-\monolength,-\secondmonolength){\line(1,2){\bimolength}}%
\put(-\monolength,-\secondmonolength){\nnehead}%
\put(\epilength,\secondepilength){\nnehead}%
\put(\Z,#3){\nnehead}%
\truex{200}\truey{800}\truez{600}%
\put(-\value{x},\value{x}){\makebox(0,\value{z})[r]{${#1}$}}%
\put(\value{x},-\value{y}){\makebox(0,\value{z})[l]{${#2}$}}%
\end{picture}}%

\newcommand{\NNEEQL}[3]{%
\Z=#3%
\divide\Z by 2%
\begin{picture}(0,0)%
\put(-\Z,-#3){\begin{picture}(0,0)%
\truex{44}\truey{89}%
\put(-\value{y},\value{x}){\line(1,2){#3}}%
\put(\value{y},-\value{x}){\line(1,2){#3}}%
\end{picture}}%
\truex{200}\truey{800}\truez{600}%
\put(-\value{x},\value{x}){\makebox(0,\value{z})[r]{${#1}$}}%
\put(\value{x},-\value{y}){\makebox(0,\value{z})[l]{${#2}$}}%
\end{picture}}%

\newcommand{\NNEBIAR}[3]{%
\Y=#3%
\divide\Y by 2%
\Z=#3%
\multiply \Z by 2%
\begin{picture}(0,0)%
\put(-\Y,-#3){\begin{picture}(0,0)%
\truex{313}\truey{156}%
\put(-\value{x},\value{y}){\line(1,2){#3}}%
\put(\value{x},-\value{y}){\line(1,2){#3}}%
\monolength=#3%
\advance\monolength by -\value{x}%
\epilength=#3%
\advance\epilength by \value{x}%
\secondmonolength=\Z%
\advance\secondmonolength by -\value{y}%
\secondepilength=\Z%
\advance\secondepilength by \value{y}%
\put(\monolength,\secondepilength){\nnehead}%
\put(\epilength,\secondmonolength){\nnehead}%
\end{picture}}
\truex{300}\truey{1000}\truez{600}%
\put(-\value{x},\value{x}){\makebox(0,\value{z})[r]{${#1}$}}%
\put(\value{x},-\value{y}){\makebox(0,\value{z})[l]{${#2}$}}%
\end{picture}}%

\newcommand{\NNEBIDIST}[3]{%
\Y=#3%
\divide\Y by 2%
\Z=#3%
\multiply \Z by 2%
\begin{picture}(0,0)%
\truex{313}\truey{156}\truez{400}%
\put(-\Y,-#3){\begin{picture}(0,0)%
\put(-\value{x},\value{y}){\line(1,2){#3}}%
\put(\value{x},-\value{y}){\line(1,2){#3}}%
\monolength=#3%
\advance\monolength by -\value{x}%
\epilength=#3%
\advance\epilength by \value{x}%
\secondmonolength=\Z%
\advance\secondmonolength by -\value{y}%
\secondepilength=\Z%
\advance\secondepilength by \value{y}%
\put(\monolength,\secondepilength){\nnehead}%
\put(\epilength,\secondmonolength){\nnehead}%
\end{picture}}
\put(-\value{x},\value{y}){\distsign}%
\put(\value{x},-\value{y}){\distsign}%
\truex{300}\truey{1000}\truez{600}%
\put(-\value{x},\value{x}){\makebox(0,\value{z})[r]{${#1}$}}%
\put(\value{x},-\value{y}){\makebox(0,\value{z})[l]{${#2}$}}%
\end{picture}}%

\newcommand{\NNEADJAR}[3]{%
\Y=#3%
\divide\Y by 2%
\Z=#3%
\multiply \Z by 2%
\begin{picture}(0,0)%
\put(-\Y,-#3){\begin{picture}(0,0)%
\truex{313}\truey{156}%
\monolength=#3%
\advance\monolength by -\value{x}%
\epilength=#3%
\advance\epilength by \value{x}%
\secondmonolength=\Z%
\advance\secondmonolength by -\value{y}%
\secondepilength=\Z%
\advance\secondepilength by \value{y}%
\put(\value{x},-\value{y}){\line(1,2){#3}}%
\put(\epilength,\secondmonolength){\nnehead}%
\put(\monolength,\secondepilength){\line(-1,-2){#3}}%
\put(-\value{x},\value{y}){\sswhead}%
\end{picture}}
\truex{300}\truey{1000}\truez{600}%
\put(-\value{x},\value{x}){\makebox(0,\value{z})[r]{${#1}$}}%
\put(\value{x},-\value{y}){\makebox(0,\value{z})[l]{${#2}$}}%
\end{picture}}%

\newcommand{\NNEADJDIST}[3]{%
\Y=#3%
\divide\Y by 2%
\Z=#3%
\multiply \Z by 2%
\begin{picture}(0,0)%
\truex{313}\truey{156}\truez{400}%
\put(-\Y,-#3){\begin{picture}(0,0)%
\monolength=#3%
\advance\monolength by -\value{x}%
\epilength=#3%
\advance\epilength by \value{x}%
\secondmonolength=\Z%
\advance\secondmonolength by -\value{y}%
\secondepilength=\Z%
\advance\secondepilength by \value{y}%
\put(\value{x},-\value{y}){\line(1,2){#3}}%
\put(\epilength,\secondmonolength){\nnehead}%
\put(\monolength,\secondepilength){\line(-1,-2){#3}}%
\put(-\value{x},\value{y}){\sswhead}%
\end{picture}}
\put(-\value{x},\value{y}){\distsign}%
\put(\value{x},-\value{y}){\distsign}%
\truex{300}\truey{1000}\truez{600}%
\put(-\value{x},\value{x}){\makebox(0,\value{z})[r]{${#1}$}}%
\put(\value{x},-\value{y}){\makebox(0,\value{z})[l]{${#2}$}}%
\end{picture}}%

\def\basicnnear[#1]{\NNEAR{}{}{#100}}%

\newcommand{\nnear}{\@ifnextchar[{\basicnnear}{\basicnnear[67]}}%

\def\basicNnear[#1]#2{\NNEAR{#2}{}{#100}}%

\newcommand{\Nnear}{\@ifnextchar[{\basicNnear}{\basicNnear[67]}}%

\def\basicnneaR[#1]#2{\NNEAR{}{#2}{#100}}%

\newcommand{\nneaR}{\@ifnextchar[{\basicnneaR}{\basicnneaR[67]}}%

\def\basicnnedist[#1]{\NNEDIST{}{}{#100}}%

\newcommand{\nnedist}{\@ifnextchar[{\basicnnedist}{\basicnnedist[67]}}%

\def\basicNnedist[#1]#2{\NNEDIST{#2}{}{#100}}%

\newcommand{\Nnedist}{\@ifnextchar[{\basicNnedist}{\basicNnedist[67]}}%

\def\basicnnedisT[#1]#2{\NNEDIST{}{#2}{#100}}%

\newcommand{\nnedisT}{\@ifnextchar[{\basicnnedisT}{\basicnnedisT[67]}}%

\def\basicnnedotar[#1]{\NNEDOTAR{}{}{#100}}%

\newcommand{\nnedotar}{\@ifnextchar[{\basicnnedotar}{\basicnnedotar[67]}}%

\def\basicNnedotar[#1]#2{\NNEDOTAR{#2}{}{#100}}%

\newcommand{\Nnedotar}{\@ifnextchar[{\basicNnedotar}{\basicNnedotar[67]}}%

\def\basicnnedotaR[#1]#2{\NNEDOTAR{}{#2}{#100}}%

\newcommand{\nnedotaR}{\@ifnextchar[{\basicnnedotaR}{\basicnnedotaR[67]}}%

\def\basicnnemono[#1]{\NNEMONO{}{}{#100}}%

\newcommand{\nnemono}{\@ifnextchar[{\basicnnemono}{\basicnnemono[67]}}%

\def\basicNnemono[#1]#2{\NNEMONO{#2}{}{#100}}%

\newcommand{\Nnemono}{\@ifnextchar[{\basicNnemono}{\basicNnemono[67]}}%

\def\basicnnemonO[#1]#2{\NNEMONO{}{#2}{#100}}%

\newcommand{\nnemonO}{\@ifnextchar[{\basicnnemonO}{\basicnnemonO[67]}}%

\def\basicnneepi[#1]{\NNEEPI{}{}{#100}}%

\newcommand{\nneepi}{\@ifnextchar[{\basicnneepi}{\basicnneepi[67]}}%

\def\basicNneepi[#1]#2{\NNEEPI{#2}{}{#100}}%

\newcommand{\Nneepi}{\@ifnextchar[{\basicNneepi}{\basicNneepi[67]}}%

\def\basicnneepI[#1]#2{\NNEEPI{}{#2}{#100}}%

\newcommand{\nneepI}{\@ifnextchar[{\basicnneepI}{\basicnneepI[67]}}%

\def\basicnnebimo[#1]{\NNEBIMO{}{}{#100}}%

\newcommand{\nnebimo}{\@ifnextchar[{\basicnnebimo}{\basicnnebimo[67]}}%

\def\basicNnebimo[#1]#2{\NNEBIMO{#2}{}{#100}}%

\newcommand{\Nnebimo}{\@ifnextchar[{\basicNnebimo}{\basicNnebimo[67]}}%

\def\basicnnebimO[#1]#2{\NNEBIMO{}{#2}{#100}}%

\newcommand{\nnebimO}{\@ifnextchar[{\basicnnebimO}{\basicnnebimO[67]}}%

\def\basicnneiso[#1]{\NNEAR{\isosign}{}{#100}}%

\newcommand{\nneiso}{\@ifnextchar[{\basicnneiso}{\basicnneiso[67]}}%

\def\basicNneiso[#1]#2{\NNEAR{#2}{\isosign}{#100}}%

\newcommand{\Nneiso}{\@ifnextchar[{\basicNneiso}{\basicNneiso[67]}}%

\def\basicnneisO[#1]#2{\NNEAR{\isosign}{#2}{#100}}%

\newcommand{\nneisO}{\@ifnextchar[{\basicnneisO}{\basicnneisO[67]}}%

\def\basicnneeql[#1]{\NNEEQL{}{}{#100}}%

\newcommand{\nneeql}{\@ifnextchar[{\basicnneeql}{\basicnneeql[67]}}%

\def\basicNneeql[#1]#2{\NNEEQL{#2}{}{#100}}%

\newcommand{\Nneeql}{\@ifnextchar[{\basicNneeql}{\basicNneeql[67]}}%

\def\basicnneeqL[#1]#2{\NNEEQL{}{#2}{#100}}%

\newcommand{\nneeqL}{\@ifnextchar[{\basicnneeqL}{\basicnneeqL[67]}}%

\def\basicnnebiar[#1]{\NNEBIAR{}{}{#100}}%

\newcommand{\nnebiar}{\@ifnextchar[{\basicnnebiar}{\basicnnebiar[67]}}%

\def\basicNnebiar[#1]#2#3{\NNEBIAR{#2}{#3}{#100}}%

\newcommand{\Nnebiar}{\@ifnextchar[{\basicNnebiar}{\basicNnebiar[67]}}%

\def\basicnnebidist[#1]{\NNEBIDIST{}{}{#100}}%

\newcommand{\nnebidist}{\@ifnextchar[{\basicnnebidist}{\basicnnebidist[67]}}%

\def\basicNnebidist[#1]#2#3{\NNEBIDIST{#2}{#3}{#100}}%

\newcommand{\Nnebidist}{\@ifnextchar[{\basicNnebidist}{\basicNnebidist[67]}}%

\def\basicnneadjar[#1]{\NNEADJAR{}{}{#100}}%

\newcommand{\nneadjar}{\@ifnextchar[{\basicnneadjar}{\basicnneadjar[67]}}%

\def\basicNneadjar[#1]#2#3{\NNEADJAR{#2}{#3}{#100}}%

\newcommand{\Nneadjar}{\@ifnextchar[{\basicNneadjar}{\basicNneadjar[67]}}%

\def\basicnneadjdist[#1]{\NNEADJDIST{}{}{#100}}%

\newcommand{\nneadjdist}{\@ifnextchar[{\basicnneadjdist}{\basicnneadjdist[67]}}%

\def\basicNneadjdist[#1]#2#3{\NNEADJDIST{#2}{#3}{#100}}%

\newcommand{\Nneadjdist}{\@ifnextchar[{\basicNneadjdist}{\basicNneadjdist[67]}}%

\newcommand{\SSEAR}[3]{%
\Z=#3%
\divide\Z by 2%
\begin{picture}(0,0)%
\put(-\Z,#3){\line(1,-2){#3}}%
\put(\Z,-#3){\ssehead}%
\truex{200}\truey{800}\truez{600}%
\put(\value{x},\value{x}){\makebox(0,\value{z})[l]{${#1}$}}%
\put(-\value{x},-\value{y}){\makebox(0,\value{z})[r]{${#2}$}}%
\end{picture}}%

\newcommand{\SSEDIST}[3]{%
\Z=#3%
\divide\Z by 2%
\begin{picture}(0,0)%
\put(-\Z,#3){\line(1,-2){#3}}%
\put(\Z,-#3){\ssehead}%
\truex{400}%
\put(0,0){\distsign}%
\truex{200}\truey{800}\truez{600}%
\put(\value{x},\value{x}){\makebox(0,\value{z})[l]{${#1}$}}%
\put(-\value{x},-\value{y}){\makebox(0,\value{z})[r]{${#2}$}}%
\end{picture}}%

\newcommand{\SSEDOTAR}[3]{%
\truex{100}\truey{268}\truez{134}%
\Z=#3%
\divide\Z by 2%
\NUMBEROFDOTS=#3%
\divide\NUMBEROFDOTS by \value{z}%
\advance\NUMBEROFDOTS by 1%
\begin{picture}(0,0)%
\multiput(-\Z,#3)(\value{z},-\value{y}){\NUMBEROFDOTS}%
{\circle*{\value{x}}}%
\put(\Z,-#3){\ssehead}%
\truex{200}\truey{800}\truez{600}%
\put(\value{x},\value{x}){\makebox(0,\value{z})[l]{${#1}$}}%
\put(-\value{x},-\value{y}){\makebox(0,\value{z})[r]{${#2}$}}%
\end{picture}}%

\newcommand{\SSEMONO}[3]{%
\Z=#3%
\divide\Z by 2%
\truetaiL%
\bimolength=#3%
\advance\bimolength by -\truemonotaiL%
\monolength=\bimolength%
\advance\monolength by -\Z%
\secondmonolength=\monolength%
\multiply\secondmonolength by 2%
\begin{picture}(0,0)%
\put(-\monolength,\secondmonolength){\line(1,-2){\bimolength}}%
\put(-\monolength,\secondmonolength){\ssehead}%
\put(\Z,-#3){\ssehead}%
\truex{200}\truey{800}\truez{600}%
\put(\value{x},\value{x}){\makebox(0,\value{z})[l]{${#1}$}}%
\put(-\value{x},-\value{y}){\makebox(0,\value{z})[r]{${#2}$}}%
\end{picture}}%

\newcommand{\SSEEPI}[3]{%
\Z=#3%
\divide\Z by 2%
\trueheaD%
\bimolength=#3%
\advance\bimolength by -\trueepiheaD%
\epilength=\bimolength%
\advance\epilength by -\Z%
\secondepilength=\epilength%
\multiply\secondepilength by 2%
\begin{picture}(0,0)%
\put(-\Z,#3){\line(1,-2){\bimolength}}%
\put(\epilength,-\secondepilength){\ssehead}%
\put(\Z,-#3){\ssehead}%
\truex{200}\truey{800}\truez{600}%
\put(\value{x},\value{x}){\makebox(0,\value{z})[l]{${#1}$}}%
\put(-\value{x},-\value{y}){\makebox(0,\value{z})[r]{${#2}$}}%
\end{picture}}%

\newcommand{\SSEBIMO}[3]{%
\Z=#3%
\divide\Z by 2%
\truetaiL\trueheaD%
\bimolength=#3%
\advance\bimolength by -\truemonotaiL%
\monolength=\bimolength%
\advance\monolength by -\Z%
\advance\bimolength by -\trueepiheaD%
\epilength=\bimolength%
\advance\epilength by -\monolength%
\secondmonolength=\monolength%
\multiply\secondmonolength by 2%
\secondepilength=\epilength%
\multiply\secondepilength by 2%
\begin{picture}(0,0)%
\put(-\monolength,\secondmonolength){\line(1,-2){\bimolength}}%
\put(-\monolength,\secondmonolength){\ssehead}%
\put(\epilength,-\secondepilength){\ssehead}%
\put(\Z,-#3){\ssehead}%
\truex{200}\truey{800}\truez{600}%
\put(\value{x},\value{x}){\makebox(0,\value{z})[l]{${#1}$}}%
\put(-\value{x},-\value{y}){\makebox(0,\value{z})[r]{${#2}$}}%
\end{picture}}%

\newcommand{\SSEEQL}[3]{%
\Z=#3%
\divide\Z by 2%
\begin{picture}(0,0)%
\put(-\Z,#3){\begin{picture}(0,0)%
\truex{44}\truey{89}%
\put(-\value{y},-\value{x}){\line(1,-2){#3}}%
\put(\value{y},\value{x}){\line(1,-2){#3}}%
\end{picture}}%
\truex{200}\truey{800}\truez{600}%
\put(\value{x},\value{x}){\makebox(0,\value{z})[l]{${#1}$}}%
\put(-\value{x},-\value{y}){\makebox(0,\value{z})[r]{${#2}$}}%
\end{picture}}%

\newcommand{\SSEBIAR}[3]{%
\Y=#3%
\divide\Y by 2%
\Z=#3%
\multiply \Z by 2%
\begin{picture}(0,0)%
\put(-\Y,#3){\begin{picture}(0,0)%
\truex{313}\truey{156}%
\put(-\value{x},-\value{y}){\line(1,-2){#3}}%
\put(\value{x},\value{y}){\line(1,-2){#3}}%
\monolength=#3%
\advance\monolength by -\value{x}%
\epilength=#3%
\advance\epilength by \value{x}%
\secondmonolength=\Z%
\advance\secondmonolength by -\value{y}%
\secondepilength=\Z%
\advance\secondepilength by \value{y}%
\put(\monolength,-\secondepilength){\ssehead}%
\put(\epilength,-\secondmonolength){\ssehead}%
\end{picture}}
\truex{400}\truey{1000}\truez{600}%
\put(\value{x},\value{x}){\makebox(0,\value{z})[l]{${#1}$}}%
\put(-\value{x},-\value{y}){\makebox(0,\value{z})[r]{${#2}$}}%
\end{picture}}%

\newcommand{\SSEBIDIST}[3]{%
\Y=#3%
\divide\Y by 2%
\Z=#3%
\multiply \Z by 2%
\begin{picture}(0,0)%
\truex{313}\truey{156}\truez{400}%
\put(-\Y,#3){\begin{picture}(0,0)%
\put(-\value{x},-\value{y}){\line(1,-2){#3}}%
\put(\value{x},\value{y}){\line(1,-2){#3}}%
\monolength=#3%
\advance\monolength by -\value{x}%
\epilength=#3%
\advance\epilength by \value{x}%
\secondmonolength=\Z%
\advance\secondmonolength by -\value{y}%
\secondepilength=\Z%
\advance\secondepilength by \value{y}%
\put(\monolength,-\secondepilength){\ssehead}%
\put(\epilength,-\secondmonolength){\ssehead}%
\end{picture}}
\put(-\value{x},-\value{y}){\distsign}%
\put(\value{x},\value{y}){\distsign}%
\truex{500}\truey{1000}\truez{600}%
\put(\value{x},\value{x}){\makebox(0,\value{z})[l]{${#1}$}}%
\put(-\value{x},-\value{y}){\makebox(0,\value{z})[r]{${#2}$}}%
\end{picture}}%

\newcommand{\SSEADJAR}[3]{%
\Y=#3%
\divide\Y by 2%
\Z=#3%
\multiply \Z by 2%
\begin{picture}(0,0)%
\put(-\Y,#3){\begin{picture}(0,0)%
\truex{313}\truey{156}%
\monolength=#3%
\advance\monolength by -\value{x}%
\epilength=#3%
\advance\epilength by \value{x}%
\secondmonolength=\Z%
\advance\secondmonolength by -\value{y}%
\secondepilength=\Z%
\advance\secondepilength by \value{y}%
\put(-\value{x},-\value{y}){\line(1,-2){#3}}%
\put(\monolength,-\secondepilength){\ssehead}%
\put(\epilength,-\secondmonolength){\line(-1,2){#3}}%
\put(\value{x},\value{y}){\nnwhead}%
\end{picture}}
\truex{400}\truey{1000}\truez{600}%
\put(\value{x},\value{x}){\makebox(0,\value{z})[l]{${#1}$}}%
\put(-\value{x},-\value{y}){\makebox(0,\value{z})[r]{${#2}$}}%
\end{picture}}%

\newcommand{\SSEADJDIST}[3]{%
\Y=#3%
\divide\Y by 2%
\Z=#3%
\multiply \Z by 2%
\begin{picture}(0,0)%
\truex{313}\truey{156}\truez{400}%
\put(-\Y,#3){\begin{picture}(0,0)%
\monolength=#3%
\advance\monolength by -\value{x}%
\epilength=#3%
\advance\epilength by \value{x}%
\secondmonolength=\Z%
\advance\secondmonolength by -\value{y}%
\secondepilength=\Z%
\advance\secondepilength by \value{y}%
\put(-\value{x},-\value{y}){\line(1,-2){#3}}%
\put(\monolength,-\secondepilength){\ssehead}%
\put(\epilength,-\secondmonolength){\line(-1,2){#3}}%
\put(\value{x},\value{y}){\nnwhead}%
\end{picture}}
\put(\value{x},\value{y}){\distsign}%
\put(-\value{x},-\value{y}){\distsign}%
\truex{500}\truey{1000}\truez{600}%
\put(\value{x},\value{x}){\makebox(0,\value{z})[l]{${#1}$}}%
\put(-\value{x},-\value{y}){\makebox(0,\value{z})[r]{${#2}$}}%
\end{picture}}%

\def\basicssear[#1]{\SSEAR{}{}{#100}}%

\newcommand{\ssear}{\@ifnextchar[{\basicssear}{\basicssear[67]}}%

\def\basicSsear[#1]#2{\SSEAR{#2}{}{#100}}%

\newcommand{\Ssear}{\@ifnextchar[{\basicSsear}{\basicSsear[67]}}%

\def\basicsseaR[#1]#2{\SSEAR{}{#2}{#100}}%

\newcommand{\sseaR}{\@ifnextchar[{\basicsseaR}{\basicsseaR[67]}}%

\def\basicssedist[#1]{\SSEDIST{}{}{#100}}%

\newcommand{\ssedist}{\@ifnextchar[{\basicssedist}{\basicssedist[67]}}%

\def\basicSsedist[#1]#2{\SSEDIST{#2}{}{#100}}%

\newcommand{\Ssedist}{\@ifnextchar[{\basicSsedist}{\basicSsedist[67]}}%

\def\basicssedisT[#1]#2{\SSEDIST{}{#2}{#100}}%

\newcommand{\ssedisT}{\@ifnextchar[{\basicssedisT}{\basicssedisT[67]}}%

\def\basicssedotar[#1]{\SSEDOTAR{}{}{#100}}%

\newcommand{\ssedotar}{\@ifnextchar[{\basicssedotar}{\basicssedotar[67]}}%

\def\basicSsedotar[#1]#2{\SSEDOTAR{#2}{}{#100}}%

\newcommand{\Ssedotar}{\@ifnextchar[{\basicSsedotar}{\basicSsedotar[67]}}%

\def\basicssedotaR[#1]#2{\SSEDOTAR{}{#2}{#100}}%

\newcommand{\ssedotaR}{\@ifnextchar[{\basicssedotaR}{\basicssedotaR[67]}}%

\def\basicssemono[#1]{\SSEMONO{}{}{#100}}%

\newcommand{\ssemono}{\@ifnextchar[{\basicssemono}{\basicssemono[67]}}%

\def\basicSsemono[#1]#2{\SSEMONO{#2}{}{#100}}%

\newcommand{\Ssemono}{\@ifnextchar[{\basicSsemono}{\basicSsemono[67]}}%

\def\basicssemonO[#1]#2{\SSEMONO{}{#2}{#100}}%

\newcommand{\ssemonO}{\@ifnextchar[{\basicssemonO}{\basicssemonO[67]}}%

\def\basicsseepi[#1]{\SSEEPI{}{}{#100}}%

\newcommand{\sseepi}{\@ifnextchar[{\basicsseepi}{\basicsseepi[67]}}%

\def\basicSseepi[#1]#2{\SSEEPI{#2}{}{#100}}%

\newcommand{\Sseepi}{\@ifnextchar[{\basicSseepi}{\basicSseepi[67]}}%

\def\basicsseepI[#1]#2{\SSEEPI{}{#2}{#100}}%

\newcommand{\sseepI}{\@ifnextchar[{\basicsseepI}{\basicsseepI[67]}}%

\def\basicssebimo[#1]{\SSEBIMO{}{}{#100}}%

\newcommand{\ssebimo}{\@ifnextchar[{\basicssebimo}{\basicssebimo[67]}}%

\def\basicSsebimo[#1]#2{\SSEBIMO{#2}{}{#100}}%

\newcommand{\Ssebimo}{\@ifnextchar[{\basicSsebimo}{\basicSsebimo[67]}}%

\def\basicssebimO[#1]#2{\SSEBIMO{}{#2}{#100}}%

\newcommand{\ssebimO}{\@ifnextchar[{\basicssebimO}{\basicssebimO[67]}}%

\def\basicsseiso[#1]{\SSEAR{\isosign}{}{#100}}%

\newcommand{\sseiso}{\@ifnextchar[{\basicsseiso}{\basicsseiso[67]}}%

\def\basicSseiso[#1]#2{\SSEAR{#2}{\isosign}{#100}}%

\newcommand{\Sseiso}{\@ifnextchar[{\basicSseiso}{\basicSseiso[67]}}%

\def\basicsseisO[#1]#2{\SSEAR{\isosign}{#2}{#100}}%

\newcommand{\sseisO}{\@ifnextchar[{\basicsseisO}{\basicsseisO[67]}}%

\def\basicsseeql[#1]{\SSEEQL{}{}{#100}}%

\newcommand{\sseeql}{\@ifnextchar[{\basicsseeql}{\basicsseeql[67]}}%

\def\basicSseeql[#1]#2{\SSEEQL{#2}{}{#100}}%

\newcommand{\Sseeql}{\@ifnextchar[{\basicSseeql}{\basicSseeql[67]}}%

\def\basicsseeqL[#1]#2{\SSEEQL{}{#2}{#100}}%

\newcommand{\sseeqL}{\@ifnextchar[{\basicsseeqL}{\basicsseeqL[67]}}%

\def\basicssebiar[#1]{\SSEBIAR{}{}{#100}}%

\newcommand{\ssebiar}{\@ifnextchar[{\basicssebiar}{\basicssebiar[67]}}%

\def\basicSsebiar[#1]#2#3{\SSEBIAR{#2}{#3}{#100}}%

\newcommand{\Ssebiar}{\@ifnextchar[{\basicSsebiar}{\basicSsebiar[67]}}%

\def\basicssebidist[#1]{\SSEBIDIST{}{}{#100}}%

\newcommand{\ssebidist}{\@ifnextchar[{\basicssebidist}{\basicssebidist[67]}}%

\def\basicSsebidist[#1]#2#3{\SSEBIDIST{#2}{#3}{#100}}%

\newcommand{\Ssebidist}{\@ifnextchar[{\basicSsebidist}{\basicSsebidist[67]}}%

\def\basicsseadjar[#1]{\SSEADJAR{}{}{#100}}%

\newcommand{\sseadjar}{\@ifnextchar[{\basicsseadjar}{\basicsseadjar[67]}}%

\def\basicSseadjar[#1]#2#3{\SSEADJAR{#2}{#3}{#100}}%

\newcommand{\Sseadjar}{\@ifnextchar[{\basicSseadjar}{\basicSseadjar[67]}}%

\def\basicsseadjdist[#1]{\SSEADJDIST{}{}{#100}}%

\newcommand{\sseadjdist}{\@ifnextchar[{\basicsseadjdist}{\basicsseadjdist[67]}}%

\def\basicSseadjdist[#1]#2#3{\SSEADJDIST{#2}{#3}{#100}}%

\newcommand{\Sseadjdist}{\@ifnextchar[{\basicSseadjdist}{\basicSseadjdist[67]}}%

\newcommand{\SSWAR}[3]{%
\Z=#3%
\divide\Z by 2%
\begin{picture}(0,0)%
\put(\Z,#3){\line(-1,-2){#3}}%
\put(-\Z,-#3){\sswhead}%
\truex{200}\truey{800}\truez{600}%
\put(-\value{x},\value{x}){\makebox(0,\value{z})[r]{${#1}$}}%
\put(\value{x},-\value{y}){\makebox(0,\value{z})[l]{${#2}$}}%
\end{picture}}%

\newcommand{\SSWDIST}[3]{%
\Z=#3%
\divide\Z by 2%
\begin{picture}(0,0)%
\put(\Z,#3){\line(-1,-2){#3}}%
\put(-\Z,-#3){\sswhead}%
\truex{400}%
\put(0,0){\distsign}%
\truex{200}\truey{800}\truez{600}%
\put(-\value{x},\value{x}){\makebox(0,\value{z})[r]{${#1}$}}%
\put(\value{x},-\value{y}){\makebox(0,\value{z})[l]{${#2}$}}%
\end{picture}}%

\newcommand{\SSWDOTAR}[3]{%
\truex{100}\truey{268}\truez{134}%
\Z=#3%
\divide\Z by 2%
\NUMBEROFDOTS=#3%
\divide\NUMBEROFDOTS by \value{z}%
\advance\NUMBEROFDOTS by 1%
\begin{picture}(0,0)%
\multiput(\Z,#3)(-\value{z},-\value{y}){\NUMBEROFDOTS}%
{\circle*{\value{x}}}%
\put(-\Z,-#3){\sswhead}%
\truex{200}\truey{800}\truez{600}%
\put(-\value{x},\value{x}){\makebox(0,\value{z})[r]{${#1}$}}%
\put(\value{x},-\value{y}){\makebox(0,\value{z})[l]{${#2}$}}%
\end{picture}}%

\newcommand{\SSWMONO}[3]{%
\Z=#3%
\divide\Z by 2%
\truetaiL%
\bimolength=#3%
\advance\bimolength by -\truemonotaiL%
\monolength=\bimolength%
\advance\monolength by -\Z%
\secondmonolength=\monolength%
\multiply\secondmonolength by 2%
\begin{picture}(0,0)%
\put(\monolength,\secondmonolength){\line(-1,-2){\bimolength}}%
\put(\monolength,\secondmonolength){\sswhead}%
\put(-\Z,-#3){\sswhead}%
\truex{200}\truey{800}\truez{600}%
\put(-\value{x},\value{x}){\makebox(0,\value{z})[r]{${#1}$}}%
\put(\value{x},-\value{y}){\makebox(0,\value{z})[l]{${#2}$}}%
\end{picture}}%

\newcommand{\SSWEPI}[3]{%
\Z=#3%
\divide\Z by 2%
\trueheaD%
\bimolength=#3%
\advance\bimolength by -\trueepiheaD%
\epilength=\bimolength%
\advance\epilength by -\Z%
\secondepilength=\epilength%
\multiply\secondepilength by 2%
\begin{picture}(0,0)%
\put(\Z,#3){\line(-1,-2){\bimolength}}%
\put(-\epilength,-\secondepilength){\sswhead}%
\put(-\Z,-#3){\sswhead}%
\truex{200}\truey{800}\truez{600}%
\put(-\value{x},\value{x}){\makebox(0,\value{z})[r]{${#1}$}}%
\put(\value{x},-\value{y}){\makebox(0,\value{z})[l]{${#2}$}}%
\end{picture}}%

\newcommand{\SSWBIMO}[3]{%
\Z=#3%
\divide\Z by 2%
\truetaiL\trueheaD%
\bimolength=#3%
\advance\bimolength by -\truemonotaiL%
\monolength=\bimolength%
\advance\monolength by -\Z%
\advance\bimolength by -\trueepiheaD%
\epilength=\bimolength%
\advance\epilength by -\monolength%
\secondmonolength=\monolength%
\multiply\secondmonolength by 2%
\secondepilength=\epilength%
\multiply\secondepilength by 2%
\begin{picture}(0,0)%
\put(\monolength,\secondmonolength){\line(-1,-2){\bimolength}}%
\put(\monolength,\secondmonolength){\sswhead}%
\put(-\epilength,-\secondepilength){\sswhead}%
\put(-\Z,-#3){\sswhead}%
\truex{200}\truey{800}\truez{600}%
\put(-\value{x},\value{x}){\makebox(0,\value{z})[r]{${#1}$}}%
\put(\value{x},-\value{y}){\makebox(0,\value{z})[l]{${#2}$}}%
\end{picture}}%

\newcommand{\SSWBIAR}[3]{%
\Y=#3%
\divide\Y by 2%
\Z=#3%
\multiply \Z by 2%
\begin{picture}(0,0)%
\put(\Y,#3){\begin{picture}(0,0)%
\truex{313}\truey{156}%
\put(-\value{x},\value{y}){\line(-1,-2){#3}}%
\put(\value{x},-\value{y}){\line(-1,-2){#3}}%
\monolength=#3%
\advance\monolength by -\value{x}%
\epilength=#3%
\advance\epilength by \value{x}%
\secondmonolength=\Z%
\advance\secondmonolength by -\value{y}%
\secondepilength=\Z%
\advance\secondepilength by \value{y}%
\put(-\monolength,-\secondepilength){\sswhead}%
\put(-\epilength,-\secondmonolength){\sswhead}%
\end{picture}}
\truex{300}\truey{1000}\truez{600}%
\put(-\value{x},\value{x}){\makebox(0,\value{z})[r]{${#1}$}}%
\put(\value{x},-\value{y}){\makebox(0,\value{z})[l]{${#2}$}}%
\end{picture}}%

\newcommand{\SSWBIDIST}[3]{%
\Y=#3%
\divide\Y by 2%
\Z=#3%
\multiply \Z by 2%
\begin{picture}(0,0)%
\truex{313}\truey{156}\truez{400}%
\put(\Y,#3){\begin{picture}(0,0)%
\put(-\value{x},\value{y}){\line(-1,-2){#3}}%
\put(\value{x},-\value{y}){\line(-1,-2){#3}}%
\monolength=#3%
\advance\monolength by -\value{x}%
\epilength=#3%
\advance\epilength by \value{x}%
\secondmonolength=\Z%
\advance\secondmonolength by -\value{y}%
\secondepilength=\Z%
\advance\secondepilength by \value{y}%
\put(-\monolength,-\secondepilength){\sswhead}%
\put(-\epilength,-\secondmonolength){\sswhead}%
\end{picture}}
\put(-\value{x},\value{y}){\distsign}%
\put(\value{x},-\value{y}){\distsign}%
\truex{300}\truey{1000}\truez{600}%
\put(-\value{x},\value{x}){\makebox(0,\value{z})[r]{${#1}$}}%
\put(\value{x},-\value{y}){\makebox(0,\value{z})[l]{${#2}$}}%
\end{picture}}%

\newcommand{\SSWADJAR}[3]{%
\Y=#3%
\divide\Y by 2%
\Z=#3%
\multiply \Z by 2%
\begin{picture}(0,0)%
\put(\Y,#3){\begin{picture}(0,0)%
\truex{313}\truey{156}%
\monolength=#3%
\advance\monolength by -\value{x}%
\epilength=#3%
\advance\epilength by \value{x}%
\secondmonolength=\Z%
\advance\secondmonolength by -\value{y}%
\secondepilength=\Z%
\advance\secondepilength by \value{y}%
\put(\value{x},-\value{y}){\line(-1,-2){#3}}%
\put(-\monolength,-\secondepilength){\sswhead}%
\put(-\epilength,-\secondmonolength){\line(1,2){#3}}%
\put(-\value{x},\value{y}){\nnehead}%
\end{picture}}
\truex{300}\truey{1000}\truez{600}%
\put(-\value{x},\value{x}){\makebox(0,\value{z})[r]{${#1}$}}%
\put(\value{x},-\value{y}){\makebox(0,\value{z})[l]{${#2}$}}%
\end{picture}}%

\newcommand{\SSWADJDIST}[3]{%
\Y=#3%
\divide\Y by 2%
\Z=#3%
\multiply \Z by 2%
\begin{picture}(0,0)%
\truex{313}\truey{156}\truez{400}%
\put(\Y,#3){\begin{picture}(0,0)%
\monolength=#3%
\advance\monolength by -\value{x}%
\epilength=#3%
\advance\epilength by \value{x}%
\secondmonolength=\Z%
\advance\secondmonolength by -\value{y}%
\secondepilength=\Z%
\advance\secondepilength by \value{y}%
\put(\value{x},-\value{y}){\line(-1,-2){#3}}%
\put(-\monolength,-\secondepilength){\sswhead}%
\put(-\epilength,-\secondmonolength){\line(1,2){#3}}%
\put(-\value{x},\value{y}){\nnehead}%
\end{picture}}
\put(-\value{x},\value{y}){\distsign}%
\put(\value{x},-\value{y}){\distsign}%
\truex{300}\truey{1000}\truez{600}%
\put(-\value{x},\value{x}){\makebox(0,\value{z})[r]{${#1}$}}%
\put(\value{x},-\value{y}){\makebox(0,\value{z})[l]{${#2}$}}%
\end{picture}}%

\def\basicsswar[#1]{\SSWAR{}{}{#100}}%

\newcommand{\sswar}{\@ifnextchar[{\basicsswar}{\basicsswar[67]}}%

\def\basicSswar[#1]#2{\SSWAR{#2}{}{#100}}%

\newcommand{\Sswar}{\@ifnextchar[{\basicSswar}{\basicSswar[67]}}%

\def\basicsswaR[#1]#2{\SSWAR{}{#2}{#100}}%

\newcommand{\sswaR}{\@ifnextchar[{\basicsswaR}{\basicsswaR[67]}}%

\def\basicsswdist[#1]{\SSWDIST{}{}{#100}}%

\newcommand{\sswdist}{\@ifnextchar[{\basicsswdist}{\basicsswdist[67]}}%

\def\basicSswdist[#1]#2{\SSWDIST{#2}{}{#100}}%

\newcommand{\Sswdist}{\@ifnextchar[{\basicSswdist}{\basicSswdist[67]}}%

\def\basicsswdisT[#1]#2{\SSWDIST{}{#2}{#100}}%

\newcommand{\sswdisT}{\@ifnextchar[{\basicsswdisT}{\basicsswdisT[67]}}%

\def\basicsswdotar[#1]{\SSWDOTAR{}{}{#100}}%

\newcommand{\sswdotar}{\@ifnextchar[{\basicsswdotar}{\basicsswdotar[67]}}%

\def\basicSswdotar[#1]#2{\SSWDOTAR{#2}{}{#100}}%

\newcommand{\Sswdotar}{\@ifnextchar[{\basicSswdotar}{\basicSswdotar[67]}}%

\def\basicsswdotaR[#1]#2{\SSWDOTAR{}{#2}{#100}}%

\newcommand{\sswdotaR}{\@ifnextchar[{\basicsswdotaR}{\basicsswdotaR[67]}}%

\def\basicsswmono[#1]{\SSWMONO{}{}{#100}}%

\newcommand{\sswmono}{\@ifnextchar[{\basicsswmono}{\basicsswmono[67]}}%

\def\basicSswmono[#1]#2{\SSWMONO{#2}{}{#100}}%

\newcommand{\Sswmono}{\@ifnextchar[{\basicSswmono}{\basicSswmono[67]}}%

\def\basicsswmonO[#1]#2{\SSWMONO{}{#2}{#100}}%

\newcommand{\sswmonO}{\@ifnextchar[{\basicsswmonO}{\basicsswmonO[67]}}%

\def\basicsswepi[#1]{\SSWEPI{}{}{#100}}%

\newcommand{\sswepi}{\@ifnextchar[{\basicsswepi}{\basicsswepi[67]}}%

\def\basicSswepi[#1]#2{\SSWEPI{#2}{}{#100}}%

\newcommand{\Sswepi}{\@ifnextchar[{\basicSswepi}{\basicSswepi[67]}}%

\def\basicsswepI[#1]#2{\SSWEPI{}{#2}{#100}}%

\newcommand{\sswepI}{\@ifnextchar[{\basicsswepI}{\basicsswepI[67]}}%

\def\basicsswbimo[#1]{\SSWBIMO{}{}{#100}}%

\newcommand{\sswbimo}{\@ifnextchar[{\basicsswbimo}{\basicsswbimo[67]}}%

\def\basicSswbimo[#1]#2{\SSWBIMO{#2}{}{#100}}%

\newcommand{\Sswbimo}{\@ifnextchar[{\basicSswbimo}{\basicSswbimo[67]}}%

\def\basicsswbimO[#1]#2{\SSWBIMO{}{#2}{#100}}%

\newcommand{\sswbimO}{\@ifnextchar[{\basicsswbimO}{\basicsswbimO[67]}}%

\def\basicsswiso[#1]{\SSWAR{\isosign}{}{#100}}%

\newcommand{\sswiso}{\@ifnextchar[{\basicsswiso}{\basicsswiso[67]}}%

\def\basicSswiso[#1]#2{\SSWAR{#2}{\isosign}{#100}}%

\newcommand{\Sswiso}{\@ifnextchar[{\basicSswiso}{\basicSswiso[67]}}%

\def\basicsswisO[#1]#2{\SSWAR{\isosign}{#2}{#100}}%

\newcommand{\sswisO}{\@ifnextchar[{\basicsswisO}{\basicsswisO[67]}}%

\def\basicsswbiar[#1]{\SSWBIAR{}{}{#100}}%

\newcommand{\sswbiar}{\@ifnextchar[{\basicsswbiar}{\basicsswbiar[67]}}%

\def\basicSswbiar[#1]#2#3{\SSWBIAR{#2}{#3}{#100}}%

\newcommand{\Sswbiar}{\@ifnextchar[{\basicSswbiar}{\basicSswbiar[67]}}%

\def\basicsswbidist[#1]{\SSWBIDIST{}{}{#100}}%

\newcommand{\sswbidist}{\@ifnextchar[{\basicsswbidist}{\basicsswbidist[67]}}%

\def\basicSswbidist[#1]#2#3{\SSWBIDIST{#2}{#3}{#100}}%

\newcommand{\Sswbidist}{\@ifnextchar[{\basicSswbidist}{\basicSswbidist[67]}}%

\def\basicsswadjar[#1]{\SSWADJAR{}{}{#100}}%

\newcommand{\sswadjar}{\@ifnextchar[{\basicsswadjar}{\basicsswadjar[67]}}%

\def\basicSswadjar[#1]#2#3{\SSWADJAR{#2}{#3}{#100}}%

\newcommand{\Sswadjar}{\@ifnextchar[{\basicSswadjar}{\basicSswadjar[67]}}%

\def\basicsswadjdist[#1]{\SSWADJDIST{}{}{#100}}%

\newcommand{\sswadjdist}{\@ifnextchar[{\basicsswadjdist}{\basicsswadjdist[67]}}%

\def\basicSswadjdist[#1]#2#3{\SSWADJDIST{#2}{#3}{#100}}%

\newcommand{\Sswadjdist}{\@ifnextchar[{\basicSswadjdist}{\basicSswadjdist[67]}}%

\newcommand{\NNWAR}[3]{%
\Z=#3%
\divide\Z by 2%
\begin{picture}(0,0)%
\put(\Z,-#3){\line(-1,2){#3}}%
\put(-\Z,#3){\nnwhead}%
\truex{200}\truey{800}\truez{600}%
\put(\value{x},\value{x}){\makebox(0,\value{z})[l]{${#1}$}}%
\put(-\value{x},-\value{y}){\makebox(0,\value{z})[r]{${#2}$}}%
\end{picture}}%

\newcommand{\NNWDIST}[3]{%
\Z=#3%
\divide\Z by 2%
\begin{picture}(0,0)%
\put(\Z,-#3){\line(-1,2){#3}}%
\put(-\Z,#3){\nnwhead}%
\truex{400}%
\put(0,0){\distsign}%
\truex{200}\truey{800}\truez{600}%
\put(\value{x},\value{x}){\makebox(0,\value{z})[l]{${#1}$}}%
\put(-\value{x},-\value{y}){\makebox(0,\value{z})[r]{${#2}$}}%
\end{picture}}%

\newcommand{\NNWDOTAR}[3]{%
\truex{100}\truey{268}\truez{134}%
\Z=#3%
\divide\Z by 2%
\NUMBEROFDOTS=#3%
\divide\NUMBEROFDOTS by \value{z}%
\advance\NUMBEROFDOTS by 1%
\begin{picture}(0,0)%
\multiput(\Z,-#3)(-\value{z},\value{y}){\NUMBEROFDOTS}%
{\circle*{\value{x}}}%
\put(-\Z,#3){\nnwhead}%
\truex{200}\truey{800}\truez{600}%
\put(\value{x},\value{x}){\makebox(0,\value{z})[l]{${#1}$}}%
\put(-\value{x},-\value{y}){\makebox(0,\value{z})[r]{${#2}$}}%
\end{picture}}%

\newcommand{\NNWMONO}[3]{%
\Z=#3%
\divide\Z by 2%
\truetaiL%
\bimolength=#3%
\advance\bimolength by -\truemonotaiL%
\monolength=\bimolength%
\advance\monolength by -\Z%
\secondmonolength=\monolength%
\multiply\secondmonolength by 2%
\begin{picture}(0,0)%
\put(\monolength,-\secondmonolength){\line(-1,2){\bimolength}}%
\put(\monolength,-\secondmonolength){\nnwhead}%
\put(-\Z,#3){\nnwhead}%
\truex{200}\truey{800}\truez{600}%
\put(\value{x},\value{x}){\makebox(0,\value{z})[l]{${#1}$}}%
\put(-\value{x},-\value{y}){\makebox(0,\value{z})[r]{${#2}$}}%
\end{picture}}%

\newcommand{\NNWEPI}[3]{%
\Z=#3%
\divide\Z by 2%
\trueheaD%
\bimolength=#3%
\advance\bimolength by -\trueepiheaD%
\epilength=\bimolength%
\advance\epilength by -\Z%
\secondepilength=\epilength%
\multiply\secondepilength by 2%
\begin{picture}(0,0)%
\put(\Z,-#3){\line(-1,2){\bimolength}}%
\put(-\epilength,\secondepilength){\nnwhead}%
\put(-\Z,#3){\nnwhead}%
\truex{200}\truey{800}\truez{600}%
\put(\value{x},\value{x}){\makebox(0,\value{z})[l]{${#1}$}}%
\put(-\value{x},-\value{y}){\makebox(0,\value{z})[r]{${#2}$}}%
\end{picture}}%

\newcommand{\NNWBIMO}[3]{%
\Z=#3%
\divide\Z by 2%
\truetaiL\trueheaD%
\bimolength=#3%
\advance\bimolength by -\truemonotaiL%
\monolength=\bimolength%
\advance\monolength by -\Z%
\advance\bimolength by -\trueepiheaD%
\epilength=\bimolength%
\advance\epilength by -\monolength%
\secondmonolength=\monolength%
\multiply\secondmonolength by 2%
\secondepilength=\epilength%
\multiply\secondepilength by 2%
\begin{picture}(0,0)%
\put(\monolength,-\secondmonolength){\line(-1,2){\bimolength}}%
\put(\monolength,-\secondmonolength){\nnwhead}%
\put(-\epilength,\secondepilength){\nnwhead}%
\put(-\Z,#3){\nnwhead}%
\truex{200}\truey{800}\truez{600}%
\put(\value{x},\value{x}){\makebox(0,\value{z})[l]{${#1}$}}%
\put(-\value{x},-\value{y}){\makebox(0,\value{z})[r]{${#2}$}}%
\end{picture}}%

\newcommand{\NNWBIAR}[3]{%
\Y=#3%
\divide\Y by 2%
\Z=#3%
\multiply \Z by 2%
\begin{picture}(0,0)%
\put(\Y,-#3){\begin{picture}(0,0)%
\truex{313}\truey{156}%
\put(-\value{x},-\value{y}){\line(-1,2){#3}}%
\put(\value{x},\value{y}){\line(-1,2){#3}}%
\monolength=#3%
\advance\monolength by -\value{x}%
\epilength=#3%
\advance\epilength by \value{x}%
\secondmonolength=\Z%
\advance\secondmonolength by -\value{y}%
\secondepilength=\Z%
\advance\secondepilength by \value{y}%
\put(-\monolength,\secondepilength){\nnwhead}%
\put(-\epilength,\secondmonolength){\nnwhead}%
\end{picture}}
\truex{400}\truey{1000}\truez{600}%
\put(\value{x},\value{x}){\makebox(0,\value{z})[l]{${#1}$}}%
\put(-\value{x},-\value{y}){\makebox(0,\value{z})[r]{${#2}$}}%
\end{picture}}%

\newcommand{\NNWBIDIST}[3]{%
\Y=#3%
\divide\Y by 2%
\Z=#3%
\multiply \Z by 2%
\begin{picture}(0,0)%
\truex{313}\truey{156}\truez{400}%
\put(\Y,-#3){\begin{picture}(0,0)%
\put(-\value{x},-\value{y}){\line(-1,2){#3}}%
\put(\value{x},\value{y}){\line(-1,2){#3}}%
\monolength=#3%
\advance\monolength by -\value{x}%
\epilength=#3%
\advance\epilength by \value{x}%
\secondmonolength=\Z%
\advance\secondmonolength by -\value{y}%
\secondepilength=\Z%
\advance\secondepilength by \value{y}%
\put(-\monolength,\secondepilength){\nnwhead}%
\put(-\epilength,\secondmonolength){\nnwhead}%
\end{picture}}
\put(-\value{x},-\value{y}){\distsign}%
\put(\value{x},\value{y}){\distsign}%
\truex{500}\truey{1000}\truez{600}%
\put(\value{x},\value{x}){\makebox(0,\value{z})[l]{${#1}$}}%
\put(-\value{x},-\value{y}){\makebox(0,\value{z})[r]{${#2}$}}%
\end{picture}}%

\newcommand{\NNWADJAR}[3]{%
\Y=#3%
\divide\Y by 2%
\Z=#3%
\multiply \Z by 2%
\begin{picture}(0,0)%
\put(\Y,-#3){\begin{picture}(0,0)%
\truex{313}\truey{156}%
\monolength=#3%
\advance\monolength by -\value{x}%
\epilength=#3%
\advance\epilength by \value{x}%
\secondmonolength=\Z%
\advance\secondmonolength by -\value{y}%
\secondepilength=\Z%
\advance\secondepilength by \value{y}%
\put(-\value{x},-\value{y}){\line(-1,2){#3}}%
\put(-\epilength,\secondmonolength){\nnwhead}%
\put(-\monolength,\secondepilength){\line(1,-2){#3}}%
\put(\value{x},\value{y}){\ssehead}%
\end{picture}}
\truex{400}\truey{1000}\truez{600}%
\put(\value{x},\value{x}){\makebox(0,\value{z})[l]{${#1}$}}%
\put(-\value{x},-\value{y}){\makebox(0,\value{z})[r]{${#2}$}}%
\end{picture}}%

\newcommand{\NNWADJDIST}[3]{%
\Y=#3%
\divide\Y by 2%
\Z=#3%
\multiply \Z by 2%
\begin{picture}(0,0)%
\truex{313}\truey{156}\truez{400}%
\put(\Y,-#3){\begin{picture}(0,0)%
\monolength=#3%
\advance\monolength by -\value{x}%
\epilength=#3%
\advance\epilength by \value{x}%
\secondmonolength=\Z%
\advance\secondmonolength by -\value{y}%
\secondepilength=\Z%
\advance\secondepilength by \value{y}%
\put(-\value{x},-\value{y}){\line(-1,2){#3}}%
\put(-\epilength,\secondmonolength){\nnwhead}%
\put(-\monolength,\secondepilength){\line(1,-2){#3}}%
\put(\value{x},\value{y}){\ssehead}%
\end{picture}}
\put(\value{x},\value{y}){\distsign}%
\put(-\value{x},-\value{y}){\distsign}%
\truex{500}\truey{1000}\truez{600}%
\put(\value{x},\value{x}){\makebox(0,\value{z})[l]{${#1}$}}%
\put(-\value{x},-\value{y}){\makebox(0,\value{z})[r]{${#2}$}}%
\end{picture}}%

\def\basicnnwar[#1]{\NNWAR{}{}{#100}}%

\newcommand{\nnwar}{\@ifnextchar[{\basicnnwar}{\basicnnwar[67]}}%

\def\basicNnwar[#1]#2{\NNWAR{#2}{}{#100}}%

\newcommand{\Nnwar}{\@ifnextchar[{\basicNnwar}{\basicNnwar[67]}}%

\def\basicnnwaR[#1]#2{\NNWAR{}{#2}{#100}}%

\newcommand{\nnwaR}{\@ifnextchar[{\basicnnwaR}{\basicnnwaR[67]}}%

\def\basicnnwdist[#1]{\NNWDIST{}{}{#100}}%

\newcommand{\nnwdist}{\@ifnextchar[{\basicnnwdist}{\basicnnwdist[67]}}%

\def\basicNnwdist[#1]#2{\NNWDIST{#2}{}{#100}}%

\newcommand{\Nnwdist}{\@ifnextchar[{\basicNnwdist}{\basicNnwdist[67]}}%

\def\basicnnwdisT[#1]#2{\NNWDIST{}{#2}{#100}}%

\newcommand{\nnwdisT}{\@ifnextchar[{\basicnnwdisT}{\basicnnwdisT[67]}}%

\def\basicnnwdotar[#1]{\NNWDOTAR{}{}{#100}}%

\newcommand{\nnwdotar}{\@ifnextchar[{\basicnnwdotar}{\basicnnwdotar[67]}}%

\def\basicNnwdotar[#1]#2{\NNWDOTAR{#2}{}{#100}}%

\newcommand{\Nnwdotar}{\@ifnextchar[{\basicNnwdotar}{\basicNnwdotar[67]}}%

\def\basicnnwdotaR[#1]#2{\NNWDOTAR{}{#2}{#100}}%

\newcommand{\nnwdotaR}{\@ifnextchar[{\basicnnwdotaR}{\basicnnwdotaR[67]}}%

\def\basicnnwmono[#1]{\NNWMONO{}{}{#100}}%

\newcommand{\nnwmono}{\@ifnextchar[{\basicnnwmono}{\basicnnwmono[67]}}%

\def\basicNnwmono[#1]#2{\NNWMONO{#2}{}{#100}}%

\newcommand{\Nnwmono}{\@ifnextchar[{\basicNnwmono}{\basicNnwmono[67]}}%

\def\basicnnwmonO[#1]#2{\NNWMONO{}{#2}{#100}}%

\newcommand{\nnwmonO}{\@ifnextchar[{\basicnnwmonO}{\basicnnwmonO[67]}}%

\def\basicnnwepi[#1]{\NNWEPI{}{}{#100}}%

\newcommand{\nnwepi}{\@ifnextchar[{\basicnnwepi}{\basicnnwepi[67]}}%

\def\basicNnwepi[#1]#2{\NNWEPI{#2}{}{#100}}%

\newcommand{\Nnwepi}{\@ifnextchar[{\basicNnwepi}{\basicNnwepi[67]}}%

\def\basicnnwepI[#1]#2{\NNWEPI{}{#2}{#100}}%

\newcommand{\nnwepI}{\@ifnextchar[{\basicnnwepI}{\basicnnwepI[67]}}%

\def\basicnnwbimo[#1]{\NNWBIMO{}{}{#100}}%

\newcommand{\nnwbimo}{\@ifnextchar[{\basicnnwbimo}{\basicnnwbimo[67]}}%

\def\basicNnwbimo[#1]#2{\NNWBIMO{#2}{}{#100}}%

\newcommand{\Nnwbimo}{\@ifnextchar[{\basicNnwbimo}{\basicNnwbimo[67]}}%

\def\basicnnwbimO[#1]#2{\NNWBIMO{}{#2}{#100}}%

\newcommand{\nnwbimO}{\@ifnextchar[{\basicnnwbimO}{\basicnnwbimO[67]}}%

\def\basicnnwiso[#1]{\NNWAR{\isosign}{}{#100}}%

\newcommand{\nnwiso}{\@ifnextchar[{\basicnnwiso}{\basicnnwiso[67]}}%

\def\basicNnwiso[#1]#2{\NNWAR{#2}{\isosign}{#100}}%

\newcommand{\Nnwiso}{\@ifnextchar[{\basicNnwiso}{\basicNnwiso[67]}}%

\def\basicnnwisO[#1]#2{\NNWAR{\isosign}{#2}{#100}}%

\newcommand{\nnwisO}{\@ifnextchar[{\basicnnwisO}{\basicnnwisO[67]}}%

\def\basicnnwbiar[#1]{\NNWBIAR{}{}{#100}}%

\newcommand{\nnwbiar}{\@ifnextchar[{\basicnnwbiar}{\basicnnwbiar[67]}}%

\def\basicNnwbiar[#1]#2#3{\NNWBIAR{#2}{#3}{#100}}%

\newcommand{\Nnwbiar}{\@ifnextchar[{\basicNnwbiar}{\basicNnwbiar[67]}}%

\def\basicnnwbidist[#1]{\NNWBIDIST{}{}{#100}}%

\newcommand{\nnwbidist}{\@ifnextchar[{\basicnnwbidist}{\basicnnwbidist[67]}}%

\def\basicNnwbidist[#1]#2#3{\NNWBIDIST{#2}{#3}{#100}}%

\newcommand{\Nnwbidist}{\@ifnextchar[{\basicNnwbidist}{\basicNnwbidist[67]}}%

\def\basicnnwadjar[#1]{\NNWADJAR{}{}{#100}}%

\newcommand{\nnwadjar}{\@ifnextchar[{\basicnnwadjar}{\basicnnwadjar[67]}}%

\def\basicNnwadjar[#1]#2#3{\NNWADJAR{#2}{#3}{#100}}%

\newcommand{\Nnwadjar}{\@ifnextchar[{\basicNnwadjar}{\basicNnwadjar[67]}}%

\def\basicnnwadjdist[#1]{\NNWADJDIST{}{}{#100}}%

\newcommand{\nnwadjdist}{\@ifnextchar[{\basicnnwadjdist}{\basicnnwadjdist[67]}}%

\def\basicNnwadjdist[#1]#2#3{\NNWADJDIST{#2}{#3}{#100}}%

\newcommand{\Nnwadjdist}{\@ifnextchar[{\basicNnwadjdist}{\basicNnwadjdist[67]}}%

\newcommand{\EENEAR}[3]{%
\Y=#3%
\divide \Y by 2%
\Z=\Y%
\divide \Z by 3%
\begin{picture}(0,0)%
\put(-\Y,-\Z){\line(3,1){#3}}%
\put(\Y,\Z){\eenehead}%
\truex{200}\truey{800}\truez{600}%
\put(-\value{x},\value{x}){\makebox(0,\value{z})[r]{${#1}$}}%
\put(\value{x},-\value{y}){\makebox(0,\value{z})[l]{${#2}$}}%
\end{picture}}%

\def\basiceenear[#1]{\EENEAR{}{}{#100}}%

\newcommand{\eenear}{\@ifnextchar[{\basiceenear}{\basiceenear[211]}}%

\def\basicEenear[#1]#2{\EENEAR{#2}{}{#100}}%

\newcommand{\Eenear}{\@ifnextchar[{\basicEenear}{\basicEenear[211]}}%

\def\basiceeneaR[#1]#2{\EENEAR{}{#2}{#100}}%

\newcommand{\eeneaR}{\@ifnextchar[{\basiceeneaR}{\basiceeneaR[211]}}%

\newcommand{\EESEAR}[3]{%
\Y=#3%
\divide \Y by 2%
\Z=\Y%
\divide \Z by 3%
\begin{picture}(0,0)%
\put(-\Y,\Z){\line(3,-1){#3}}%
\put(\Y,-\Z){\eesehead}%
\truex{200}\truey{800}\truez{600}%
\put(\value{x},\value{x}){\makebox(0,\value{z})[l]{${#1}$}}%
\put(-\value{x},-\value{y}){\makebox(0,\value{z})[r]{${#2}$}}%
\end{picture}}%

\def\basiceesear[#1]{\EESEAR{}{}{#100}}%

\newcommand{\eesear}{\@ifnextchar[{\basiceesear}{\basiceesear[211]}}%

\def\basicEesear[#1]#2{\EESEAR{#2}{}{#100}}%

\newcommand{\Eesear}{\@ifnextchar[{\basicEesear}{\basicEesear[211]}}%

\def\basiceeseaR[#1]#2{\EESEAR{}{#2}{#100}}%

\newcommand{\eeseaR}{\@ifnextchar[{\basiceeseaR}{\basiceeseaR[211]}}%

\newcommand{\WWNWAR}[3]{%
\Y=#3%
\divide \Y by 2%
\Z=\Y%
\divide \Z by 3%
\begin{picture}(0,0)%
\put(\Y,-\Z){\line(-3,1){#3}}%
\put(-\Y,\Z){\wwnwhead}%
\truex{200}\truey{800}\truez{600}%
\put(\value{x},\value{x}){\makebox(0,\value{z})[l]{${#1}$}}%
\put(-\value{x},-\value{y}){\makebox(0,\value{z})[r]{${#2}$}}%
\end{picture}}%

\def\basicwwnwar[#1]{\WWNWAR{}{}{#100}}%

\newcommand{\wwnwar}{\@ifnextchar[{\basicwwnwar}{\basicwwnwar[211]}}%

\def\basicWwnwar[#1]#2{\WWNWAR{#2}{}{#100}}%

\newcommand{\Wwnwar}{\@ifnextchar[{\basicWwnwar}{\basicWwnwar[211]}}%

\def\basicwwnwaR[#1]#2{\WWNWAR{}{#2}{#100}}%

\newcommand{\wwnwaR}{\@ifnextchar[{\basicwwnwaR}{\basicwwnwaR[211]}}%

\newcommand{\WWSWAR}[3]{%
\Y=#3%
\divide \Y by 2%
\Z=\Y%
\divide \Z by 3%
\begin{picture}(0,0)%
\put(\Y,\Z){\line(-3,-1){#3}}%
\put(-\Y,-\Z){\wwswhead}%
\truex{200}\truey{800}\truez{600}%
\put(-\value{x},\value{x}){\makebox(0,\value{z})[r]{${#1}$}}%
\put(\value{x},-\value{y}){\makebox(0,\value{z})[l]{${#2}$}}%
\end{picture}}%

\def\basicwwswar[#1]{\WWSWAR{}{}{#100}}%

\newcommand{\wwswar}{\@ifnextchar[{\basicwwswar}{\basicwwswar[211]}}%

\def\basicWwswar[#1]#2{\WWSWAR{#2}{}{#100}}%

\newcommand{\Wwswar}{\@ifnextchar[{\basicWwswar}{\basicWwswar[211]}}%

\def\basicwwswaR[#1]#2{\WWSWAR{}{#2}{#100}}%

\newcommand{\wwswaR}{\@ifnextchar[{\basicwwswaR}{\basicwwswaR[211]}}%

\newcommand{\NNNEAR}[3]{%
\Y=#3%
\divide \Y by 2%
\Z=\Y%
\multiply \Z by 3%
\begin{picture}(0,0)%
\put(-\Y,-\Z){\line(1,3){#3}}%
\put(\Y,\Z){\nnnehead}%
\truex{100}\truez{600}%
\put(-\value{x},\value{x}){\makebox(0,\value{z})[r]{${#1}$}}%
\put(\value{x},-\value{z}){\makebox(0,\value{z})[l]{${#2}$}}%
\end{picture}}%

\def\basicnnnear[#1]{\NNNEAR{}{}{#100}}%

\newcommand{\nnnear}{\@ifnextchar[{\basicnnnear}{\basicnnnear[71]}}%

\def\basicNnnear[#1]#2{\NNNEAR{#2}{}{#100}}%

\newcommand{\Nnnear}{\@ifnextchar[{\basicNnnear}{\basicNnnear[71]}}%

\def\basicnnneaR[#1]#2{\NNNEAR{}{#2}{#100}}%

\newcommand{\nnneaR}{\@ifnextchar[{\basicnnneaR}{\basicnnneaR[71]}}%

\newcommand{\SSSWAR}[3]{%
\Y=#3%
\divide \Y by 2%
\Z=\Y%
\multiply \Z by 3%
\begin{picture}(0,0)%
\put(\Y,\Z){\line(-1,-3){#3}}%
\put(-\Y,-\Z){\ssswhead}%
\truex{100}\truez{600}%
\put(-\value{x},\value{x}){\makebox(0,\value{z})[r]{${#1}$}}%
\put(\value{x},-\value{z}){\makebox(0,\value{z})[l]{${#2}$}}%
\end{picture}}%

\def\basicssswar[#1]{\SSSWAR{}{}{#100}}%

\newcommand{\ssswar}{\@ifnextchar[{\basicssswar}{\basicssswar[71]}}%

\def\basicSsswar[#1]#2{\SSSWAR{#2}{}{#100}}%

\newcommand{\Ssswar}{\@ifnextchar[{\basicSsswar}{\basicSsswar[71]}}%

\def\basicssswaR[#1]#2{\SSSWAR{}{#2}{#100}}%

\newcommand{\ssswaR}{\@ifnextchar[{\basicssswaR}{\basicssswaR[71]}}%

\newcommand{\SSSEAR}[3]{%
\Y=#3%
\divide \Y by 2%
\Z=\Y%
\multiply \Z by 3%
\begin{picture}(0,0)%
\put(-\Y,\Z){\line(1,-3){#3}}%
\put(\Y,-\Z){\sssehead}%
\truex{200}\truez{600}%
\put(\value{x},\value{x}){\makebox(0,\value{z})[l]{${#1}$}}%
\put(-\value{x},-\value{z}){\makebox(0,\value{z})[r]{${#2}$}}%
\end{picture}}%

\def\basicsssear[#1]{\SSSEAR{}{}{#100}}%

\newcommand{\sssear}{\@ifnextchar[{\basicsssear}{\basicsssear[71]}}%

\def\basicSssear[#1]#2{\SSSEAR{#2}{}{#100}}%

\newcommand{\Sssear}{\@ifnextchar[{\basicSssear}{\basicSssear[71]}}%

\def\basicssseaR[#1]#2{\SSSEAR{}{#2}{#100}}%

\newcommand{\ssseaR}{\@ifnextchar[{\basicssseaR}{\basicssseaR[71]}}%

\newcommand{\NNNWAR}[3]{%
\Y=#3%
\divide \Y by 2%
\Z=\Y%
\multiply \Z by 3%
\begin{picture}(0,0)%
\put(\Y,-\Z){\line(-1,3){#3}}%
\put(-\Y,\Z){\nnnwhead}%
\truex{200}\truez{600}%
\put(\value{x},\value{x}){\makebox(0,\value{z})[l]{${#1}$}}%
\put(-\value{x},-\value{z}){\makebox(0,\value{z})[r]{${#2}$}}%
\end{picture}}%

\def\basicnnnwar[#1]{\NNNWAR{}{}{#100}}%

\newcommand{\nnnwar}{\@ifnextchar[{\basicnnnwar}{\basicnnnwar[71]}}%

\def\basicNnnwar[#1]#2{\NNNWAR{#2}{}{#100}}%

\newcommand{\Nnnwar}{\@ifnextchar[{\basicNnnwar}{\basicNnnwar[71]}}%

\def\basicnnnwaR[#1]#2{\NNNWAR{}{#2}{#100}}%

\newcommand{\nnnwaR}{\@ifnextchar[{\basicnnnwaR}{\basicnnnwaR[71]}}%

\newcommand{\NEENEAR}[3]{%
\Y=#3%
\divide \Y by 2%
\Z=#3%
\divide \Z by 3%
\begin{picture}(0,0)%
\put(-\Y,-\Z){\line(3,2){#3}}%
\put(\Y,\Z){\neenehead}%
\truex{200}\truey{800}\truez{600}%
\put(-\value{x},\value{x}){\makebox(0,\value{z})[r]{${#1}$}}%
\put(\value{x},-\value{y}){\makebox(0,\value{z})[l]{${#2}$}}%
\end{picture}}%

\def\basicneenear[#1]{\NEENEAR{}{}{#100}}%

\newcommand{\neenear}{\@ifnextchar[{\basicneenear}{\basicneenear[215]}}%

\def\basicNeenear[#1]#2{\NEENEAR{#2}{}{#100}}%

\newcommand{\Neenear}{\@ifnextchar[{\basicNeenear}{\basicNeenear[215]}}%

\def\basicneeneaR[#1]#2{\NEENEAR{}{#2}{#100}}%

\newcommand{\neeneaR}{\@ifnextchar[{\basicneeneaR}{\basicneeneaR[215]}}%

\newcommand{\SEESEAR}[3]{%
\Y=#3%
\divide \Y by 2%
\Z=#3%
\divide \Z by 3%
\begin{picture}(0,0)%
\put(-\Y,\Z){\line(3,-2){#3}}%
\put(\Y,-\Z){\seesehead}%
\truex{200}\truey{800}\truez{600}%
\put(\value{x},\value{x}){\makebox(0,\value{z})[l]{${#1}$}}%
\put(-\value{x},-\value{y}){\makebox(0,\value{z})[r]{${#2}$}}%
\end{picture}}%

\def\basicseesear[#1]{\SEESEAR{}{}{#100}}%

\newcommand{\seesear}{\@ifnextchar[{\basicseesear}{\basicseesear[215]}}%

\def\basicSeesear[#1]#2{\SEESEAR{#2}{}{#100}}%

\newcommand{\Seesear}{\@ifnextchar[{\basicSeesear}{\basicSeesear[215]}}%

\def\basicseeseaR[#1]#2{\SEESEAR{}{#2}{#100}}%

\newcommand{\seeseaR}{\@ifnextchar[{\basicseeseaR}{\basicseeseaR[215]}}%

\newcommand{\NWWNWAR}[3]{%
\Y=#3%
\divide \Y by 2%
\Z=#3%
\divide \Z by 3%
\begin{picture}(0,0)%
\put(\Y,-\Z){\line(-3,2){#3}}%
\put(-\Y,\Z){\nwwnwhead}%
\truex{200}\truey{800}\truez{600}%
\put(\value{x},\value{x}){\makebox(0,\value{z})[l]{${#1}$}}%
\put(-\value{x},-\value{y}){\makebox(0,\value{z})[r]{${#2}$}}%
\end{picture}}%

\def\basicnwwnwar[#1]{\NWWNWAR{}{}{#100}}%

\newcommand{\nwwnwar}{\@ifnextchar[{\basicnwwnwar}{\basicnwwnwar[215]}}%

\def\basicNwwnwar[#1]#2{\NWWNWAR{#2}{}{#100}}%

\newcommand{\Nwwnwar}{\@ifnextchar[{\basicNwwnwar}{\basicNwwnwar[215]}}%

\def\basicnwwnwaR[#1]#2{\NWWNWAR{}{#2}{#100}}%

\newcommand{\nwwnwaR}{\@ifnextchar[{\basicnwwnwaR}{\basicnwwnwaR[215]}}%

\newcommand{\SWWSWAR}[3]{%
\Y=#3%
\divide \Y by 2%
\Z=#3%
\divide \Z by 3%
\begin{picture}(0,0)%
\put(\Y,\Z){\line(-3,-2){#3}}%
\put(-\Y,-\Z){\swwswhead}%
\truex{200}\truey{800}\truez{600}%
\put(-\value{x},\value{x}){\makebox(0,\value{z})[r]{${#1}$}}%
\put(\value{x},-\value{y}){\makebox(0,\value{z})[l]{${#2}$}}%
\end{picture}}%

\def\basicswwswar[#1]{\SWWSWAR{}{}{#100}}%

\newcommand{\swwswar}{\@ifnextchar[{\basicswwswar}{\basicswwswar[215]}}%

\def\basicSwwswar[#1]#2{\SWWSWAR{#2}{}{#100}}%

\newcommand{\Swwswar}{\@ifnextchar[{\basicSwwswar}{\basicSwwswar[215]}}%

\def\basicswwswaR[#1]#2{\SWWSWAR{}{#2}{#100}}%

\newcommand{\swwswaR}{\@ifnextchar[{\basicswwswaR}{\basicswwswaR[215]}}%

\newcommand{\NENNEAR}[3]{%
\Y=#3%
\divide \Y by 2%
\Z=#3%
\multiply \Z by 3%
\divide \Z by 4%
\begin{picture}(0,0)%
\put(-\Y,-\Z){\line(2,3){#3}}%
\put(\Y,\Z){\nennehead}%
\truex{100}\truez{600}%
\put(-\value{x},\value{x}){\makebox(0,\value{z})[r]{${#1}$}}%
\put(\value{x},-\value{z}){\makebox(0,\value{z})[l]{${#2}$}}%
\end{picture}}%

\def\basicnennear[#1]{\NENNEAR{}{}{#100}}%

\newcommand{\nennear}{\@ifnextchar[{\basicnennear}{\basicnennear[143]}}%

\def\basicNennear[#1]#2{\NENNEAR{#2}{}{#100}}%

\newcommand{\Nennear}{\@ifnextchar[{\basicNennear}{\basicNennear[143]}}%

\def\basicnenneaR[#1]#2{\NENNEAR{}{#2}{#100}}%

\newcommand{\nenneaR}{\@ifnextchar[{\basicnenneaR}{\basicnenneaR[143]}}%

\newcommand{\SWSSWAR}[3]{%
\Y=#3%
\divide \Y by 2%
\Z=#3%
\multiply \Z by 3%
\divide \Z by 4%
\begin{picture}(0,0)%
\put(\Y,\Z){\line(-2,-3){#3}}%
\put(-\Y,-\Z){\swsswhead}%
\truex{100}\truez{600}%
\put(-\value{x},\value{x}){\makebox(0,\value{z})[r]{${#1}$}}%
\put(\value{x},-\value{z}){\makebox(0,\value{z})[l]{${#2}$}}%
\end{picture}}%

\def\basicswsswar[#1]{\SWSSWAR{}{}{#100}}%

\newcommand{\swsswar}{\@ifnextchar[{\basicswsswar}{\basicswsswar[143]}}%

\def\basicSwsswar[#1]#2{\SWSSWAR{#2}{}{#100}}%

\newcommand{\Swsswar}{\@ifnextchar[{\basicSwsswar}{\basicSwsswar[143]}}%

\def\basicswsswaR[#1]#2{\SWSSWAR{}{#2}{#100}}%

\newcommand{\swsswaR}{\@ifnextchar[{\basicswsswaR}{\basicswsswaR[143]}}%

\newcommand{\SESSEAR}[3]{%
\Y=#3%
\divide \Y by 2%
\Z=#3%
\multiply \Z by 3%
\divide \Z by 4%
\begin{picture}(0,0)%
\put(-\Y,\Z){\line(2,-3){#3}}%
\put(\Y,-\Z){\sessehead}%
\truex{200}\truez{600}%
\put(\value{x},\value{x}){\makebox(0,\value{z})[l]{${#1}$}}%
\put(-\value{x},-\value{z}){\makebox(0,\value{z})[r]{${#2}$}}%
\end{picture}}%

\def\basicsessear[#1]{\SESSEAR{}{}{#100}}%

\newcommand{\sessear}{\@ifnextchar[{\basicsessear}{\basicsessear[143]}}%

\def\basicSessear[#1]#2{\SESSEAR{#2}{}{#100}}%

\newcommand{\Sessear}{\@ifnextchar[{\basicSessear}{\basicSessear[143]}}%

\def\basicsesseaR[#1]#2{\SESSEAR{}{#2}{#100}}%

\newcommand{\sesseaR}{\@ifnextchar[{\basicsesseaR}{\basicsesseaR[143]}}%

\newcommand{\NWNNWAR}[3]{%
\Y=#3%
\divide \Y by 2%
\Z=#3%
\multiply \Z by 3%
\divide \Z by 4%
\begin{picture}(0,0)%
\put(\Y,-\Z){\line(-2,3){#3}}%
\put(-\Y,\Z){\nwnnwhead}%
\truex{200}\truez{600}%
\put(\value{x},\value{x}){\makebox(0,\value{z})[l]{${#1}$}}%
\put(-\value{x},-\value{z}){\makebox(0,\value{z})[r]{${#2}$}}%
\end{picture}}%

\def\basicnwnnwar[#1]{\NWNNWAR{}{}{#100}}%

\newcommand{\nwnnwar}{\@ifnextchar[{\basicnwnnwar}{\basicnwnnwar[143]}}%

\def\basicNwnnwar[#1]#2{\NWNNWAR{#2}{}{#100}}%

\newcommand{\Nwnnwar}{\@ifnextchar[{\basicNwnnwar}{\basicNwnnwar[143]}}%

\def\basicnwnnwaR[#1]#2{\NWNNWAR{}{#2}{#100}}%

\newcommand{\nwnnwaR}{\@ifnextchar[{\basicnwnnwaR}{\basicnwnnwaR[143]}}%

\newcommand{\Necurve}[2]%
{\begin{picture}(0,0)%
\truex{1300}\truey{2000}\truez{200}%
\put(0,\value{x}){\oval(#200,\value{y})[t]}%
\put(0,\value{x}){\makebox(0,0){\begin{picture}(#200,0)%
\put(#200,0){\line(0,-1){\value{z}}}%
\put(#200,-\value{z}){\shead}%
\put(0,0){\line(0,-1){\value{z}}}\end{picture}}}%
\truex{2500}%
\put(0,\value{x}){\makebox(0,0)[b]{${#1}$}}%
\end{picture}}%

\def\basicnecurvar[#1]{\Necurve{}{#1}}

\newcommand{\necurvar}{\@ifnextchar[{\basicnecurvar}{\basicnecurvar[160]}}%

\def\basicNecurvar[#1]#2{\Necurve{#2}{#1}}%

\newcommand{\Necurvar}{\@ifnextchar[{\basicNecurvar}{\basicNecurvar[160]}}%

\newcommand{\Nwcurve}[2]%
{\begin{picture}(0,0)%
\truex{1300}\truey{2000}\truez{200}%
\put(0,\value{x}){\oval(#200,\value{y})[t]}%
\put(0,\value{x}){\makebox(0,0){\begin{picture}(#200,0)%
\put(#200,0){\line(0,-1){\value{z}}}%
\put(0,0){\line(0,-1){\value{z}}}%
\put(0,-\value{z}){\shead}%
\end{picture}}}%
\truex{2500}%
\put(0,\value{x}){\makebox(0,0)[b]{${#1}$}}%
\end{picture}}%

\def\basicnwcurvar[#1]{\Nwcurve{}{#1}}

\newcommand{\nwcurvar}{\@ifnextchar[{\basicnwcurvar}{\basicnwcurvar[160]}}%

\def\basicNwcurvar[#1]#2{\Nwcurve{#2}{#1}}%

\newcommand{\Nwcurvar}{\@ifnextchar[{\basicNwcurvar}{\basicNwcurvar[160]}}%

\newcommand{\Securve}[2]%
{\begin{picture}(0,0)%
\truex{1300}\truey{2000}\truez{200}%
\put(0,-\value{x}){\oval(#200,\value{y})[b]}%
\put(0,-\value{x}){\makebox(0,0){\begin{picture}(#200,0)%
\put(#200,0){\line(0,1){\value{z}}}%
\put(0,0){\line(0,1){\value{z}}}%
\put(#200,\value{z}){\nhead}%
\end{picture}}}%
\truex{2500}%
\put(0,-\value{x}){\makebox(0,0)[t]{${#1}$}}%
\end{picture}}%

\def\basicsecurvar[#1]{\Securve{}{#1}}

\newcommand{\securvar}{\@ifnextchar[{\basicsecurvar}{\basicsecurvar[160]}}%

\def\basicSecurvar[#1]#2{\Securve{#2}{#1}}%

\newcommand{\Securvar}{\@ifnextchar[{\basicSecurvar}{\basicSecurvar[160]}}%

\newcommand{\Swcurve}[2]%
{\begin{picture}(0,0)%
\truex{1300}\truey{2000}\truez{200}%
\put(0,-\value{x}){\oval(#200,\value{y})[b]}%
\put(0,-\value{x}){\makebox(0,0){\begin{picture}(#200,0)%
\put(#200,0){\line(0,1){\value{z}}}%
\put(0,0){\line(0,1){\value{z}}}%
\put(0,\value{z}){\nhead}%
\end{picture}}}%
\truex{2500}%
\put(0,-\value{x}){\makebox(0,0)[t]{${#1}$}}%
\end{picture}}%

\def\basicswcurvar[#1]{\Swcurve{}{#1}}

\newcommand{\swcurvar}{\@ifnextchar[{\basicswcurvar}{\basicswcurvar[160]}}%

\def\basicSwcurvar[#1]#2{\Swcurve{#2}{#1}}%

\newcommand{\Swcurvar}{\@ifnextchar[{\basicSwcurvar}{\basicSwcurvar[160]}}%

\newcommand{\Escurve}[2]%
{\begin{picture}(0,0)%
\truex{1400}\truey{2000}\truez{200}%
\put(\value{x},0){\oval(\value{y},#200)[r]}%
\put(\value{x},0){\makebox(0,0){\begin{picture}(0,#200)%
\put(0,0){\line(-1,0){\value{z}}}%
\put(0,#200){\line(-1,0){\value{z}}}%
\put(-\value{z},0){\whead}%
\end{picture}}}%
\truex{2500}%
\put(\value{x},0){\makebox(0,0)[l]{${#1}$}}%
\end{picture}}%

\def\basicescurvar[#1]{\Escurve{}{#1}}

\newcommand{\escurvar}{\@ifnextchar[{\basicescurvar}{\basicescurvar[160]}}%

\def\basicEscurvar[#1]#2{\Escurve{#2}{#1}}%

\newcommand{\Escurvar}{\@ifnextchar[{\basicEscurvar}{\basicEscurvar[160]}}%

\newcommand{\Encurve}[2]%
{\begin{picture}(0,0)%
\truex{1400}\truey{2000}\truez{200}%
\put(\value{x},0){\oval(\value{y},#200)[r]}%
\put(\value{x},0){\makebox(0,0){\begin{picture}(0,#200)%
\put(0,0){\line(-1,0){\value{z}}}%
\put(0,#200){\line(-1,0){\value{z}}}%
\put(-\value{z},#200){\whead}%
\end{picture}}}%
\truex{2500}%
\put(\value{x},0){\makebox(0,0)[l]{${#1}$}}%
\end{picture}}%

\def\basicencurvar[#1]{\Encurve{}{#1}}

\newcommand{\encurvar}{\@ifnextchar[{\basicencurvar}{\basicencurvar[160]}}%

\def\basicEncurvar[#1]#2{\Encurve{#2}{#1}}%

\newcommand{\Encurvar}{\@ifnextchar[{\basicEncurvar}{\basicEncurvar[160]}}%

\newcommand{\Wscurve}[2]%
{\begin{picture}(0,0)%
\truex{1300}\truey{2000}\truez{200}%
\put(-\value{x},0){\oval(\value{y},#200)[l]}%
\put(-\value{x},0){\makebox(0,0){\begin{picture}(0,#200)%
\put(0,0){\line(1,0){\value{z}}}%
\put(0,#200){\line(1,0){\value{z}}}%
\put(\value{z},0){\ehead}%
\end{picture}}}%
\truex{2400}%
\put(-\value{x},0){\makebox(0,0)[r]{${#1}$}}%
\end{picture}}%

\def\basicwscurvar[#1]{\Wscurve{}{#1}}

\newcommand{\wscurvar}{\@ifnextchar[{\basicwscurvar}{\basicwscurvar[160]}}%

\def\basicWscurvar[#1]#2{\Wscurve{#2}{#1}}%

\newcommand{\Wscurvar}{\@ifnextchar[{\basicWscurvar}{\basicWscurvar[160]}}%

\newcommand{\Wncurve}[2]%
{\begin{picture}(0,0)%
\truex{1300}\truey{2000}\truez{200}%
\put(-\value{x},0){\oval(\value{y},#200)[l]}%
\put(-\value{x},0){\makebox(0,0){\begin{picture}(0,#200)%
\put(0,0){\line(1,0){\value{z}}}%
\put(\value{z},#200){\ehead}%
\put(0,#200){\line(1,0){\value{z}}}%
\end{picture}}}%
\truex{2400}%
\put(-\value{x},0){\makebox(0,0)[r]{${#1}$}}%
\end{picture}}%

\def\basicwncurvar[#1]{\Wncurve{}{#1}}

\newcommand{\wncurvar}{\@ifnextchar[{\basicwncurvar}{\basicwncurvar[160]}}%

\def\basicWncurvar[#1]#2{\Wncurve{#2}{#1}}%

\newcommand{\Wncurvar}{\@ifnextchar[{\basicWncurvar}{\basicWncurvar[160]}}%

\catcode`\@=12

\begin{document}

\maketitle

\begin{abstract}
The notion of a group $G$ acting on a group $X$ is well-known. Fixing $X$, the corresponding functor $\Act(-,X)$ is representable by the group $[X]$ of automorphisms of $X$. The notion of $G$-action on $X$ has been generalized to the context of a semi-abelian category, but in this general context, the functor $\Act(-,X)$ is generally not representable. We investigate the representability of the functor $\Act(-,X)$ for the semi-abelian category $\EStarOp$, dual of the category of pointed objects of a topos $\calE$. The representability holds in particular when $\calE$ is a Boolean topos or a topos of presheaves of sets, but does not hold in general, not even for a Grothendieck topos $\calE$.
\end{abstract} 

\section{Introduction}

An action of the group $(G,\cdot)$ on the group $(X,+)$ is a mapping
$$
G\times X \ar X,~~~(g,x)\mapsto gx
$$
satisfying the axioms
$$
1x=x,~~(gg')x=g(g'x),~~g(x+x')=gx+gx'.
$$
This is the same as giving a group homomorphism
$G\to\Aut(X)$, where $\Aut(X)$ is the group of automorphisms of $(X,+)$. In other words, fixing the group $X$, the functor $\Act(-,X)$ mapping a group $G$ on the set of $G$-actions on $X$, is representable by the group $\Aut(X)$ of automorphisms of $X$. Moreover every $G$-action on $X$ induces a right-split short exact sequence
$$
\One \ar X \ar X\rtimes G \adjar G \ar \One
$$
where $X\rtimes G$ is the so-called {\em semi-direct product} of $G$ and $X$ for the given action. (These classical results are recalled in \cite{BB}, Section 5.2.)

The notion of semi-abelian category has been introduced in \cite{SAC} to provide an elegant non-abelian context where to develop -- in particular -- homological algebra. The most celebrated example of a semi-abelian category is precisely that of groups. It is shown in \cite{SDP} that semi-direct products, and via these the notion of action, can be defined in an arbitrary semi-abelian category. In \cite{SES}, a full proof is given of the bijection, in a semi-abelian category, between $G$-actions on $X$ and isomorphism classes of right-split exact sequences
$$
\One \ar X \ar A \adjar G \ar \One
$$
In the present paper, it will be convenient to rely at once on this alternative description of the functor $\Act(-,X)$.

Given an elementary topos $\calE$, the dual  $\EStarOp$ of the category $\EStar$ of pointed objects in $\calE$ is semi-abelian (see \cite{BT}). During his talk at his 60th birthday meeting in Coimbra, on July 13, 2012, George Janelidze asked the question of the representability of actions in this semi-abelian category $\EStarOp$. He showed that $X$ itself represents the actions on $X$ when the topos $\calE$ is Boolean. Referring to a more general (unpublished) result of James Gray, he showed also the representability of actions when $\calE$ is the topos of arrows in $\Set$ and produced a counter-example where actions are not representable: namely, the so-called {\em Par\'e topos} described in \cite{PARE}.

Through this paper, $\calE$ will remain an elementary topos and to avoid any confusion, we shall always work in the category $\EStar$, not in its dual. Except otherwise specified, every object  $A$ and every morphism $f$ considered in this paper  always lie in $\EStar$, without any mention of the base point, which will always be written ${\star}\colon\One\to A$. The zero morphisms are written $0$ while coproducts, kernels and cokernels in $\calE_{\star}$ are written $+$, $\Ker$ and $\Coker$.  Fixing $X\in\EStar$, the problem is thus to decide if the functor $\Act(X,-)$ mapping an object $G$ on the set of isomorphism classes of left-split exact sequences
$$
\One \ar G \adjar A \ar X \ar \One
$$
is representable by some object $[X]\in\EStar$.

Consider the full subcategory $\calN/X$ of $\EStar/X$ whose objects are the normal epimorphisms. We prove that the existence of an initial object $[X] \epi X$ in $\calN/X$ is sufficient to force $[X]$ to represent the actions on $X$. This sufficient condition is in particular satisfied when $\calE$ is a Boolean topos or when $\calE$ is a topos of presheaves over the topos of sets. We also exhibit a necessary condition for the representability of actions in $\EStar$ and use it to exhibit a Grothendieck topos $\calE$ for which actions are generally not representable in $\EStar$.

Of course I thank George Janelidze for having brought this problem to my attention. I thank him also for having proposed an elegant simplification of my original proof of Theorem~\ref{THM}.

\section{Initial normal covers}

The key assumption in this note is:

\begin{definition}\label{1}
Let $\calE$ be a topos.
An object $X\in\EStar$ {\em admits an initial normal cover} when the full subcategory $\calN/X$ of $\EStar/X$, whose objects are the normal epimorphisms, admits an initial object.
\end{definition}

When it exists, let us write $\chi\colon[X]\epi X$ for the initial normal cover of $X$.

\begin{example}\label{2}
When $\calE$ is a presheaf topos over $\Set$, every object of $\EStar$ admits an initial normal cover.
\end{example}

\proof
Consider the topos $\calE$ of presheaves on a small category $\calC$. As far as possible, the proof is developed in an arbitrary topos, just to underline clearly the points where the presheaf assumption is crucial.

It is standard basic category theory that a category $\calB$ has an initial object if and only if the identity functor admits a limit. And for this limit to exist, it suffices to prove that $\calB$ is complete and admits a small full subcategory $\calS$ such that
$$
\forall B\in\calB~~\exists S\in\calS~~\exists f\colon S\ar B
$$
($\cal S$ is called an {\em initial} or a {\em final} subcategory, according to the authors). We shall apply these results to prove the existence of an initial object in $\calN/X$.

To prove the completeness of $\calN_X$, let us show that $\calN_X$ is stable under limits in $\EStar/X$. Of course the identity on $X$ is the terminal object of $\calN/X$.

For equalizers, consider two parallel morphisms $u$, $v$ in $\calN/X$ and their equalizer $w$ in $\EStar$. With the notation of the following diagram, where $K$ and $L$ are the kernels of the normal epimorphisms $f$ and $g$, we must prove that $h=fw$ is a normal epimorphism as well.
\begin{diagram}[80]
K\cap E ? ? K ? ? L ??
\smono ? ? \smono ? ? \smono ??
E ? \Emono w ? A ? \Ebiar uv ? B ??
? ? \eseaR{h=fw} ? \Seepi f ? \sepI g ??
? ? ? ? X ??
\end{diagram}
The kernel of $h$ is of course $K\cap E$.

Let us prove first that $h$ is an epimorphism, using the internal logic of the topos. Given $x\in X$ there exists $a\in A$ such that $f(a)=x$. We have
$$
gu(a)=f(a)=x,~~~gv(a)=f(a)=x.
$$
Thus $gu(a)=gv(a)$ which implies, by normality of $g$
$$
\bigl(u(a),v(a)\bigr)\in L\times L~~\Or~~u(a)=v(a).
$$
But $u(a)\in L$ forces $x=f(a)=gu(a)={\star}$ while $u(a)=v(a)$ means $a\in E$. So we have
$$
x={\star}=fw({\star})=h({\star})~~\Or~~a\in E
$$
and since ${\star}$ itself is an element of $E$, in both cases we have
$$
\exists a\in E~~h(a)=x
$$
proving that $h$ is an epimorphism.

To prove that $h$ is the cokernel of $K\cap E$, consider two elements $e,e'\in E$ such that $h(e)=h(e')$. By normality of $f$, $fw(e)=fw(e')$ implies
$$
\bigl(w(e),w'(e')\bigr)\in K\times K~~\Or~~w(e)=w(e').
$$
And since $w$ is a monomorphism, this yields at once
$$
(e,e')\in (K\cap E)\times (K\cap E)~~\Or~~e=e'.
$$
This concludes the proof that $h$ is the cokernel of $K\cap E$.

For arbitrary products in $\calN_X$, consider a family $f_i\colon A_i\epi X$ of normal epimorphisms, for some indexing set $I$. Compute the corresponding generalized pullback $P$ in $\EStar$, which is the product in $\EStar/X$; write $f$ for the common composite $f_ip_i=f=f_jp_j$. We must prove that $f$ is a normal epimorphism.
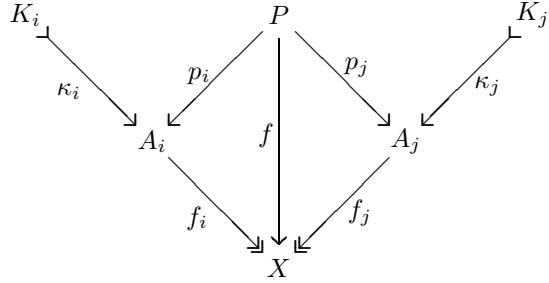
\begin{diagram}[60]
K_i ? ? ? ? P ? ? ? ? K_j ??
? \semonO{\kappa_i} ? ? \Swar{p_i} ? ? \Sear{p_j} ? ? \swmonO{\kappa_j} ??
? ? A_i ? ? \Sar[130]f ? ? A_j ??
? ? ? \seepI{f_i} ? ? \swepI{f_j} ??
? ? ? ? X ??
\end{diagram}

Fix an object $C\in\calC$ and an element $b\in X(C)$. Since we are in a topos of presheaves, each component $(f_i)_C$ is surjective, so that we can choose $a_i\in A_i(C)$ mapped on $b$ by $(f_i)_C$. Then $(a_i)_{i\in I}$ lies in the pullback $P(C)$ and is mapped on $b$ by $f_C$. So each $f_C$ is surjective and $f$ is an epimorphism.

Let us recall that in a topos of presheaves, the locales of subobjects are also co-locales, that is, finite joins distribute over arbitrary meets. Indeed this property holds pointwise in the category of sets. 

It remains to prove that $f$ is a normal epimorphism. Trivially, $\prod_{i\in I}K_i\subseteq P$, since every element of $K_i$ is mapped on $\star$ by $f_i$. Let us prove that $f$ is the cokernel of that inclusion, using the infinitary internal logic of the topos (or if you prefer, interpret what follows as occurring at an arbitrary level $C\in\calC$). Choose $(a_i)_{i\in I}$ and $(a'_i)_{i\in I}$ in $P$  identified by $f$. This means $f_i(a_i)=f_i(a'_i)$ for each index $i\in I$. By normality of $f_i$, 
$$
(a_i,a'_i)\in K_i\times K_i~~\Or~~a_i=a'_i.
$$
But $(a_i,a'_i)\in K_i\times K_i$ for some fixed index $i$ implies, for every other index $j$,
$$
f_j(a_j)=f_i(a_i)={\star}=f_i(a'_i)=f_j(a'_j)
$$
thus $(a_j,a'_j)\in K_j\times K_j$. So if one pair $(a_i,a'_i)$ is in $K_i\times K_i$, all pairs $(a_j,a'_j)$ are in $K_j\times K_j$. We have thus at each level $j\in I$
$$
\forall j\in I~~
\bigl((a_i)_{i\in I},(a'_i)_{i\in I}\bigr)
\in \left(\prod_{i\in I}K_i\right)
\times \left(\prod_{i\in I}K_i\right)
~~
\Or
~~
a_j=a'_j.
$$
Since finite joins distribute over arbitrary meets in a topos of presheaves, this is the same as
$$
\bigl((a_i)_{i\in I},(a'_i)_{i\in I}\bigr)
\in \left(\prod_{i\in I}K_i\right)
\times \left(\prod_{i\in I}K_i\right)
~~
\Or
~~
(a_j)_{j\in I}=(a'_j)_{j\in I}
.
$$
This proves that $f$ is a normal epimorphism.

It remains to check the ``solution set condition''. For this consider a normal epimorphism $g\colon A\epi X$. Since we are in a presheaf category, we have a normal epimorphism $g_C$ in $\Set$ at each level $C\in\calC$. Thus at each level $C$, $A(C)$ is constituted of
\begin{itemize}
\item an arbitrary number of elements mapped by $g_C$ on ${\star}\in X(C)$;
\item another piece, in bijection with $X(C)\setminus\{{\star}\}$ via $g_C$.
\end{itemize}
Write $B$ for the subpresheaf of $A$ generated by, for all levels $C$, the element ${\star}\in A(C)$ and the elements $a\in A(C)$ such that $g_C(a)\not={\star}$. There are thus at most $\card\coprod_{C\in\calC}X(C)$ such elements.
Trivially the composite
$$
B \mono A \Epi g X
$$
is still a normal epimorphism whose kernel is $K\cap B$, where $K$ is the kernel of $g$. 

The full subcategory $\calS\subseteq\calN/X$ we are looking for admits thus as objects those normal epimorphisms $B\epi X$ where $B$ is a presheaf generated by at most $\card\coprod_{C\in\calC}X(C)$ elements. Since the category $\calC$ is small, $\calS$ is small.
\qed

\begin{example}\label{3}
When $\calE$ is a Boolean topos , every object $X\in\EStar$ admits the identity on $X$ as initial normal cover.
\end{example}

\proof
Consider a normal epimorphism $f\colon A\epi X$ and its kernel $K$. As a subobject in the topos $\calE$, $K$ admits a complement $L'$ and putting $L=L'\cup\{{\star}\}$, we obtain $A=K+L$ in $\calE_{\star}$.
We have then the trivial exact sequence
$$
\One \ar K \Adjar{p_K}{i_K} K+L \Adjar{i_L}{p_L} L \ar \One
$$
where $K+L\cong A$ and $L=X$, $p_L\cong f$, since both are the cokernel of $K\subseteq A$.

We have thus already a commutative triangle 
\begin{diagram}[80]
\movevertexleft{X\cong L} ? ? \Ear[130]{i_L} ? ? \movevertexright{K+L \cong A} ??
? \seeql ? ? \swepI{p_L=f} ??
? ? X ??
\end{diagram}
It remains to prove that $i_L$ is the unique possible morphism making the diagram commutative. If $t$ is another morphism such that $ft=\id_X$, working in the internal logic of the topos, choose $l\in L$. Then $t(l)\in K$ or $t(l)\in L$. When $t(l)\in K$, we have $p_L\bigl(t(l)\bigr)={\star}$ and thus $l={\star}$ since $p_Lt=\id_L$; and since $l={\star}$,  $i_L(l)={\star}=t(l)$. On the other hand when $t(l)\in L$, we have $p_L\bigl(t(l)\bigr)=t(l)$ and again since $p_Lt=\id_L$, we have $i_L(l)=l=t(l)$. This proves that $t=i_L$ and thus the identity on $X$ is the initial normal cover of $X$.
\qed

Just to give a flavor of the possible form of an initial normal cover, let us recall an elementary example, also exhibited independently by James Gray.

\begin{example}\label{4}
The initial normal covers in the topos of arrows of sets.
\end{example}

\proof
Let $\calE$ be the topos of arrows in $\Set$. Fix an object $X=(f\colon X_1\ar X_0)$. Consider 
$$\setdisplayarrowlength{20}
[X]=\bigl(h\colon X_1\ar X_0+f^{-1}(\star)\bigr)
$$\setdisplayarrowlength{40}
where $+$ is the coproduct in the category of pointed sets. The morphism $h$ restricts as the identity from  $f^{-1}(\star)\subseteq X_1$ to $f^{-1}(\star)\subseteq X_0$ and as $f$ elsewhere. The morphism $\chi\colon [X]\ar X$ is the identity at the upper level while at the lower level, it is the identity on $X_0$ and maps $f^{-1}(\star)$ on $\star$. 

It is immediate that $\chi$ is a normal epimorphism, with kernel 
$$\setdisplayarrowlength{20}
K=\bigl(\{\star\}\ar f^{-1}(\star)\bigr).
\setdisplayarrowlength{40}$$
To prove that $\chi$ is an initial normal cover, consider a normal epimorphism $\varphi\colon A \epi X$, where $A=(g\colon A_1\ar A_0)$. By normality of $\chi$ and $\varphi$ at both levels, these morphisms induce at each level bijections
$$\setdisplayarrowlength{20}
[X]_i\setminus K_i
=
X_i\setminus\{\star\}
\cong
A_i\setminus (\Ker\varphi)_i,~~~i=1,0.
$$\setdisplayarrowlength{40}
This forces the definition on $[X]_i\setminus K_i$ of a possible factorization $\psi\colon [X] \ar A$ yielding $\varphi\psi=\chi$. This definition extends uniquely as a morphism $\psi\colon [X] \ar A$, namely, $\psi_1(\star)=\star$ and if $f(a)\not=\star$, $\psi_0(a)=g\psi_1(a)$.
\qed

In Example~\ref{4}, we have in particular the short exact sequence
$$
\One \ar K \Ar k [X] \Ar{\chi} X \ar \One.
$$
It should be noticed that in general, this is not a right-split exact sequence. Indeed choose $\star\not=a\in f^{-1}(\star)$. A retraction $[X]\ar K$ maps $a$ on $\star$ at the upper level, thus cannot map $a\in K_0$ on itself at the lower level.

The observant reader will have noticed that our proof of Example~\ref{4} remains valid when $\calE$ is replaced by an arbitrary Boolean topos. It remains an open problem to investigate an analogous internal generalization of Example~\ref{2}.

\section{A criterion of representability}

This section contains the main result of this paper. The proof given here is an elegant simplification of my original proof: a simplification suggested by George Janelidze. 

\begin{theorem}\label{THM}
Let $\calE$ be a topos and assume that $X\in\EStar$ admits the initial normal cover $\chi\colon [X]\epi X$. Then $[X]$ represents the actions on $X$ in the semi-abelian category $\EStarOp$.
\end{theorem}

\proof
Let us first construct the two pieces of the expected bijection. One of them is trivial: the mapping
$$
\alpha\colon\Act(X,G)\ar\EStar\bigl([X],G\bigr).
$$
Given a left-split exact sequence $(u,v,w)$ in $\EStar$, the assumption implies the existence of a unique morphism $t$ such that $wt=\chi$. 
\begin{diagram}[80]
? ? ? ? ? ? [X] ??
? ? ? ? ? \Swdotar{t~} ? \sepI{\chi} ??
\One ? \ear ? G ? \Eadjar vu ? A ? \eaR w ? X ? \ear ? \One ??
\end{diagram}
This yields at once a composite morphism $vt\colon[X]\ar G$ which we choose as $\alpha(u,v,w)$. The naturality of $\alpha$ in the variable $G$ is obvious.

Conversely, to construct 
$$
\beta\colon \EStar\bigl([X],G\bigr)\ar\Act(X,G)
$$
let us start with a morphism $s\colon [X] \ar G$. We write $k\colon K\mono [X]$ for the kernel of the initial normal cover $\chi\colon[X]\epi X$. We have in particular the morphisms
$$
K \Mono k [X] \Ar s G
$$
so that it makes sense to speak of the pushouts $[X]+_K[X]$ and $[X]+_KG$ over $K$.
The construction refers to the three upper rows in Diagram~1 and the morphisms between them. In this diagram, the morphisms $i_j$ are the canonical morphisms of the pushouts and $\nabla$ is the codiagonal. The part of the diagram so described is trivially commutative, with $\nabla$ and $(s,\id_G)$ retractions of the corresponding morphisms $i_2$. Moreover $(1)$ and $(1)+(2)$ are pushouts by construction, thus $(2)$ is a pushout as well.

\begin{floatingdiagram}[80]
\One ? \ear ? K ? \Ear k ? [X] ? \Ear{\chi} ? X ? \ear ? \One ??
 ? ? \Sar k ? (1) ? \saR{i_1} ? ? \seql ??
\One ? \ear ? [X] ? \Eadjar{\nabla}{i_2} ? [X]+_K[X] ? \Ear{(\chi,0)} ? X ? \ear ? \One ??
 ? ? \Sar s ? (2) ? \saR{\id_{[X]}+_Ks} ? ? \seql ??
\One ? \ear ? G ? \Eadjar{(s,\id_G)}{i_2} ? [X]+_KG ? \Ear{(\chi,0)} ? X ? \ear ? \One ??
 ? ? \seql ? ? \saR{(t,u)} ? ? \seql ??
\One ? \ear ? G ? \Eadjar vu ? A ? \Ear w ? X ? \ear ? \One ??
\end{floatingdiagram}

The top row of Diagram~1 is a short exact sequence, because $\chi$ is a normal monomorphism with kernel $k$. Therefore, the second and the third rows, obtained by pushouts from the first one, are short exact sequences as well; they are thus left-split short exact sequences. We define $\beta(s)$ to be the third row in Diagram~1.

Proving that $\alpha\beta(s)=s$ is easy. Indeed $i_1\colon [X]\ar [X]+_KG$ is such that $(\chi,0)i_1=\chi$. By initiality of $\chi$, $i_1$ is the unique such morphism and by definition of $\alpha$, $\alpha\beta(s)=(s,\id_G)i_1=s$.

Conversely start with the left-split short exact sequence $(u,v,w)$, bottom row in Diagram~1. Consider the unique morphism $t\colon [X]\ar A$ such that $wt=\chi$ and the composite $vt=\alpha(u,v,w)$. Choose $s=vt$ in Diagram~1: we must prove that $\beta(s)$, the third row of the diagram, is isomorphic to $(u,v,w)$, the bottom row of the diagram. For this we observe first that $tk=usk$. Indeed, since $wt=\chi$, we have $wtk=\chi k =0$, thus $tk$ factors as $tk=ut'$ through the kernel $u=\Ker w$. And then
$$
usk=uvtk=uvut'=ut'=tk.
$$
This equality implies the existence of the factorization $(t,u)$ through the pushout $[X]+_KG$. This factorization makes trivially the bottom part of Diagram~1 commutative, thus it is an isomorphism by the short five lemma (see \cite{PROTO}). This proves $\beta\alpha(u,v,w)\cong(u,v,w)$.
\qed 

Obviously, the proof of Theorem~\ref{THM}  does not use the full strength of the assumption on the existence of an initial normal cover. As in other questions related to the representability of actions (see \cite{REP}), the existence of a unique morphism $\bigl([X],\chi\bigr)\ar(A,w)$ in $\calN_X$ (notation of Definition~\ref{1}) is needed only for those normal epimorphisms $w\colon A\epi X$ whose kernel is a retract. 

One could decide to weaken accordingly the assumption of the Theorem. But when considering such a possible weakening, it is important to remain aware that the kernel of $\chi\colon [X] \epi X$ is generally not a retract, as the comment following Example~\ref{4} indicates.

The interested reader will also observe that the dual proof of Theorem~\ref{THM} carries over to the case of a pointed protomodular category with cokernels (see \cite{PROTO}). It yields then a criterion for the existence of a classifier for short right-split exact sequences, provided the existence of terminal normal extensions. But in the more general context of protomodular categories, one looses the connection with the notion of action. Moreover, as a matter of fact, in the well-known many ``algebraic-like'' examples of pointed protomodular categories or even semi-abelian categories, it is almost never the case that terminal normal extensions exist. Nevertheless \cite{REP} has exhibited many ``algebraic-like'' examples of semi-abelian categories where actions are representable.

\begin{corollary}
When $\calE$ is a topos of presheaves of sets, actions are representable in the semi-abelian category $\EStarOp$.
\end{corollary}

\proof
By Theorem~\ref{THM} and  Example~\ref{2}.
\qed

The following Corollary was first observed by George Janelidze, who gave a direct proof of it.

\begin{corollary}
When $\calE$ is a Boolean topos, actions on an object $X$ of the semi-abelian category $\EStarOp$ are representable by $X$ itself.
\end{corollary}

\proof
By Theorem~\ref{THM} and  Example~\ref{3}.
\qed

\section{Counter-examples}

It remains to give evidence that the representability of actions in the semi-abelian category $\EStarOp$, for a given topos $\calE$, is not a general fact. The first known (unpublished) counter-example has been produced by George Janelidze.

\begin{counterexample}
When $\calE$ is the {\em Par\'e topos} defined in \cite{PARE}, actions in $\EStarOp$ are generally not representable.
\end{counterexample}

\proof
Let $\calE$ be the topos of arrows in $\Set$. The full subcategory $\calP\subseteq\calE$ generated by those arrows with a finite codomain is stable in $\calE$ under finite limits, finite colimits, exponentiation and the $\Omega$-object. $\calP$ is thus a topos: let us call it the {\em Par\'e topos} (see \cite{PARE}). 

Consider a right-split short exact sequence in $\calE$.
$$
\One\ar G \adjar A \ar X \ar \One
$$
If both $G=(G_1\ar G_0)$ and $X=(X_1\ar X_0)$ are in $\calP$, then $A$ is in $\calP$ as well since $A_0$ is the coproduct as pointed sets of $G_0$ and $X_0$. This shows that for $X$ and $G$ in $\calP$, the $G$-actions on $X$ are the same in $\calE$ and $\calP$.

Fix $X\in\calP$. By Example~\ref{4}, actions on $X$ in $\calE$ are representable by some object $[X]\in\calE$. Suppose that actions on $X$ in $\calP$ are representable as well by some object $\langle X \rangle\in\calP$. These two representabilities, together with the observation concerning $G$-actions on $X$, with $G,X\in\calP$, force an isomorphism of functors
$$
\calP\bigl(\langle X \rangle,-\bigr)
\cong\calE\bigl([X]\bigr)
\colon \calP \ar \Set.
$$
Via this isomorphism, the identity on $\langle X \rangle$ corresponds to a morphism 
$$
\gamma\colon [X] \ar \langle X \rangle
$$
and the isomorphism of functors is then simply composition with $\gamma$.

Choose now $X=\bigl(\BbbN \ar \{\star\}\bigr)$ with $\star=0$ as base point at the top level. As proved in Example~\ref{4}, $[X]=\bigl(\BbbN\eql\BbbN\bigr)$. Let us show that the bottom component of the morphism 
$\gamma\colon [X] \ar \langle X \rangle$ is injective, which will contradict the fact that $\langle X \rangle$ lies in $\calP$.

Given $n\not=m$ in $\BbbN$, consider the morphism
$$
f\colon [X]=\bigl(\BbbN 
\setdisplayarrowlength{20}
\eql 
\setdisplayarrowlength{40}
\BbbN \bigr)
\ar 
\bigl(\BbbN
\setdisplayarrowlength{20}
\ar
\setdisplayarrowlength{40}
\{\star,n,m\}\bigr)
=Y.
$$
In $Y$, the elements $n,m\in Y_1$ are mapped respectively on $n,m\in Y_0$ and the rest on $\star$; analogously $f_1$ is the identity on $\BbbN$ while $f_0$ maps $n,m\in [X]_0$ respectively on $n,m\in Y_0$ and the rest on $\star$. Since $Y\in\calP$, we know by the isomorphism above that $f$ factors as $f=g\circ \gamma$. Since $f_0(n)\not=f_0(m)$, we have $\gamma_0(n)\not=\gamma_0(m)$. This proves the injectivity of $\gamma_0$ and thus the expected contradiction.
\qed

The {\em Par\'e topos} has been constructed to provide an example of a topos $\calE$ such that there does not exist any geometric morphism to a Boolean topos; it is thus definitely not a Grothendieck topos. But what about the representability of actions in $\EStarOp$, when $\calE$ is a Grothendieck topos?
 
The proof in the case of a topos of presheaves has shown that the bad behavior of infinite pullbacks of normal epimorphisms is the point which can possibly prevent the existence of initial normal covers in an arbitrary Grothendieck topos. This can be made more precise in terms of a necessary condition for the representability of actions.

\begin{proposition}
Let $\calE$ be a Grothendieck topos and $X\in\EStar$. For the actions on $X$ being representable in $\EStarOp$, it is necessary that given a family $f_i\colon A_i\epi X$ of normal epimorphisms in $\EStar$ whose kernels are retracts, the generalized pullback of these epimorphisms is still an epimorphism.
\end{proposition}

\proof
Suppose that the actions on $X$ are represented by an object $[X]$. The identity on $[X]$ is the initial object of the category of all arrows with domain $[X]$ in $\EStar$. By the representability assumption, this category is equivalent to that of left-split exact sequences with quotient object $X$. That category has thus an initial object as well, which we can write as
$$
\One \ar G \Adjar uv A \Epi w X \ar \One.
$$
Consider now a family $f_i\colon A_i \epi X$ of normal epimorphisms whose respective kernels are retracts. Choosing arbitrarily a retract for each index $i\in I$, by initiality of the split exact sequence above, we obtain commutative diagrams
\begin{diagram}[80]
\One ? \ear ? G ? \Eadjar vu ? A ? \Eepi w ? X ? \ear ? \One ??
? ? \sdotar ? ? \sdotaR{g_i} ? ? \seql ??
\One ? \ear ? K_i ? \eadjar ? A_i ? \Eepi{f_i} ? X ? \ear ? \One ??
\end{diagram}
Next consider the generalized pullback $P$ of the various $f_i$ and the corresponding composite $f\colon P \ar X$. Since $f_ig_i=w$ for all indices $i\in I$, there is a corresponding factorization $g\colon A\ar P$ through the pullback, making the following diagram commutative:
\begin{diagram}[80]
A ??
\Sdotar g ? \seaR{g_i} ? \Eseepi w ??
P ? \Ear{p_i} ? 
\cross{A_i}{\Securvar f} ? \Ear{f_i} ? X ??
\spacing(0,0,30)
\end{diagram}
Since $w$ is an epimorphism, so is $f$.
\qed

\begin{counterexample}
A Grothendieck topos for which actions are not always representable.
\end{counterexample}

\proof
Consider the poset $(\BbbN,\leq)$. For a non-empty set $I$ of indices, one has trivially
$$
a\vee\left(\bigwedge_{i\in I} b_i\right)
=
\bigwedge_{i\in I}(a\vee b_i)
$$
simply because $(\BbbN,\leq)$ is well-ordered, thus $\bigwedge_{i\in I}b_i$ is some $b_{i_0}$ of the family. As a consequence,
$$
\calL=(\BbbN\times\BbbN)\cup\{(\infty,\infty)\}
$$
provided with the reversed ordering, is a locale.
Let us draw this locale $\calL$ as in Diagram~2, with thus $(\infty,\infty)$ as bottom element. And let us freely speak of the elementary ``South-East'' and ``South-West'' restrictions in the case of a presheaf on $\calL$, as well as of the various elementary ``diamonds'' in the locale $\calL$.

\begin{floatingdiagram}[40]
 ? ? ? ? ? ? ? (0,0) ??
 ? ? ? ? ? ?
\nedotar[45]? ? \nwdotar[45] ??
 ? ? ? ? ? (0,1) ? ? ? ? (1,0)  ??
 ? ? ? ? 
\nedotar[45]? ? \nwdotar[45] ? ?
\nedotar[45]? ? \nwdotar[45] ? ?
??
 ? ? ? (0,2) ? ? ? ? (1,1) ? ? ? ? (2,0) ??
 ? ?  
\nedotar[45]? ? \nwdotar[45] ? ?
\nedotar[45]? ? \nwdotar[45] ? ?
\nedotar[45]? ? \nwdotar[45] ? ?
??
 ? (0,3) ? ? ? ? (1,2) ? ? ?  ? (2,1) ? ? ? ? (3,0) ??
\nedotar[45]? ? \nwdotar[45] ? ? 
\nedotar[45]? ? \nwdotar[45] ? ?
\nedotar[45]? ? \nwdotar[45] ? ?
\nedotar[45]? ? \nwdotar[45] ??
\end{floatingdiagram}

In the locale $\calL$, every non-trivial covering
$$
(i,j)=\bigvee_{s\in S}(i_s,j_s)
$$
has thus in $\BbbN\times\BbbN$ the form
$$
(i,j)
=
\left(\bigwedge_{s\in S}i_s,\bigwedge_{s\in S}j_s\right)
=
(i_{s_0},j_{s_1})
$$
since an infimum in $\BbbN$ is a smallest element. The covering contains in particular the two elements
$$
(i_{s_0},j_{s_0})=(i,j_{s_0}),~~~
(i_{s_1},j_{s_1})=(i_{s_1},j).
$$
So the covering is redundant since it contains already an {\em elementary covering}
$$
(i,j)=(i,j_{s_0})\vee(i_{s_1},j).
$$

Given an ``elementary covering'' as above
$$
(i,j)=(i,j')\vee(i',j)
$$
a compatible family 
$\bigl(a\in A(i,j'),b\in A(i',j)\bigr)$ 
along this covering in a presheaf $A$ over $\calL$ will simply be written $\langle a,b\rangle$, with thus $a$ in the ``West corner'' and $b$ in the ``East corner''.

Let us now define a pointed presheaf $T$ on $\calL$ by
$$
T(i,j)=\{\star,x\}.
$$
The South-West restrictions are the constant mappings on $\star$ while the South-East restrictions are the identities. This situation is pictured in Diagram~3. Given a diamond in $\calL$, the only possible compatible families are then
\begin{itemize}
\item
$\langle\star,\star\rangle$ which glues uniquely as $\star$;
\item
$\langle\star,x\rangle$ which glues uniquely as $x$.
\end{itemize}
This implies at once that $T$ is a sheaf on $\calL$, pointed by $\star$.

\begin{floatingdiagram}[40]
 ? ? ? ? ? ? ? \star x ??
 ? ? ? ? ? ?
\swdotar[45]? ? \sedotar[45] ??
 ? ? ? ? ? \star x ? ? ? ? \star x  ??
 ? ? ? ? 
\swdotar[45]? ? \sedotar[45] ? ?
\swdotar[45]? ? \sedotar[45] ? ?
??
 ? ? ? \star x ? ? ? ? \star x ? ? ? ? \star x ??
 ? ?  
\swdotar[45]? ? \sedotar[45] ? ?
\swdotar[45]? ? \sedotar[45] ? ?
\swdotar[45]? ? \sedotar[45] ? ?
??
 ? \star x ? ? ? ? \star x ? ? ?  ? \star x ? ? ? ? \star x ??
\swdotar[45]? ? \sedotar[45] ? ? 
\swdotar[45]? ? \sedotar[45] ? ?
\swdotar[45]? ? \sedotar[45] ? ?
\swdotar[45]? ? \sedotar[45] ??
\end{floatingdiagram}

We construct next a sequence of sheaves $A_n$, for each $n\in\BbbN$.
\begin{itemize}
\item
$A_n(i,j)=\{\star\}$ when $i< n$;
\item
$A_n(i,j)=\{\star,k, a\}$ when $i\geq n$.
\end{itemize}
All South-West restriction are such that
$$
\star\mapsto\star,~~~
k\mapsto k,~~~
a\mapsto k
$$
while all South-East restriction are given by
$$
\star\mapsto\star,~~~
k\mapsto k,~~~
a\mapsto a.
$$
This situation is pictured in Diagram~4 in the case $n=2$.
The composite of a South-West and a South-East restrictions, in whatever order, yields
$$
\star\mapsto\star,~~~
k\mapsto k,~~~
a\mapsto k
$$
thus $A_n$ is a presheaf. Observe further that the only possible compatible families along a diamond are
\begin{itemize}
\item 
$\langle\star,\star\rangle$ which glues uniquely as $\star$;
\item 
$\langle k,k\rangle$ which glues uniquely as $k$;
\item
$\langle k,a\rangle$ which glues uniquely as $a$.
\end{itemize}
This proves that $A_n$ is a sheaf, pointed by $\star$.  

\begin{floatingdiagram}[40]
 ? ? ? ? ? ? ? ? ? \star ??
 ? ? ? ? ? ? ? ?
\swdotar[45]? ? \sedotar[45] ??
 ? ? ? ? ? ? ? \star ? ? ? ? \star ??
 ? ? ? ? ? ? 
\swdotar[45]? ? \sedotar[45] ? ?
\swdotar[45]? ? \sedotar[45] ? ?
??
 ? ? ? ? ? \star ? ? ? ? \star ? ? ? ? \star ka ??
 ? ? ? ?  
\swdotar[45]? ? \sedotar[45] ? ?
\swdotar[45]? ? \sedotar[45] ? ?
\swdotar[45]? ? \sedotar[45] ? ?
??
 ? ? ? \star  ? ? ? ? \star  ? ? ?  ? \star ka ? ? ? ? \star ka ??
 ? ?
\swdotar[45]? ? \sedotar[45] ? ? 
\swdotar[45]? ? \sedotar[45] ? ?
\swdotar[45]? ? \sedotar[45] ? ?
\swdotar[45]? ? \sedotar[45] ??
 ? \star ? ? ? ? \star ? ? ? ? \star ka ? ? ?  ? \star ka ? ? ? ? \star ka ??
\swdotar[45]? ? \sedotar[45] ? ?
\swdotar[45]? ? \sedotar[45] ? ? 
\swdotar[45]? ? \sedotar[45] ? ?
\swdotar[45]? ? \sedotar[45] ? ?
\swdotar[45]? ? \sedotar[45] ??
\end{floatingdiagram}

We arrive at the family of normal epimorphisms. The morphism $f_n\colon A_n\to T$ is defined, at each level $(i,j)$, by
$$
\star\mapsto\star,~~~
k\mapsto\star,~~~
a\mapsto x.
$$
It is immediate that $f_n$ is a morphism of sheaves. 

The kernel $K_n$ of $f_n$ is given by
\begin{itemize}
\item
$K_n(i,j)=\{\star\}$ when $i< n$;
\item
$K_n(i,j)=\{\star, k\}$ when $i\geq n$.
\end{itemize}
Observe that $K_n$ becomes a retract of $A_n$ when defining the retraction by
$$
\star\mapsto\star,~~~
k\mapsto k,~~~
a\mapsto k.
$$

Let us prove that $f_n$ is an epimorphism of sheaves. Of course each element $\star\in T(i,j)$ is the image of the element $\star\in A_n(i,j)$. It remains to prove that given $x\in T(i,j)$, there exists a covering of $(i,j)$, thus as we have seen there exists an elementary covering
$$
(i,j)=(i,j')\vee(i',j),
$$
together with elements in $A_n(i,j')$, $A_n(i',j)$ mapped on the corresponding restrictions of $x$. Let us choose $i'$ and $j'$ strictly greater than $i$, $j$, $n$.  Then $x|_{(i,j')}=\star$ and $x|_{(i',j)}=x$. Of course $\star\in A_n(i,j')$ is mapped on $\star\in T(i,j')$ while $a\in A_n(i',j)$ exists and is mapped on $x\in T(i',j)$. This proves that $f_n$ is an epimorphism of sheaves.

Let us prove further that $f_n$ is the cokernel of the inclusion $K_n\subseteq T$. We must prove that when two elements $u,v\in A_n(i,j)$ are mapped on the same element of $T(i,j)$, there exists a covering of $(i,j)$ -- thus an elementary covering --  such that on one piece of the covering both restrictions of $u$ and $v$ are in $K_n$, while on the other piece of the covering the two restrictions of $u$ and $v$ are equal. Of course if $u$ and $v$ are mapped on $\star$, they are at once both in $K_n$. And if they are mapped on $x$, both $u$ and $v$ are equal to $a$ and thus are at once equal. So the required property holds at once ``set theoretically'' at each level $(i,j)$. This proves that the epimorphism $f_n$ is the cokernel of the inclusion $K_n\subseteq A_n$.

It remains to compute the infinite pullback $P$ of all the morphisms $f_n$ and consider the corresponding morphism $f\colon P \ar T$. This morphism $f$ cannot possibly be an epimorphism. Indeed choose $x\in T(0,0)$. If $f$ is an epimorphism, there exists a covering of $(0,0)$, thus an elementary covering, on which the restrictions of $x\in T(0,0)$ are images of elements in $P$. But such an elementary covering has necessarily the form
$$
(0,0) = (0,j)\vee(i,0).
$$
An element of $P(i,0)$ has the form $(a_n)_{n\in\BbbN}$ with $a_n\in A_n(i,0)$. The restriction of $x\in T(0,0)$ in $T(i,0)$ is still $x$. Thus $a_n\in A_n(i,0)$ is mapped by $f_n$ on $x\in T(i,0)$. By definition of $A_n$ and $f_n$, this forces $i\geq n$ and $a_n=a$. Thus the index $i$ should be greater than all the natural numbers, which is impossible. This concludes the proof.
\qed

\end{document}